


\catcode`\@=11
\font\tensmc=cmcsc10      
\def\smc{\tensmc}

\def\hcorrection#1{\advance\hoffset by #1 }
\def\vcorrection#1{\advance\voffset by #1 }
\def\wlog#1{}
\newif\iftitle@
\outer\def\title{\title@true\vglue 24\p@ plus 12\p@ minus 12\p@
   \bgroup\let\\=\cr\tabskip\centering
   \halign to \hsize\bgroup\tenbf\hfill\ignorespaces##\unskip\hfill\cr}
\def\endtitle{\cr\egroup\egroup\vglue 18\p@ plus 12\p@ minus 6\p@}
\outer\def\author{\iftitle@\vglue -18\p@ plus -12\p@ minus -6\p@\fi\vglue
    12\p@ plus 6\p@ minus 3\p@\bgroup\let\\=\cr\tabskip\centering
    \halign to \hsize\bgroup\smc\hfill\ignorespaces##\unskip\hfill\cr}
\def\endauthor{\cr\egroup\egroup\vglue 18\p@ plus 12\p@ minus 6\p@}
\outer\def\heading{\bigbreak\bgroup\let\\=\cr\tabskip\centering
    \halign to \hsize\bgroup\smc\hfill\ignorespaces##\unskip\hfill\cr}
\def\endheading{\cr\egroup\egroup\nobreak\medskip}

\outer\def\proclaim#1{\medbreak\noindent\smc\ignorespaces
    #1\unskip.\enspace\sl\ignorespaces}
\outer\def\endproclaim{\par\ifdim\lastskip<\medskipamount\removelastskip
  \penalty 55 \fi\medskip\rm}
\outer\def\demo#1{\par\ifdim\lastskip<\smallskipamount\removelastskip
    \smallskip\fi\noindent{\smc\ignorespaces#1\unskip:\enspace}\rm
      \ignorespaces}

\newcount\footmarkcount@
\footmarkcount@=1
\def\makefootnote@#1#2{\insert\footins{\interlinepenalty=100
  \splittopskip=\ht\strutbox \splitmaxdepth=\dp\strutbox
  \floatingpenalty=\@MM
  \leftskip=\z@\rightskip=\z@\spaceskip=\z@\xspaceskip=\z@
  \noindent{#1}\footstrut\rm\ignorespaces #2\strut}}
\def\footnote{\let\@sf=\empty\ifhmode\edef\@sf{\spacefactor
   =\the\spacefactor}\/\fi\futurelet\next\footnote@}
\def\footnote@{\ifx"\next\let\next\footnote@@\else
    \let\next\footnote@@@\fi\next}
\def\footnote@@"#1"#2{#1\@sf\relax\makefootnote@{#1}{#2}}
\def\footnote@@@#1{$^{\number\footmarkcount@}$\makefootnote@
   {$^{\number\footmarkcount@}$}{#1}\global\advance\footmarkcount@ by 1 }

\hyphenation{man-u-script man-u-scripts ap-pen-dix ap-pen-di-ces}
\hyphenation{data-base data-bases}
\ifx\amstexloaded@\relax\catcode`\@=13
  \endinput\else\let\amstexloaded@=\relax\fi
\newlinechar=`\^^J
\def\eat@#1{}
\def\Space@.{\futurelet\Space@\relax}
\Space@. %
\newhelp\athelp@
{Only certain combinations beginning with @ make sense to me.^^J
Perhaps you wanted \string\@\space for a printed @?^^J
I've ignored the character or group after @.}
\def\futureletnextat@{\futurelet\next\at@}
{\catcode`\@=\active
\lccode`\Z=`\@ \lowercase
{\gdef@{\expandafter\csname futureletnextatZ\endcsname}
\expandafter\gdef\csname atZ\endcsname
   {\ifcat\noexpand\next a\def\next{\csname atZZ\endcsname}\else
   \ifcat\noexpand\next0\def\next{\csname atZZ\endcsname}\else
    \def\next{\csname atZZZ\endcsname}\fi\fi\next}
\expandafter\gdef\csname atZZ\endcsname#1{\expandafter
   \ifx\csname #1Zat\endcsname\relax\def\next
     {\errhelp\expandafter=\csname athelpZ\endcsname
      \errmessage{Invalid use of \string@}}\else
       \def\next{\csname #1Zat\endcsname}\fi\next}
\expandafter\gdef\csname atZZZ\endcsname#1{\errhelp
    \expandafter=\csname athelpZ\endcsname
      \errmessage{Invalid use of \string@}}}}
\def\atdef@#1{\expandafter\def\csname #1@at\endcsname}
\newhelp\defahelp@{If you typed \string\define\space cs instead of
\string\define\string\cs\space^^J
I've substituted an inaccessible control sequence so that your^^J
definition will be completed without mixing me up too badly.^^J
If you typed \string\define{\string\cs} the inaccessible control sequence^^J
was defined to be \string\cs, and the rest of your^^J
definition appears as input.}
\newhelp\defbhelp@{I've ignored your definition, because it might^^J
conflict with other uses that are important to me.}
\def\define{\futurelet\next\define@}
\def\define@{\ifcat\noexpand\next\relax
  \def\next{\define@@}%
  \else\errhelp=\defahelp@
  \errmessage{\string\define\space must be followed by a control
     sequence}\def\next{\def\garbage@}\fi\next}
\def\undefined@{}
\def\preloaded@{}
\def\define@@#1{\ifx#1\relax\errhelp=\defbhelp@
   \errmessage{\string#1\space is already defined}\def\next{\def\garbage@}%
   \else\expandafter\ifx\csname\expandafter\eat@\string
	 #1@\endcsname\undefined@\errhelp=\defbhelp@
   \errmessage{\string#1\space can't be defined}\def\next{\def\garbage@}%
   \else\expandafter\ifx\csname\expandafter\eat@\string#1\endcsname\relax
     \def\next{\def#1}\else\errhelp=\defbhelp@
     \errmessage{\string#1\space is already defined}\def\next{\def\garbage@}%
      \fi\fi\fi\next}
\def\famzero{\fam\z@}

\def\lim{\mathop{\famzero lim}}

\def\max{\mathop{\famzero max}}
\def\min{\mathop{\famzero min}}

\def\sup{\mathop{\famzero sup}}

\def\textfont@#1#2{\def#1{\relax\ifmmode
    \errmessage{Use \string#1\space only in text}\else#2\fi}}
\textfont@\rm\tenrm
\textfont@\it\tenit
\textfont@\sl\tensl
\textfont@\bf\tenbf
\textfont@\smc\tensmc
\let\ic@=\/
\def\/{\unskip\ic@}
\def\textfonti{\the\textfont1 }
\def\t#1#2{{\edef\next{\the\font}\textfonti\accent"7F \next#1#2}}
\let\B=\=
\let\D=\.
\def~{\unskip\nobreak\ \ignorespaces}
{\catcode`\@=\active
\gdef\@{\char'100 }}
\atdef@-{\leavevmode\futurelet\next\athyph@}
\def\athyph@{\ifx\next-\let\next=\athyph@@
  \else\let\next=\athyph@@@\fi\next}
\def\athyph@@@{\hbox{-}}
\def\athyph@@#1{\futurelet\next\athyph@@@@}
\def\athyph@@@@{\if\next-\def\next##1{\hbox{---}}\else
    \def\next{\hbox{--}}\fi\next}
\def\.{.\spacefactor=\@m}
\atdef@.{\null.}
\atdef@,{\null,}
\atdef@;{\null;}
\atdef@:{\null:}
\atdef@?{\null?}
\atdef@!{\null!}
\def\srdr@{\thinspace}
\def\drsr@{\kern.02778em}
\def\sldl@{\kern.02778em}
\def\dlsl@{\thinspace}
\atdef@"{\unskip\futurelet\next\atqq@}
\def\atqq@{\ifx\next\Space@\def\next. {\atqq@@}\else
	 \def\next.{\atqq@@}\fi\next.}
\def\atqq@@{\futurelet\next\atqq@@@}
\def\atqq@@@{\ifx\next`\def\next`{\atqql@}\else\def\next'{\atqqr@}\fi\next}
\def\atqql@{\futurelet\next\atqql@@}
\def\atqql@@{\ifx\next`\def\next`{\sldl@``}\else\def\next{\dlsl@`}\fi\next}
\def\atqqr@{\futurelet\next\atqqr@@}
\def\atqqr@@{\ifx\next'\def\next'{\srdr@''}\else\def\next{\drsr@'}\fi\next}
\def\flushpar{\par\noindent}
\def\textfontii{\the\textfont2 }
\def\{{\relax\ifmmode\lbrace\else
    {\textfontii f}\spacefactor=\@m\fi}
\def\}{\relax\ifmmode\rbrace\else
    \let\@sf=\empty\ifhmode\edef\@sf{\spacefactor=\the\spacefactor}\fi
      {\textfontii g}\@sf\relax\fi}
\def\nonhmodeerr@#1{\errmessage
     {\string#1\space allowed only within text}}
\def\linebreak{\relax\ifhmode\unskip\break\else
    \nonhmodeerr@\linebreak\fi}
\def\allowlinebreak{\relax
   \ifhmode\allowbreak\else\nonhmodeerr@\allowlinebreak\fi}
\newskip\saveskip@
\def\nolinebreak{\relax\ifhmode\saveskip@=\lastskip\unskip
  \nobreak\ifdim\saveskip@>\z@\hskip\saveskip@\fi
   \else\nonhmodeerr@\nolinebreak\fi}
\def\newline{\relax\ifhmode\null\hfil\break
    \else\nonhmodeerr@\newline\fi}
\def\nonmathaerr@#1{\errmessage
     {\string#1\space is not allowed in display math mode}}
\def\nonmathberr@#1{\errmessage{\string#1\space is allowed only in math mode}}
\def\mathbreak{\relax\ifmmode\ifinner\break\else
   \nonmathaerr@\mathbreak\fi\else\nonmathberr@\mathbreak\fi}
\def\nomathbreak{\relax\ifmmode\ifinner\nobreak\else
    \nonmathaerr@\nomathbreak\fi\else\nonmathberr@\nomathbreak\fi}
\def\allowmathbreak{\relax\ifmmode\ifinner\allowbreak\else
     \nonmathaerr@\allowmathbreak\fi\else\nonmathberr@\allowmathbreak\fi}
\def\pagebreak{\relax\ifmmode
   \ifinner\errmessage{\string\pagebreak\space
     not allowed in non-display math mode}\else\postdisplaypenalty-\@M\fi
   \else\ifvmode\penalty-\@M\else\edef\spacefactor@
       {\spacefactor=\the\spacefactor}\vadjust{\penalty-\@M}\spacefactor@
	\relax\fi\fi}
\def\nopagebreak{\relax\ifmmode
     \ifinner\errmessage{\string\nopagebreak\space
    not allowed in non-display math mode}\else\postdisplaypenalty\@M\fi
    \else\ifvmode\nobreak\else\edef\spacefactor@
	{\spacefactor=\the\spacefactor}\vadjust{\penalty\@M}\spacefactor@
	 \relax\fi\fi}
\def\newpage{\relax\ifvmode\vfill\penalty-\@M\else\nonvmodeerr@\newpage\fi}
\def\nonvmodeerr@#1{\errmessage
    {\string#1\space is allowed only between paragraphs}}
\def\smallpagebreak{\relax\ifvmode\smallbreak
      \else\nonvmodeerr@\smallpagebreak\fi}
\def\medpagebreak{\relax\ifvmode\medbreak
       \else\nonvmodeerr@\medpagebreak\fi}
\def\bigpagebreak{\relax\ifvmode\bigbreak
      \else\nonvmodeerr@\bigpagebreak\fi}
\newdimen\captionwidth@
\captionwidth@=\hsize
\advance\captionwidth@ by -1.5in
\def\caption#1{}
\def\topspace#1{\gdef\thespace@{#1}\ifvmode\def\next
    {\futurelet\next\topspace@}\else\def\next{\nonvmodeerr@\topspace}\fi\next}
\def\topspace@{\ifx\next\Space@\def\next. {\futurelet\next\topspace@@}\else
     \def\next.{\futurelet\next\topspace@@}\fi\next.}
\def\topspace@@{\ifx\next\caption\let\next\topspace@@@\else
    \let\next\topspace@@@@\fi\next}
 \def\topspace@@@@{\topinsert\vbox to
       \thespace@{}\endinsert}
\def\topspace@@@\caption#1{\topinsert\vbox to
    \thespace@{}\nobreak
      \smallskip
    \setbox\z@=\hbox{\noindent\ignorespaces#1\unskip}%
   \ifdim\wd\z@>\captionwidth@
   \centerline{\vbox{\hsize=\captionwidth@\noindent\ignorespaces#1\unskip}}%
   \else\centerline{\box\z@}\fi\endinsert}
\def\midspace#1{\gdef\thespace@{#1}\ifvmode\def\next
    {\futurelet\next\midspace@}\else\def\next{\nonvmodeerr@\midspace}\fi\next}
\def\midspace@{\ifx\next\Space@\def\next. {\futurelet\next\midspace@@}\else
     \def\next.{\futurelet\next\midspace@@}\fi\next.}
\def\midspace@@{\ifx\next\caption\let\next\midspace@@@\else
    \let\next\midspace@@@@\fi\next}
 \def\midspace@@@@{\midinsert\vbox to
       \thespace@{}\endinsert}
\def\midspace@@@\caption#1{\midinsert\vbox to
    \thespace@{}\nobreak
      \smallskip
      \setbox\z@=\hbox{\noindent\ignorespaces#1\unskip}%
      \ifdim\wd\z@>\captionwidth@
    \centerline{\vbox{\hsize=\captionwidth@\noindent\ignorespaces#1\unskip}}%
    \else\centerline{\box\z@}\fi\endinsert}
\mathchardef\prime@="0230
\def\prime{{{}\prime@{}}}
\def\prim@s{\prime@\futurelet\next\pr@m@s}

\def\,{\relax\ifmmode\mskip\thinmuskip\else\thinspace\fi}
\def\!{\relax\ifmmode\mskip-\thinmuskip\else\negthinspace\fi}
\def\frac#1#2{{#1\over#2}}

\def\:{\nobreak\hskip.1111em{:}\hskip.3333em plus .0555em\relax}
\def\intic@{\mathchoice{\hskip5\p@}{\hskip4\p@}{\hskip4\p@}{\hskip4\p@}}
\def\negintic@
 {\mathchoice{\hskip-5\p@}{\hskip-4\p@}{\hskip-4\p@}{\hskip-4\p@}}
\def\intkern@{\mathchoice{\!\!\!}{\!\!}{\!\!}{\!\!}}
\def\intdots@{\mathchoice{\cdots}{{\cdotp}\mkern1.5mu
    {\cdotp}\mkern1.5mu{\cdotp}}{{\cdotp}\mkern1mu{\cdotp}\mkern1mu
      {\cdotp}}{{\cdotp}\mkern1mu{\cdotp}\mkern1mu{\cdotp}}}
\newcount\intno@
\def\iint{\intno@=\tw@\futurelet\next\ints@}
\def\iiint{\intno@=\thr@@\futurelet\next\ints@}
\def\iiiint{\intno@=4 \futurelet\next\ints@}
\def\idotsint{\intno@=\z@\futurelet\next\ints@}
\def\ints@{\findlimits@\ints@@}
\newif\iflimtoken@
\newif\iflimits@
\def\findlimits@{\limtoken@false\limits@false\ifx\next\limits
 \limtoken@true\limits@true\else\ifx\next\nolimits\limtoken@true\limits@false
    \fi\fi}
\def\multintlimits@{\intop\ifnum\intno@=\z@\intdots@
  \else\intkern@\fi
    \ifnum\intno@>\tw@\intop\intkern@\fi
     \ifnum\intno@>\thr@@\intop\intkern@\fi\intop}
\def\multint@{\int\ifnum\intno@=\z@\intdots@\else\intkern@\fi
   \ifnum\intno@>\tw@\int\intkern@\fi
    \ifnum\intno@>\thr@@\int\intkern@\fi\int}
\def\ints@@{\iflimtoken@\def\ints@@@{\iflimits@
   \negintic@\mathop{\intic@\multintlimits@}\limits\else
    \multint@\nolimits\fi\eat@}\else
     \def\ints@@@{\multint@\nolimits}\fi\ints@@@}
\def\Sb{_\bgroup\vspace@
	\baselineskip=\fontdimen10 \scriptfont\tw@
	\advance\baselineskip by \fontdimen12 \scriptfont\tw@
	\lineskip=\thr@@\fontdimen8 \scriptfont\thr@@
	\lineskiplimit=\thr@@\fontdimen8 \scriptfont\thr@@
	\Let@\vbox\bgroup\halign\bgroup \hfil$\scriptstyle
	    {##}$\hfil\cr}
\def\endSb{\crcr\egroup\egroup\egroup}
\def\Sp{^\bgroup\vspace@
	\baselineskip=\fontdimen10 \scriptfont\tw@
	\advance\baselineskip by \fontdimen12 \scriptfont\tw@
	\lineskip=\thr@@\fontdimen8 \scriptfont\thr@@
	\lineskiplimit=\thr@@\fontdimen8 \scriptfont\thr@@
	\Let@\vbox\bgroup\halign\bgroup \hfil$\scriptstyle
	    {##}$\hfil\cr}
\def\endSp{\crcr\egroup\egroup\egroup}
\def\Let@{\relax\iffalse{\fi\let\\=\cr\iffalse}\fi}
\def\vspace@{\def\vspace##1{\noalign{\vskip##1 }}}
\def\aligned{\,\vcenter\bgroup\vspace@\Let@\openup\jot\m@th\ialign
  \bgroup \strut\hfil$\displaystyle{##}$&$\displaystyle{{}##}$\hfil\crcr}
\def\endaligned{\crcr\egroup\egroup}
\def\matrix{\,\vcenter\bgroup\Let@\vspace@
    \normalbaselines
  \m@th\ialign\bgroup\hfil$##$\hfil&&\quad\hfil$##$\hfil\crcr
    \mathstrut\crcr\noalign{\kern-\baselineskip}}
\def\endmatrix{\crcr\mathstrut\crcr\noalign{\kern-\baselineskip}\egroup
		\egroup\,}
\newtoks\hashtoks@
\hashtoks@={#}
\def\format{\crcr\egroup\iffalse{\fi\ifnum`}=0 \fi\format@}
\def\format@#1\\{\def\preamble@{#1}%
  \def\c{\hfil$\the\hashtoks@$\hfil}%
  \def\r{\hfil$\the\hashtoks@$}%
  \def\l{$\the\hashtoks@$\hfil}%
  \setbox\z@=\hbox{\xdef\Preamble@{\preamble@}}\ifnum`{=0 \fi\iffalse}\fi
   \ialign\bgroup\span\Preamble@\crcr}

\def\cases{\left\{\,\vcenter\bgroup\vspace@
     \normalbaselines\openup\jot\m@th
       \Let@\ialign\bgroup$##$\hfil&\quad$##$\hfil\crcr
      \mathstrut\crcr\noalign{\kern-\baselineskip}}

\newif\iftagsleft@
\tagsleft@true
\def\TagsOnRight{\global\tagsleft@false}
\def\tag#1$${\iftagsleft@\leqno\else\eqno\fi
 \hbox{\def\pagebreak{\global\postdisplaypenalty-\@M}%
 \def\nopagebreak{\global\postdisplaypenalty\@M}\rm(#1\unskip)}%
  $$\postdisplaypenalty\z@\ignorespaces}
\interdisplaylinepenalty=\@M
\def\allowdisplaybreak@{\def\allowdisplaybreak{\noalign{\allowbreak}}}
\def\displaybreak@{\def\displaybreak{\noalign{\break}}}
\def\align#1\endalign{\def\tag{&}\vspace@\allowdisplaybreak@\displaybreak@
  \iftagsleft@\lalign@#1\endalign\else
   \ralign@#1\endalign\fi}
\def\ralign@#1\endalign{\displ@y\Let@\tabskip\centering\halign to\displaywidth
     {\hfil$\displaystyle{##}$\tabskip=\z@&$\displaystyle{{}##}$\hfil
       \tabskip=\centering&\llap{\hbox{(\rm##\unskip)}}\tabskip\z@\crcr
	     #1\crcr}}
\def\lalign@
 #1\endalign{\displ@y\Let@\tabskip\centering\halign to \displaywidth
   {\hfil$\displaystyle{##}$\tabskip=\z@&$\displaystyle{{}##}$\hfil
   \tabskip=\centering&\kern-\displaywidth
	\rlap{\hbox{(\rm##\unskip)}}\tabskip=\displaywidth\crcr
	       #1\crcr}}
\def\overrightarrow{\mathpalette\overrightarrow@}
\def\overrightarrow@#1#2{\vbox{\ialign{$##$\cr
    #1{-}\mkern-6mu\cleaders\hbox{$#1\mkern-2mu{-}\mkern-2mu$}\hfill
     \mkern-6mu{\to}\cr
     \noalign{\kern -1\p@\nointerlineskip}
     \hfil#1#2\hfil\cr}}}
\def\overleftarrow{\mathpalette\overleftarrow@}
\def\overleftarrow@#1#2{\vbox{\ialign{$##$\cr
     #1{\leftarrow}\mkern-6mu\cleaders\hbox{$#1\mkern-2mu{-}\mkern-2mu$}\hfill
      \mkern-6mu{-}\cr
     \noalign{\kern -1\p@\nointerlineskip}
     \hfil#1#2\hfil\cr}}}
\def\overleftrightarrow{\mathpalette\overleftrightarrow@}
\def\overleftrightarrow@#1#2{\vbox{\ialign{$##$\cr
     #1{\leftarrow}\mkern-6mu\cleaders\hbox{$#1\mkern-2mu{-}\mkern-2mu$}\hfill
       \mkern-6mu{\to}\cr
    \noalign{\kern -1\p@\nointerlineskip}
      \hfil#1#2\hfil\cr}}}
\def\underrightarrow{\mathpalette\underrightarrow@}
\def\underrightarrow@#1#2{\vtop{\ialign{$##$\cr
    \hfil#1#2\hfil\cr
     \noalign{\kern -1\p@\nointerlineskip}
    #1{-}\mkern-6mu\cleaders\hbox{$#1\mkern-2mu{-}\mkern-2mu$}\hfill
     \mkern-6mu{\to}\cr}}}
\def\underleftarrow{\mathpalette\underleftarrow@}
\def\underleftarrow@#1#2{\vtop{\ialign{$##$\cr
     \hfil#1#2\hfil\cr
     \noalign{\kern -1\p@\nointerlineskip}
     #1{\leftarrow}\mkern-6mu\cleaders\hbox{$#1\mkern-2mu{-}\mkern-2mu$}\hfill
      \mkern-6mu{-}\cr}}}
\def\underleftrightarrow{\mathpalette\underleftrightarrow@}
\def\underleftrightarrow@#1#2{\vtop{\ialign{$##$\cr
      \hfil#1#2\hfil\cr
    \noalign{\kern -1\p@\nointerlineskip}
     #1{\leftarrow}\mkern-6mu\cleaders\hbox{$#1\mkern-2mu{-}\mkern-2mu$}\hfill
       \mkern-6mu{\to}\cr}}}
\def\sqrt#1{\radical"270370 {#1}}
\def\dots{\relax\ifmmode\let\next=\ldots\else\let\next=\tdots@\fi\next}
\def\tdots@{\unskip\ \tdots@@}
\def\tdots@@{\futurelet\next\tdots@@@}
\def\tdots@@@{$\mathinner{\ldotp\ldotp\ldotp}\,
   \ifx\next,$\else
   \ifx\next.\,$\else
   \ifx\next;\,$\else
   \ifx\next:\,$\else
   \ifx\next?\,$\else
   \ifx\next!\,$\else
   $ \fi\fi\fi\fi\fi\fi}
\def\text{\relax\ifmmode\let\next=\text@\else\let\next=\text@@\fi\next}
\def\text@@#1{\hbox{#1}}
\def\text@#1{\mathchoice
 {\hbox{\everymath{\displaystyle}\def\textfonti{\the\textfont1 }%
    \def\textfontii{\the\textfont2 }\textdef@@ T#1}}
 {\hbox{\everymath{\textstyle}\def\textfonti{\the\textfont1 }%
    \def\textfontii{\the\textfont2 }\textdef@@ T#1}}
 {\hbox{\everymath{\scriptstyle}\def\textfonti{\the\scriptfont1 }%
   \def\textfontii{\the\scriptfont2 }\textdef@@ S\rm#1}}
 {\hbox{\everymath{\scriptscriptstyle}\def\textfonti{\the\scriptscriptfont1 }%
   \def\textfontii{\the\scriptscriptfont2 }\textdef@@ s\rm#1}}}
\def\textdef@@#1{\textdef@#1\rm \textdef@#1\bf
   \textdef@#1\sl \textdef@#1\it}

\def\textdef@#1#2{\def\next{\csname\expandafter\eat@\string#2fam\endcsname}%
\if S#1\edef#2{\the\scriptfont\next\relax}%
 \else\if s#1\edef#2{\the\scriptscriptfont\next\relax}%
 \else\edef#2{\the\textfont\next\relax}\fi\fi}
\scriptfont\itfam=\tenit \scriptscriptfont\itfam=\tenit
\scriptfont\slfam=\tensl \scriptscriptfont\slfam=\tensl
\mathcode`\0="0030
\mathcode`\1="0031
\mathcode`\2="0032
\mathcode`\3="0033
\mathcode`\4="0034
\mathcode`\5="0035
\mathcode`\6="0036
\mathcode`\7="0037
\mathcode`\8="0038
\mathcode`\9="0039
\def\Cal{\relax\ifmmode\let\next=\Cal@\else
     \def\next{\errmessage{Use \string\Cal\space only in math mode}}\fi\next}
\def\Cal@#1{{\fam2 #1}}
\def\bold{\relax\ifmmode\let\next=\bold@\else
   \def\next{\errmessage{Use \string\bold\space only in math
      mode}}\fi\next}\def\bold@#1{{\fam\bffam #1}}
\mathchardef\Gamma="0000
\mathchardef\Delta="0001
\mathchardef\Theta="0002
\mathchardef\Lambda="0003
\mathchardef\Xi="0004
\mathchardef\Pi="0005
\mathchardef\Sigma="0006
\mathchardef\Upsilon="0007
\mathchardef\Phi="0008
\mathchardef\Psi="0009
\mathchardef\Omega="000A
\mathchardef\varGamma="0100
\mathchardef\varDelta="0101
\mathchardef\varTheta="0102
\mathchardef\varLambda="0103
\mathchardef\varXi="0104
\mathchardef\varPi="0105
\mathchardef\varSigma="0106
\mathchardef\varUpsilon="0107
\mathchardef\varPhi="0108
\mathchardef\varPsi="0109
\mathchardef\varOmega="010A
\font\dummyft@=dummy
\fontdimen1 \dummyft@=\z@
\fontdimen2 \dummyft@=\z@
\fontdimen3 \dummyft@=\z@
\fontdimen4 \dummyft@=\z@
\fontdimen5 \dummyft@=\z@
\fontdimen6 \dummyft@=\z@
\fontdimen7 \dummyft@=\z@
\fontdimen8 \dummyft@=\z@
\fontdimen9 \dummyft@=\z@
\fontdimen10 \dummyft@=\z@
\fontdimen11 \dummyft@=\z@
\fontdimen12 \dummyft@=\z@
\fontdimen13 \dummyft@=\z@
\fontdimen14 \dummyft@=\z@
\fontdimen15 \dummyft@=\z@
\fontdimen16 \dummyft@=\z@
\fontdimen17 \dummyft@=\z@
\fontdimen18 \dummyft@=\z@
\fontdimen19 \dummyft@=\z@
\fontdimen20 \dummyft@=\z@
\fontdimen21 \dummyft@=\z@
\fontdimen22 \dummyft@=\z@
\def\fontlist@{\\{\tenrm}\\{\sevenrm}\\{\fiverm}\\{\teni}\\{\seveni}%
 \\{\fivei}\\{\tensy}\\{\sevensy}\\{\fivesy}\\{\tenex}\\{\tenbf}\\{\sevenbf}%
 \\{\fivebf}\\{\tensl}\\{\tenit}\\{\tensmc}}
\def\dodummy@{{\def\\##1{\global\let##1=\dummyft@}\fontlist@}}
\newif\ifsyntax@
\newcount\countxviii@
\def\newtoks@{\alloc@5\toks\toksdef\@cclvi}
\def\nopages@{\output={\setbox\z@=\box\@cclv \deadcycles=\z@}\newtoks@\output}
\def\syntax{\syntax@true\dodummy@\countxviii@=\count18
\loop \ifnum\countxviii@ > \z@ \textfont\countxviii@=\dummyft@
   \scriptfont\countxviii@=\dummyft@ \scriptscriptfont\countxviii@=\dummyft@
     \advance\countxviii@ by-\@ne\repeat
\dummyft@\tracinglostchars=\z@
  \nopages@\frenchspacing\hbadness=\@M}
\def\magstep#1{\ifcase#1 1000\or
 1200\or 1440\or 1728\or 2074\or 2488\or
 \errmessage{\string\magstep\space only works up to 5}\fi\relax}
{\lccode`\2=`\p \lccode`\3=`\t
 \lowercase{\gdef\tru@#123{#1truept}}}

\def\scaletype#1{\mag=#1\relax
 \hsize=\expandafter\tru@\the\hsize
 \vsize=\expandafter\tru@\the\vsize
 \dimen\footins=\expandafter\tru@\the\dimen\footins}

\def\scalefont#1#2\andcallit#3{\edef\font@{\the\font}#1\font#3=
  \fontname\font\space scaled #2\relax\font@}
\def\Mag@#1#2{\ifdim#1<1pt\multiply#1 #2\relax\divide#1 1000 \else
  \ifdim#1<10pt\divide#1 10 \multiply#1 #2\relax\divide#1 100\else
  \divide#1 100 \multiply#1 #2\relax\divide#1 10 \fi\fi}
\def\scalelinespacing#1{\Mag@\baselineskip{#1}\Mag@\lineskip{#1}%
  \Mag@\lineskiplimit{#1}}
\def\wlog#1{\immediate\write-1{#1}}
\catcode`\@=\active

\catcode`@=11
\def\binrel@#1{\setbox\z@\hbox{\thinmuskip0mu
\medmuskip\m@ne mu\thickmuskip\@ne mu$#1\m@th$}%
\setbox\@ne\hbox{\thinmuskip0mu\medmuskip\m@ne mu\thickmuskip
\@ne mu${}#1{}\m@th$}%
\setbox\tw@\hbox{\hskip\wd\@ne\hskip-\wd\z@}}
\def\overset#1\to#2{\binrel@{#2}\ifdim\wd\tw@<\z@
\mathbin{\mathop{\kern\z@#2}\limits^{#1}}\else\ifdim\wd\tw@>\z@
\mathrel{\mathop{\kern\z@#2}\limits^{#1}}\else
{\mathop{\kern\z@#2}\limits^{#1}}{}\fi\fi}
\def\underset#1\to#2{\binrel@{#2}\ifdim\wd\tw@<\z@
\mathbin{\mathop{\kern\z@#2}\limits_{#1}}\else\ifdim\wd\tw@>\z@
\mathrel{\mathop{\kern\z@#2}\limits_{#1}}\else
{\mathop{\kern\z@#2}\limits_{#1}}{}\fi\fi}
\def\circle#1{\leavevmode\setbox0=\hbox{h}\dimen@=\ht0
\advance\dimen@ by-1ex\rlap{\raise1.5\dimen@\hbox{\char'27}}#1}
\def\sqr#1#2{{\vcenter{\hrule height.#2pt
     \hbox{\vrule width.#2pt height#1pt \kern#1pt
       \vrule width.#2pt}
     \hrule height.#2pt}}}
\def\square{\mathchoice\sqr34\sqr34\sqr{2.1}3\sqr{1.5}3}
\def\force{\hbox{$\|\hskip-2pt\hbox{--}$\hskip2pt}}
 
\catcode`@=\active
 
\catcode`\@=11
\def\bold{\relaxnext@\ifmmode\let\next\bold@\else
 \def\next{\Err@{Use \string\bold\space only in math mode}}\fi\next}
\def\bold@#1{{\bold@@{#1}}}
\def\bold@@#1{\fam\bffam#1}
\def\hexnumber@#1{\ifnum#1<10 \number#1\else
 \ifnum#1=10 A\else\ifnum#1=11 B\else\ifnum#1=12 C\else
 \ifnum#1=13 D\else\ifnum#1=14 E\else\ifnum#1=15 F\fi\fi\fi\fi\fi\fi\fi}
\def\bffam@{\hexnumber@\bffam}
 
 
\font\tenmsx=msam10
\font\sevenmsx=msam7
\font\fivemsx=msam5
\font\tenmsy=msbm10
\font\sevenmsy=msbm7
\font\fivemsy=msbm7
 
\newfam\msxfam
\newfam\msyfam
\textfont\msxfam=\tenmsx
\scriptfont\msxfam=\sevenmsx
\scriptscriptfont\msxfam=\fivemsx
\textfont\msyfam=\tenmsy
\scriptfont\msyfam=\sevenmsy
\scriptscriptfont\msyfam=\fivemsy
\def\msx@{\hexnumber@\msxfam}
\def\msy@{\hexnumber@\msyfam}
\mathchardef\boxdot="2\msx@00
\mathchardef\boxplus="2\msx@01
\mathchardef\boxtimes="2\msx@02
\mathchardef\square="0\msx@03
\mathchardef\blacksquare="0\msx@04
\mathchardef\centerdot="2\msx@05
\mathchardef\lozenge="0\msx@06
\mathchardef\blacklozenge="0\msx@07
\mathchardef\circlearrowright="3\msx@08
\mathchardef\circlearrowleft="3\msx@09
\mathchardef\rightleftharpoons="3\msx@0A
\mathchardef\leftrightharpoons="3\msx@0B
\mathchardef\boxminus="2\msx@0C
\mathchardef\Vdash="3\msx@0D
\mathchardef\Vvdash="3\msx@0E
\mathchardef\vDash="3\msx@0F
\mathchardef\twoheadrightarrow="3\msx@10
\mathchardef\twoheadleftarrow="3\msx@11
\mathchardef\leftleftarrows="3\msx@12
\mathchardef\rightrightarrows="3\msx@13
\mathchardef\upuparrows="3\msx@14
\mathchardef\downdownarrows="3\msx@15
\mathchardef\upharpoonright="3\msx@16

\mathchardef\downharpoonright="3\msx@17
\mathchardef\upharpoonleft="3\msx@18
\mathchardef\downharpoonleft="3\msx@19
\mathchardef\rightarrowtail="3\msx@1A
\mathchardef\leftarrowtail="3\msx@1B
\mathchardef\leftrightarrows="3\msx@1C
\mathchardef\rightleftarrows="3\msx@1D
\mathchardef\Lsh="3\msx@1E
\mathchardef\Rsh="3\msx@1F
\mathchardef\rightsquigarrow="3\msx@20
\mathchardef\leftrightsquigarrow="3\msx@21
\mathchardef\looparrowleft="3\msx@22
\mathchardef\looparrowright="3\msx@23
\mathchardef\circeq="3\msx@24
\mathchardef\succsim="3\msx@25
\mathchardef\gtrsim="3\msx@26
\mathchardef\gtrapprox="3\msx@27
\mathchardef\multimap="3\msx@28
\mathchardef\therefore="3\msx@29
\mathchardef\because="3\msx@2A
\mathchardef\doteqdot="3\msx@2B

\mathchardef\triangleq="3\msx@2C
\mathchardef\precsim="3\msx@2D
\mathchardef\lesssim="3\msx@2E
\mathchardef\lessapprox="3\msx@2F
\mathchardef\eqslantless="3\msx@30
\mathchardef\eqslantgtr="3\msx@31
\mathchardef\curlyeqprec="3\msx@32
\mathchardef\curlyeqsucc="3\msx@33
\mathchardef\preccurlyeq="3\msx@34
\mathchardef\leqq="3\msx@35
\mathchardef\leqslant="3\msx@36
\mathchardef\lessgtr="3\msx@37
\mathchardef\backprime="0\msx@38
\mathchardef\risingdotseq="3\msx@3A
\mathchardef\fallingdotseq="3\msx@3B
\mathchardef\succcurlyeq="3\msx@3C
\mathchardef\geqq="3\msx@3D
\mathchardef\geqslant="3\msx@3E
\mathchardef\gtrless="3\msx@3F
\mathchardef\sqsubset="3\msx@40
\mathchardef\sqsupset="3\msx@41
\mathchardef\vartriangleright="3\msx@42
\mathchardef\vartriangleleft ="3\msx@43
\mathchardef\trianglerighteq="3\msx@44
\mathchardef\trianglelefteq="3\msx@45
\mathchardef\bigstar="0\msx@46
\mathchardef\between="3\msx@47
\mathchardef\blacktriangledown="0\msx@48
\mathchardef\blacktriangleright="3\msx@49
\mathchardef\blacktriangleleft="3\msx@4A
\mathchardef\vartriangle="3\msx@4D
\mathchardef\blacktriangle="0\msx@4E
\mathchardef\triangledown="0\msx@4F
\mathchardef\eqcirc="3\msx@50
\mathchardef\lesseqgtr="3\msx@51
\mathchardef\gtreqless="3\msx@52
\mathchardef\lesseqqgtr="3\msx@53
\mathchardef\gtreqqless="3\msx@54
\mathchardef\Rrightarrow="3\msx@56
\mathchardef\Lleftarrow="3\msx@57
\mathchardef\veebar="2\msx@59
\mathchardef\barwedge="2\msx@5A
\mathchardef\doublebarwedge="2\msx@5B
\mathchardef\angle="0\msx@5C
\mathchardef\measuredangle="0\msx@5D
\mathchardef\sphericalangle="0\msx@5E
\mathchardef\varpropto="3\msx@5F
\mathchardef\smallsmile="3\msx@60
\mathchardef\smallfrown="3\msx@61
\mathchardef\Subset="3\msx@62
\mathchardef\Supset="3\msx@63
\mathchardef\Cup="2\msx@64

\mathchardef\Cap="2\msx@65

\mathchardef\curlywedge="2\msx@66
\mathchardef\curlyvee="2\msx@67
\mathchardef\leftthreetimes="2\msx@68
\mathchardef\rightthreetimes="2\msx@69
\mathchardef\subseteqq="3\msx@6A
\mathchardef\supseteqq="3\msx@6B
\mathchardef\bumpeq="3\msx@6C
\mathchardef\Bumpeq="3\msx@6D
\mathchardef\lll="3\msx@6E

\mathchardef\ggg="3\msx@6F

\mathchardef\circledS="0\msx@73
\mathchardef\pitchfork="3\msx@74
\mathchardef\dotplus="2\msx@75
\mathchardef\backsim="3\msx@76
\mathchardef\backsimeq="3\msx@77
\mathchardef\complement="0\msx@7B
\mathchardef\intercal="2\msx@7C
\mathchardef\circledcirc="2\msx@7D
\mathchardef\circledast="2\msx@7E
\mathchardef\circleddash="2\msx@7F
\def\ulcorner{\delimiter"4\msx@70\msx@70 }
\def\urcorner{\delimiter"5\msx@71\msx@71 }
\def\llcorner{\delimiter"4\msx@78\msx@78 }
\def\lrcorner{\delimiter"5\msx@79\msx@79 }
\def\yen{{\mathhexbox@\msx@55 }}
\def\checkmark{{\mathhexbox@\msx@58 }}
\def\circledR{{\mathhexbox@\msx@72 }}
\def\maltese{{\mathhexbox@\msx@7A }}
\mathchardef\lvertneqq="3\msy@00
\mathchardef\gvertneqq="3\msy@01
\mathchardef\nleq="3\msy@02
\mathchardef\ngeq="3\msy@03
\mathchardef\nless="3\msy@04
\mathchardef\ngtr="3\msy@05
\mathchardef\nprec="3\msy@06
\mathchardef\nsucc="3\msy@07
\mathchardef\lneqq="3\msy@08
\mathchardef\gneqq="3\msy@09
\mathchardef\nleqslant="3\msy@0A
\mathchardef\ngeqslant="3\msy@0B
\mathchardef\lneq="3\msy@0C
\mathchardef\gneq="3\msy@0D
\mathchardef\npreceq="3\msy@0E
\mathchardef\nsucceq="3\msy@0F
\mathchardef\precnsim="3\msy@10
\mathchardef\succnsim="3\msy@11
\mathchardef\lnsim="3\msy@12
\mathchardef\gnsim="3\msy@13
\mathchardef\nleqq="3\msy@14
\mathchardef\ngeqq="3\msy@15
\mathchardef\precneqq="3\msy@16
\mathchardef\succneqq="3\msy@17
\mathchardef\precnapprox="3\msy@18
\mathchardef\succnapprox="3\msy@19
\mathchardef\lnapprox="3\msy@1A
\mathchardef\gnapprox="3\msy@1B
\mathchardef\nsim="3\msy@1C
\mathchardef\napprox="3\msy@1D
\mathchardef\varsubsetneq="3\msy@20
\mathchardef\varsupsetneq="3\msy@21
\mathchardef\nsubseteqq="3\msy@22
\mathchardef\nsupseteqq="3\msy@23
\mathchardef\subsetneqq="3\msy@24
\mathchardef\supsetneqq="3\msy@25
\mathchardef\varsubsetneqq="3\msy@26
\mathchardef\varsupsetneqq="3\msy@27
\mathchardef\subsetneq="3\msy@28
\mathchardef\supsetneq="3\msy@29
\mathchardef\nsubseteq="3\msy@2A
\mathchardef\nsupseteq="3\msy@2B
\mathchardef\nparallel="3\msy@2C
\mathchardef\nmid="3\msy@2D
\mathchardef\nshortmid="3\msy@2E
\mathchardef\nshortparallel="3\msy@2F
\mathchardef\nvdash="3\msy@30
\mathchardef\nVdash="3\msy@31
\mathchardef\nvDash="3\msy@32
\mathchardef\nVDash="3\msy@33
\mathchardef\ntrianglerighteq="3\msy@34
\mathchardef\ntrianglelefteq="3\msy@35
\mathchardef\ntriangleleft="3\msy@36
\mathchardef\ntriangleright="3\msy@37
\mathchardef\nleftarrow="3\msy@38
\mathchardef\nrightarrow="3\msy@39
\mathchardef\nLeftarrow="3\msy@3A
\mathchardef\nRightarrow="3\msy@3B
\mathchardef\nLeftrightarrow="3\msy@3C
\mathchardef\nleftrightarrow="3\msy@3D
\mathchardef\divideontimes="2\msy@3E
\mathchardef\varnothing="0\msy@3F
\mathchardef\nexists="0\msy@40
\mathchardef\mho="0\msy@66
\mathchardef\thorn="0\msy@67
\mathchardef\beth="0\msy@69
\mathchardef\gimel="0\msy@6A
\mathchardef\daleth="0\msy@6B
\mathchardef\lessdot="3\msy@6C
\mathchardef\gtrdot="3\msy@6D
\mathchardef\ltimes="2\msy@6E
\mathchardef\rtimes="2\msy@6F
\mathchardef\shortmid="3\msy@70
\mathchardef\shortparallel="3\msy@71
\mathchardef\smallsetminus="2\msy@72
\mathchardef\thicksim="3\msy@73
\mathchardef\thickapprox="3\msy@74
\mathchardef\approxeq="3\msy@75
\mathchardef\succapprox="3\msy@76
\mathchardef\precapprox="3\msy@77
\mathchardef\curvearrowleft="3\msy@78
\mathchardef\curvearrowright="3\msy@79
\mathchardef\digamma="0\msy@7A
\mathchardef\varkappa="0\msy@7B
\mathchardef\hslash="0\msy@7D
\mathchardef\hbar="0\msy@7E
\mathchardef\backepsilon="3\msy@7F
\def\Bbb{\relaxnext@\ifmmode\let\next\Bbb@\else
 \def\next{\Err@{Use \string\Bbb\space only in math mode}}\fi\next}
\def\Bbb@#1{{\Bbb@@{#1}}}
\def\Bbb@@#1{\noaccents@\fam\msyfam#1}
\catcode`\@=12

\font\tenmsy=msbm10
\font\sevenmsy=msbm7
\font\fivemsy=msbm5
\font\tenmsx=msam10
\font\sevenmsx=msam7
\font\fivemsx=msam5
\newfam\msyfam

\textfont\msyfam=\tenmsy
\scriptfont\msyfam=\sevenmsy
\scriptscriptfont\msyfam=\fivemsy
\newfam\msxfam

\textfont\msxfam=\tenmsx
\scriptfont\msxfam=\sevenmsx
\scriptscriptfont\msxfam=\fivemsx

%
 
\magnification 1200
\def\today{\ifcase\month\or January\or February\or March\or April\or
May\or June\or July\or August\or September\or October\or November\or
December \fi\space\number \day, \number\year}
\define\a{\alpha}

\define\egg{\Relbar\kern-.30em\Relbar\kern-.30em\Relbar}

\define\p1{P^1_{\delta, \lambda}}
\define\dom{\hbox{\rm dom}}
\define\pbf{\par\bigpagebreak\flushpar}
\define\g0{G^0_{\delta, \lambda}}
\define\cof{\hbox{\rm cof}}

\define\dell{\delta, \lambda}

\define\pj{^{p_j}}


\def\b{\beta}
\def\a{\alpha}
\def\l{\lambda}



\def\k{\kappa}
\pageno=1
\magnification 1200
\baselineskip=14pt 

\define\cU{\Cal U}



\def\b{\beta}
\def\a{\alpha}
\def\l{\lambda}


\def\k{\kappa}
\def\d{\delta}
\def\no{\noindent}
\def\g{\gamma}
\def\la{\langle}
\def\ra{\rangle}
\def\G{\Gamma}
\def\A{{\cal A}}
\def\B{{\cal B}}
\def\C{{\cal C}}
\def\K{{\cal K}}

\def\dom{{\hbox{\rm dom}}}
\def\max{{\hbox{\rm max}}}
\def\min{{\hbox{\rm min}}}

\def\pic{\underset i \in C \to{\prod} G_i}
\def\pic0{\underset i \in C_0 \to{\prod} G_i}

\def\sqr#1#2{{\vcenter{\vbox{\hrule height.#2pt
       \hbox{\vrule width.#2pt height#1pt \kern#1pt
         \vrule width.#2pt}
        \hrule height.#2pt}}}}
\def\square{\mathchoice\sqr34\sqr34\sqr{2.1}3\sqr{1.5}3}
\def\finpf{\no\hfill $\square$ }

\voffset=-.5in\vsize=7.5in

\vskip .3in
\centerline{``Menas' Result is Best Possible''}
\vskip .25in
\centerline{by}
\vskip .25in
\centerline{Arthur W. Apter*}
\centerline{Department of Mathematics}
\centerline{Baruch College of CUNY}
\centerline{New York, New York 10010}
\vskip .125in
\centerline{and}
\vskip .125in
\centerline{Saharon Shelah**}
\centerline{Department of Mathematics}
\centerline{The Hebrew University}
\centerline{Jerusalem, Israel}
\vskip .125in\centerline{and}\vskip .125in
\centerline{Department of Mathematics}
\centerline{Rutgers University}
\centerline{New Brunswick, New Jersey 08904}
\vskip .25in
\centerline{December 11, 1995}
\vskip .25in
\noindent Abstract: Generalizing some earlier techniques due to
the second author, we show that
Menas' theorem which states that the least cardinal
$\kappa$ which is a measurable limit of supercompact or
strongly compact cardinals is strongly compact but not
$2^\kappa$ supercompact is best possible.
Using these same techniques, we also extend and give a new proof
of a theorem of Woodin and extend and give a new proof of an
unpublished theorem due to the first author.
\hfil
\vskip .25in

\noindent *The research of the first author was partially supported
by PSC-CUNY Grant 662341 and a salary grant from Tel Aviv
University.  In addition, the first author wishes to thank the
Mathematics Departments of Hebrew University and Tel Aviv
University for the hospitality shown him during his
sabbatical in Israel.
\hfil\vskip .125in\noindent
**Publication 496.  The second author wishes to thank the
Basic Research Fund of the Israeli Academy of Sciences for
partially supporting this research.
\hfil\break
\eject

\S0 Introduction and Preliminaries
 
It is well known that if $\k$ is $2^\k$ supercompact, then
$\k$ is quite large in both size and consistency strength.
As an example of the former, if $\k$ is $2^\k$ supercompact,
then $\k$ has a normal measure concentrating on measurable
cardinals. The key to the proof of this fact and many other
similar ones is the existence of an elementary embedding
$j : V \to M$ with critical point $\k$ so that
$M^{2^\k} \subseteq M$. Thus, if $2^\k > \k^+$, one can
ask whether $\k$ must be large in size if $\k$ is merely
$\d$ supercompact for some $\k < \d < 2^\k$.
 
A natural question of the above venue to ask is whether a
cardinal $\k$ can be both the least measurable cardinal
and $\d$ supercompact for some $\k < \d < 2^\k$ if
$2^\k > \k^+$. Indeed, the first author posed this very
question to Woodin in the spring of 1983. In response, using
Radin forcing, Woodin (see [CW]) proved the following
\proclaim{Theorem} Suppose $V \models ``$ZFC + GCH +
$\k < \l$ are such that $\k$ is $\l^+$ supercompact and $\l$ is
regular''.
There is then a generic extension
$V[G]$ so that $V[G] \models ``$ZFC + $2^\k = \l$ +
$\k$ is $\d$ supercompact for all regular $\d < \l$ +
$\k$ is the least measurable cardinal''.
\endproclaim
 
The purpose of this paper is to extend the techniques of
[AS] and show that they can be used to
demonstrate that Menas' result of [Me] which says that the
least measurable cardinal $\k$ which is a limit of
supercompact or strongly compact cardinals is strongly
compact but not $2^\k$ supercompact is best possible.
Along the way, we generalize and strengthen Woodin's
result above, and we also produce a model in which, on a proper
class, the notions of measurability, $\d$ supercompactness,
and $\d$ strong compactness are all the same.
Specifically, we prove the following theorems.
\proclaim{Theorem 1} Suppose $V \models ``$ZFC + GCH +
$\k < \l$ are such that $\k$ is $< \l$ supercompact, $\l > \k^+$
is a regular cardinal which is either inaccessible or is the
successor of a cardinal of cofinality $> \k$,
and $h : \k \to \k$ is a function
so that for some elementary embedding $j : V \to M$
witnessing the $< \l$ supercompactness of $\k$,
$j(h)(\k) = \l$''. There is then a cardinal and cofinality
preserving generic extension $V[G] \models ``$ZFC +
For every inaccessible $\d < \k$ and every cardinal
$\g \in [\d,  h(\d))$, $2^\g = h(\d)$ +
For every cardinal $\g \in [\k, \l)$,
$2^\g = \l$ +
$\k$ is $<\l$ supercompact + $\k$
is the least measurable cardinal''.
\endproclaim
\proclaim{Theorem 2} Let $\l$ be a (class) function such
that for any infinite cardinal $\d$, $\l(\d) > \d^+$ is a
regular cardinal which is either inaccessible or is the
successor of a cardinal of cofinality $> \d$,
$\l(0) = 0$, and $\l(\d)$ is below
the least inaccessible $> \d$ if $\d$ is singular.
Suppose $V \models ``$ZFC + GCH +
$A$ is a proper class of cardinals so that for
each $\k \in A$, $h_\k : \k \to \k$ is a function and
$j_\k : V \to M$ is an elementary embedding witnessing the
$< \l(\k)$ supercompactness of $\k$ with
$j_\k(h_\k)(\k) = \l(\k) < \k^*$ for $\k^*$ the least
element of $A$ $> \k$''. There is then a cardinal and
cofinality preserving generic extension $V[G]
\models ``$ZFC + $2^\g = \l(\k)$ if $\k \in A$
and $\g \in [\k, \l(\k))$ is a cardinal +
There is a proper class of measurable cardinals +
$\forall \k                             [\k$ is
measurable iff $\k$ is $< \l(\k)$ strongly compact iff $\k$ is
$< \l(\k)$ supercompact] + No cardinal $\k$ is $\l(\k)$
strongly compact''.
\endproclaim
\proclaim{Theorem 3} Suppose $V \models ``$ZFC + GCH +
$\k$ is the least supercompact limit of supercompact
cardinals + $        \l > \k^+$ is a regular cardinal
which is either inaccessible or is the successor of
a cardinal of cofinality $> \k$ and
$h : \k \to \k$ is a function so that for some
elementary embedding $j : V \to M$ witnessing the
$< \l$ supercompactness of $\k$, $j(h)(\k) =
\l$''. There is then a
generic extension $V[G] \models ``$ZFC + For every
cardinal $\d < \k$ which is an inaccessible limit of
supercompact cardinals and every cardinal $\g \in
[\d, h(\d))$, $2^\g = h(\d)$ + For every cardinal
$\g \in [\k, \l)$, $2^\k = \l$ +
$\k$ is $<\l$ supercompact +
$\forall \d < \k[\d$ is strongly compact iff
$\d$ is supercompact] +
$\k$ is the least measurable limit of strongly compact or
supercompact cardinals''.
\endproclaim
 
Let us take this opportunity to make several remarks
concerning Theorems 1, 2, and 3. Note that we use a
weaker supercompactness hypothesis in the proof of Theorem 1 than 
Woodin does in the proof of his Theorem. Also, since
Woodin uses Radin forcing in the proof of his Theorem,
cofinalities are not preserved in his generic extension
(cardinals may or may not be preserved in Woodin's
Theorem, depending upon the proof used),
although they are in our Theorem 1 when the appropriate
forcing conditions are used. Further, in Theorem 2,
the model constructed will be so that
on the proper class composed of all cardinals possessing
some non-trivial degree of strong compactness or
supercompactness,
$\k$ is $\g$
strongly compact iff $\k$ is $\g$ supercompact, although
there won't be any (fully) strongly compact or (fully)
supercompact cardinals in this model. (This is the
generalized version of the theorem
that inspired the work of this paper and of [AS]. The
original theorem was
initially proven using an iteration of Woodin's version
of Radin forcing used to prove his above mentioned Theorem.)
Finally, Theorem 3 illustrates the flexibility of our
forcing as compared to Radin forcing. Since iterating a
Radin, Prikry, or Magidor [Ma1] forcing (for changing
the cofinality of $\k$ to some uncountable $\d < \k$)
above a strongly compact or supercompact cardinal $\k$
destroys the strong compactness or supercompactness of
$\k$, it is impossible to use any of these forcings in
the proof of Theorem 3. Our forcing for Theorem 3,
however, has been
designed so that, if $\k$ is a supercompact cardinal
which is Laver [L] indestructible,
then we can force above $\k$,
destroy measurability, yet preserve the supercompactness of
$\k$.
 
The structure of this paper is as follows. Section 0
contains our introductory comments and preliminary material
concerning notation, terminology, etc. Section 1 defines
and discusses the basic properties of the forcing notion
used in the iterations we employ to construct our models.
Section 2 gives a proof of Theorem 1. Section 3 contains
a proof of Theorems 2 and 3.
Section 4 concludes the paper by giving an alternate
forcing that can be used in the proofs of Theorems
1 and 2.
 
Before beginning the material of Section 1, we briefly mention some
preliminary information.
Essentially, our notation and terminology are standard, and when this is not the
case, this will be clearly noted.
For $\a < \b$ ordinals, $[\a, \b], [\a, \b), (\a, \b]$, and $(\a, \b) $ are as
 in
standard interval notation.
If $f$ is the characteristic function of a set
$x \subseteq \a$, then $x = \{\b < \a : f(\b) = 1\}$.
If $\a < \a'$, $f$ is a characteristic function having
domain $\a$, and $f'$ is a characteristic function having
domain $\a'$, we will when the context is clear abuse
notation somewhat and write $f \subseteq f'$, $f = f'$,
and $f \neq f'$ when we actually mean that the sets
defined by these functions satisfy these properties.
 
When forcing, $q \ge p$ will mean that $q$ is stronger than $p$.
  For $P$  a partial ordering, $\varphi$ a formula in the forcing language
with respect
 to $P$, and $ p \in P$, $ p \| \varphi$ will  mean
$p$ decides $\varphi$.
For $G$ $V$-generic over $P$, we will use both $V[G]$ and $V^{P}$ to indicate
 the universe obtained by forcing with $P$.
If $x \in V[G]$, then $\dot x$ will be a term in $V$ for $x$.
We may, from time to time, confuse terms with the sets they denote and write $x$
 when we actually mean $\dot x$, especially
 when $x$ is some variant of the generic set $G$, or $x$ is
in the ground model $V$.
 
If $\k$ is a cardinal and $P$ is
a partial ordering, $P$ is $\k$-closed if given a sequence
$\langle p_\a: \a < \k \rangle$ of elements of $P$ so that
$\beta < \gamma < \k$ implies $p_\beta \le p_\gamma$ (an increasing chain of
 length
 $\k$), then there is some $p \in P$ (an upper bound to this chain) so that
$p_\a \le p$ for all $\a < \k$.
$P$ is $<\k$-closed if $P$ is $\delta$-closed for all cardinals $\delta <
\k$.
$P$ is $(\k, \infty)$-distributive if for any sequence
$\la D_\a : \a < \k \ra$ of dense open subsets of $P$,
$D = \underset \a < \k \to{\cap} D_\a$ is a
dense open subset of $P$. $P$ is
$(< \k, \infty)$-distributive if $P$ is
$(\d, \infty)$-distributive for all cardinals $\d < \k$.
$P$ is $\k$-directed closed if for every cardinal $\delta < \k$ and every
 directed
 set $\langle p_\a : \a < \delta \rangle $ of elements of $P$
(where $\langle p_\alpha : \alpha < \delta \rangle$ is directed if
for every two distinct elements $p_\rho, p_\nu \in
\langle p_\alpha : \alpha < \delta \rangle$, $p_\rho$ and
$p_\nu$ have a common upper bound) there is an
upper bound $p \in P$. $P$ is $\k$-strategically closed if in the
two person game in which the players construct an increasing sequence
 $\langle p_\a: \a \le\k\rangle$, where player I plays odd stages and player
II plays even and limit stages, then player II has a strategy which ensures the
 game
 can always be continued.
$ P$ is $< \k$-strategically closed if $P$ is $\delta$-strategically
 closed for all cardinals $\delta < \k$.
$P$ is $\prec   \k$-strategically closed if in the two
person game in which the players construct an increasing
sequence $\langle p_\alpha : \alpha < \k \rangle$, where
player I plays odd stages and player II plays even and limit
stages, then player II has a strategy which ensures the game
can always be continued.
Note that trivially, if $P$ is $<\k$-closed, then $P$ is $<\k$-strategically
closed and $\prec \k  $-strategically closed. The converse of
both of these facts is false.
 
For $\k \le   \l$ regular cardinals, two partial orderings to  which we will refer
quite a bit are the standard partial orderings $\C(\l)$ for adding a Cohen
 subset to
$\l  $ using conditions having support $<\l$ and $\C(\k,\l)$ for adding
$\l  $ many Cohen subsets to $\k$ using conditions having
support $<\k$.
The basic properties and explicit definitions of these partial orderings
may be found in [J].
 
We mention that we are assuming complete familiarity with the notions
 of
measurability, strong compactness, and supercompactness.
Interested readers may consult [SRK], [Ka], or [KaM] for
further details.
We note only that all elementary embeddings witnessing the $\lambda$
supercompactness  of $\k$ are presumed to come from some
fine, $\k$-complete, normal
ultrafilter $\cU$ over $P_\k (\l) = \{ x \subseteq \l: | x| < \k \}$,
and all elementary embeddings witnessing the $< \l$
supercompactness of $\k$ for $\l$ a limit cardinal are
presumed to be generated by the appropriate system
$\la {\cal U}_\d : \d < \l \ra$ of ultrafilters over
$P_\k(\d)$ for $\d \in [\k, \l)$ a cardinal.
Also, where appropriate, all ultrapowers
will be confused with their transitive isomorphs.
 
Finally, we remark that a good deal of the notions and techniques
used in this paper are quite similar to those used in
[AS]. Since we desire this paper to be as comprehensible as
possible, regardless if readers have read [AS], many of
the arguments of [AS] will be repeated here in the
appropriately modified form.

\S 1 The Forcing Conditions
 
In this section, we describe and prove the basic properties of the forcing
 conditions
 we shall use in  our later iteration.
Let $\delta < \l, \, \l  >  \d^+    $ be regular cardinals in our ground model
 $V$, with $\d$ inaccessible and $\l$ either inaccessible or
the successor of a cardinal of cofinality $> \d$. We
assume throughout this section
that GCH holds for all cardinals $\k \ge \d$ in $V$, and
we define three notions of forcing.
Our first notion of forcing $P^0_{\delta, \l}$ is just the standard notion of
 forcing
 for adding a non-reflecting stationary set of ordinals of cofinality
$\delta $ to $\l$.
Specifically, $P^0_{\delta,\l} = \{ p$ : For some
$\a < \l$, $p : \a \to \{0,1\}$ is a characteristic
function of $S_p$, a subset of $\a$ not stationary at its
supremum nor having any initial segment which is stationary
at its supremum, so that $\b \in S_p$ implies
$\b > \d$ and cof$(\b) = \d \}$,
ordered by $q \ge p$ iff $q \supseteq p$ and $S_p = S_q \cap
 \sup (S_p)$, i.e., $S_q$ is an end extension of $S_p$. It is well-known
that for $G$ $V$-generic over $P^0_{\delta, \l}$ (see
[Bu] or   [KiM]), in $V[G]$,
since GCH holds in $V$ for all cardinals
$\k \ge \d$, a non-reflecting stationary
set $S=S[G]=\cup\{S_p:p\in G        \} \subseteq \l$ of ordinals of cofinality $\delta $ has been introduced,
 the
bounded subsets of $\l$ are the same as those in $V$,
and cardinals, cofinalities, and GCH at cardinals
$\k \ge \d$ have been preserved.
It is also virtually immediate that $P^0_{\d, \l}$
is $\d$-directed closed.
 
Work now in $V_1 = V^{P^0_{\delta, \l}}$, letting $\dot S$ be a term always
 forced
 to denote the above set $S$. $P^2_{\d, \l}[S]$ is the standard notion
 of forcing
 for introducing a club set $C$ which is disjoint to $S$ (and therefore
 makes $S$ non-stationary).
Specifically, $P^2_{\delta, \l} [S] = \{ p$ : For
some successor ordinal $\a < \l$,
$p : \a \to \{0,1\}$ is a characteristic function of
$C_p$, a club subset of $\a$, so that
$C_p \cap S = \emptyset \}$,
ordered by $ q \ge p $ iff $C_q$ is an end extension of $C_p$.
It is again well-known (see [MS]) that for $H$
$V_1$-generic over $P^2_{\delta, \l}[S]$, a club set
$C = C[H]         = \cup \{C_p : p \in H        \}
\subseteq \l$ which is disjoint to $S$ has been introduced, the bounded subsets
 of $\l$
are the same as those in $V_1$,
and cardinals, cofinalities, and GCH
for cardinals $\k \ge \d$ have been preserved.
 
More will be said about $P^0_{\d, \l}$ and
$P^2_{\d, \l}[S]$ in Lemmas 4, 6, and 7. In the meantime,
before defining in $V_1$ the partial ordering
$P^1_{\delta, \l}[S]$ which will be used to destroy
measurability, we first prove two preliminary lemmas.
 
\proclaim{Lemma 1} $\force_{P^0_{\delta, \l}} ``
\clubsuit({\dot S})$'', i.e., $V_1 \models ``$There is a
sequence $\langle x_\alpha : \alpha \in S \rangle$ so that
for each $\alpha \in S$, $ x_\alpha \subseteq \alpha$ is
cofinal in
$\alpha$, and for any $A \in
{[\l]}^{\l}$, $\{\alpha \in S : x_\alpha
\subseteq A      \}$ is stationary''.
\endproclaim
 
\demo{Proof of Lemma 1} Since GCH holds in $V$
for cardinals $\k \ge \d$ and $V$ and
$V_1$ contain the same bounded subsets of $\l$, we can let
$\langle y_\alpha : \alpha < \l \rangle \in V$ be a listing
of all elements $x \in {({[\l]}^\delta)}^V =
{({[\l]}^\delta)}^{V_1}$ so that each $x \in
{[\l]}^\delta$ appears on this list $\l$ times at
ordinals of cofinality $\d$, i.e., for any $x \in
{[\l]}^\d$, $\l = \sup\{\a < \l$ : cof$(\a) = \d$ and $y_\a
= x\}$. This then
allows us to define $\langle x_\alpha : \alpha \in
S \rangle$ by letting $x_\alpha$ be $y_\beta$ for the least
$\beta \in S - (\alpha + 1)$ so that $y_\beta \subseteq
\alpha$ and $y_\beta$ is unbounded in $\alpha$. By genericity,
each $x_\alpha$ is well-defined.
 
Let now $p \in P^0_{\delta, \l}$ be so that $p \force
``{\dot A} \in {[\l]}^{\l}$ and ${\dot K} \subseteq
\l$ is club''. We show that for some $r \ge p$ and some
$\zeta < \l$, $r \force ``\zeta \in {\dot K} \cap
{\dot S}$ and ${\dot x}_\zeta \subseteq {\dot A}$''. To do
this, we inductively define an increasing sequence $\langle
p_\alpha : \alpha < \delta \rangle$ of elements of
$P^0_{\delta, \l}$ and increasing sequences $\langle
\beta_\alpha : \alpha < \delta \rangle$ and $\langle
\gamma_\alpha : \alpha < \delta \rangle$ of ordinals
$< \l$ so that $\beta_0 \le \gamma_0 \le \beta_1 \le
\gamma_1 \le \cdots \le \beta_\alpha \le \gamma_\alpha
\le \cdots$ $(\alpha < \delta)$.
We begin by letting $p_0 = p$ and $\beta_0 = \gamma_0 =
0$. For $\eta = \alpha + 1 < \delta$ a successor, let
$p_\eta \ge p_\alpha$ and $\beta_\eta \le \gamma_\eta$,
$\beta_\eta \ge$max$(\beta_\alpha, \gamma_\alpha,
\sup($dom$(p_\alpha))) + 1$ be so that $p_\eta \force
``\beta_\eta \in {\dot A}$ and $\gamma_\eta \in
{\dot K}$''. For $\rho < \delta$ a limit, let $p_\rho =
\underset \alpha < \rho \to {\cup} p_\alpha$,
$\beta_\rho = \underset \alpha < \rho \to {\cup} \beta_\alpha$,
and $\gamma_\rho = \underset \alpha < \rho \to {\cup}
\gamma_\alpha$. Note that since $\rho < \delta$, $p_\rho$ is
well-defined, and since $\delta < \l$, $\beta_\rho,
\gamma_\rho < \l$. Also, by construction,
$\underset \alpha < \delta \to {\cup} \beta_\alpha =
\underset \alpha < \delta \to {\cup} \gamma_\alpha =
\underset \alpha < \delta \to {\cup} \sup($dom$(p_\alpha))
< \l$.
Call $\zeta$ this common sup. We thus have that $q =
\underset \alpha < \delta \to {\bigcup} p_\alpha \cup
\{\zeta\}$ is a well-defined condition so that $q \force
``\{\beta_\alpha : \alpha \in \delta - \{0\} \}
\subseteq \dot A$ and $\zeta \in {\dot K} \cap
{\dot S}$''.
 
To complete the proof of Lemma 1, we know that as $\langle
\beta_\alpha : \alpha \in \delta - \{0\} \rangle \in V$ and
as each $y \in \langle y_\alpha : \alpha < \l \rangle$ must
appear $\l$ times at ordinals of cofinality $\d$,
we can find some $\eta \in (\zeta,
\l)$ so that cof$(\eta) = \d$ and
$\langle \beta_\alpha : \alpha \in \delta -
\{0\} \rangle = y_\eta$. If we let $r \ge q$ be so that
$r \force ``{\dot S} \cap [\zeta, \eta] = \{\zeta, \eta\}$'',
then $r \force ``{\dot x}_\zeta = y_\eta =
\langle \beta_\alpha : \alpha \in \delta - \{0\}
\rangle$''. This proves Lemma 1.
 
\noindent \hfill $\square$ Lemma 1
 
We fix now in $V_1$ a $\clubsuit(S)$ sequence
$X = \la x_\a : \a \in S \ra$.
 
\proclaim{Lemma 2}
Let $S'$ be an initial segment of $S$ so that $S'$
is not stationary at its supremum nor has any initial
segment which is stationary at its supremum. There is then
a sequence $\la y_\a : \a \in S' \ra$ so that for every
$\a \in S'$, $y_\a \subseteq x_\a$, $x_\a - y_\a$ is
bounded in $\a$, and if $\a_1 \neq \a_2 \in S'$, then
$y_{\a_1} \cap y_{\a_2} = \emptyset$.
\endproclaim
 
\demo{Proof of Lemma 2}
We define by induction on $\a \le \a_0 = \sup S' + 1$
a function $h_\a$ so that $\dom(h_\a) = S' \cap \a$,
$h_\a(\b) < \b$, and $\la x_\b - h_\a(\b) :
\b \in S' \cap \a \ra$ is pairwise disjoint. The sequence
$\la x_\b - h_{\a_0}(\b) : \b \in S' \ra$ will be our
desired sequence.
 
If $\a = 0$, then we take $h_\a$ to be the empty function.
If $\a = \b + 1$ and $\b \not\in S'$, then we take
$h_\a = h_\b$. If $\a = \b + 1$ and $\b \in S'$, then
we notice that since each $x_\g \in X$ has order type
$\d$ and is cofinal in $\g$, for all $\g \in S' \cap
\b$, $x_\b \cap \g$ is bounded in $\g$. This allows us to
define a function $h_\a$ having domain $S' \cap \a$ by
$h_\a(\b) = 0$, and for $\g \in S' \cap \b$,
$h_\a(\g) = \min(\{\rho : \rho < \g$, $\rho \ge
h_\b(\g)$, and $x_\b \cap \g \subseteq \rho \})$. By the
next to last sentence and the induction hypothesis on
$h_\b$, $h_\a(\g) < \g$. And, if $\g_1 < \g_2 \in
S' \cap \a$, then if $\g_2 < \b$, $(x_{\g_1} -
h_\a(\g_1)) \cap (x_{\g_2} - h_\a(\g_2)) \subseteq
(x_{\g_1} - h_\b(\g_1)) \cap (x_{\g_2} -
h_\b(\g_2)) = \emptyset$ by the induction
hypothesis on $h_\b$.
If $\g_2 = \b$, then $(x_{\g_1} - h_\a(\g_1)) \cap
(x_{\g_2} - h_\a(\g_2)) = (x_{\g_1} - h_\a(\g_1)) \cap
x_{\g_2} = \emptyset$ by the definition of $h_\a(\g_1)$.
The sequence $\la x_\g - h_\a(\g) : \g \in S' \cap \a \ra$
is thus as desired.
 
If $\a$ is a limit ordinal, then as $S'$ is non-stationary
at its supremum nor has any initial segment which is
stationary at its supremum, we can let
$\la \b_\g : \g < \cof(\a) \ra$ be a strictly increasing,
continuous sequence having sup $\a$ so that for all
$\g < \cof(\a)$, $\b_\g \not\in S'$. Thus, if $\rho \in
S' \cap \a$, then $\{\b_\g : \b_\g < \rho \}$ is
bounded in $\rho$, meaning we can find some largest
$\g$ so that $\b_\g < \rho$. It is also the case that
$\rho < \b_{\g + 1}$. This allows us to define
$h_\a(\rho) = \max(\{h_{\b_{\g + 1}}(\rho),
\b_\g \})$ for the $\g$ just described. It is still the
case that $h_\a(\rho) < \rho$. And, if $\rho_1,
\rho_2 \in (\b_\g, \b_{\g + 1})$, then
$(x_{\rho_1} - h_\a(\rho_1)) \cap
 (x_{\rho_2} - h_\a(\rho_2)) \subseteq
(x_{\rho_1} - h_{\b_{\g + 1}}(\rho_1)) \cap
(x_{\rho_2} - h_{\b_{\g + 1}}(\rho_2))
= \emptyset$ by the definition of $h_{\b_{\g + 1}}$.
If $\rho_1 \in (\b_\g, \b_{\g + 1})$,
$\rho_2 \in (\b_\sigma, \b_{\sigma + 1})$ with
$\g < \sigma$, then
$(x_{\rho_1} - h_\a(\rho_1)) \cap
 (x_{\rho_2} - h_\a(\rho_2)) \subseteq
x_{\rho_1} \cap (x_{\rho_2} - \b_\sigma) \subseteq
\rho_1 - \b_\sigma \subseteq
\rho_1 - \b_{\g + 1} = \emptyset$.
Thus, the sequence $\la x_\rho - h_\a(\rho) :
\rho \in S' \cap \a \ra$ is again as desired.
This proves Lemma 2.
 
\no \hfill \finpf Lemma 2
 
At this point, we are in a position to define in $V_1$
the partial ordering $P^1_{\delta, \l}
[S] $ which will be used to destroy measurability.
$P^1_{\delta, \l} [S]$ is the set of all $5$-tuples
$\langle w, \a, \bar r, Z, \G \ra $ satisfying the following properties.
 
\item{1.} $w \subseteq \l$ is so that $|w| = \d$.
\item{2.} $\a < \d$.
\item{3.} $\bar r = \langle r_i : i \in w \rangle $ is a sequence of
functions from $\a $ to $\{ 0, 1 \}$, i.e., a sequence
of subsets of $\a$.
\item{4.} $Z$ is a function so that:
\item{a)} dom$(Z) \subseteq            \{x_\beta : \beta \in
S \}$ and range$(Z) \subseteq \{0, 1\}$.
\item{b)} If $z \in$ dom$(Z)$, then for some
$y \in {[w]}^\delta$, $y \subseteq z$ and $z - y$ is
bounded in the $\b$ so that $z = x_\b$.
\item{5.} $\G$ is a function so that:
\item{a)} dom$(\G) =$ dom$(Z)$.
\item{b)} If $z \in$ dom$(\G)$, then $\G(z)$ is a
closed, bounded subset of $\a$ such that if $\g$ is
inaccessible, $\g \in \G(z)$, and $\b$ is the
$\g$th element of $z$, then $\b \in w$, and for some
$\b' \in \b \cap w \cap z$, $r_{\b'}(\g) = Z(z)$.
 
\noindent Note that the definitions of $Z$ and $\G$ imply
$|$dom$(Z)| = |$dom$(\G)| \le \d$.
 
The ordering on $P^1_{\delta, \l} [S]$ is given by $\langle w^1, \alpha^1, \bar
 r^1,
Z^1, \G^1 \ra\le \langle w^2, \a^2, \bar r^2, Z^2,\G^2\ra$ iff the following
hold.
\item{1.} $w^1 \subseteq w^2$.
\item{2.} $\alpha^1 \le \a^2$.
\item{3.} If $i \in w^1$, then $r^1_i \subseteq r^2_i$
and $|\{i \in w^1 : r^2_i \vert (\a_2 - \a_1)
$ is not constantly $0 \}| < \d$.
\item{4.} $Z^1 \subseteq Z^2.$
\item{5.} dom$(\G^1) \subseteq$ dom$(\G^2)$.
\item{6.} If $z \in$ dom$(\G^1)$, then $\G^1(z)$ is
an initial segment of $\G^2(z)$ and
$|\{z \in {\hbox{\rm dom}}(\G^1) :
\G^1(z) \neq \G^2(z) \}| < \d$.
 
At this point, a few intuitive remarks are in order.
If $\d$ is measurable, then $\d$ must carry a normal measure.
The forcing $P^1_{\d, \l}[S]$ has specifically been designed
to destroy this fact. It has been designed, however, to
destroy the measurability of $\d$ ``as lightly as possible'',
making little damage, assuming $\d$ is $< \l$
supercompact. Specifically, if $\d$ is $< \l$
supercompact, then the non-reflecting stationary set $S$,
having been added to $\l$, does not kill the $< \l$
supercompactness of $\d$ by itself. The additional forcing
$P^1_{\d, \l}[S]$ is necessary to do the job and has been
designed so as not only to destroy the $< \l$
supercompactness of $\d$ but to destroy the measurability of
$\d$ as well. The forcing $P^1_{\d, \l}[S]$, however, has
been designed so that if necessary, we can resurrect the
$< \l$ supercompactness of $\d$ by forcing further with
$P^2_{\d, \l}[S]$.
 
\proclaim{Lemma 3} $V^{P^1_{\delta, \l}[S]}_1 \models ``\d$ is not
measurable''.
\endproclaim
 
\demo{Proof of Lemma 3}  Assume to the contrary that
$V^{P^1_{\delta, \l}[S]}_1 \models ``\d$ is measurable''.
Let $p \force `` \dot{\cal D}$ is a normal measure over $\d$''.
We show that $p$ can be extended to a condition $q$ so that
$q \force ``\dot {\cal D}$ is non-normal'', an
immediate contradiction.
 
We use a $\Delta$-system argument to establish this.
First, for $G_1$ $V_1$-generic over $P^1_{\delta, \l} [S]$ and $i < \l$,
let $r'_i = \cup \{ r^p_i: \exists p = \langle w^p, \a^p, \bar r^p, Z^p
, \G^p \ra \in$
$G_1          [r^p_i \in \bar r^p ] \}$.
An easy density argument shows $ \force_{P^1_{\delta, \l}[S]} ``\dot r'_i: \d \to \{
0, 1 \}$ is a function whose domain is all of $\d$".
Thus, we can let $ r^\ell_i = \{ \a < \d: r'_i (\a) = \ell \} $
for $\ell \in \{0,1\}$.
 
For each $i < \l  $, pick $p_i = \langle w^{p_i}, \a^{p_i}, \bar r^{p_i},
 Z^{p_i}, \G^{p_i} \ra
 \ge p$ so that $p_i \force ``\dot r^{\ell(i)}_i \not\in \dot {\cal D} "$
for some $\ell(i) \in \{0,1\}$. This is possible since
$\force_{P^1_{\delta, \lambda}[S]} ``$For each $i < \lambda  $,
${\dot r}^i_0 \cap {\dot r}^i_1 = \emptyset$ and
$\dot r^0_i \cup \dot r^1_i = \d       $''.
Without loss of generality, by extending $p_i$ if necessary,
since clause 4b) of the definition of the forcing implies
$|\dom(Z^{p_i})| \le \d$, we can assume that
$i \in w^{p_i}$ and $z \subseteq w^{p_i}$ for every
$z \in \dom(Z^{p_i})$.
Thus, since each $w^{p_i} \in [\l  ]^{< \d^+}$,
$\l > \d^+$, $\l$ is either inaccessible or is the successor of
a cardinal of
cofinality $> \d$, and GCH holds in $V_1$ for cardinals
$\k \ge \d$, we can find some
$A \in {[\l]}^\l$
so that $\{ w^{p_i} : i \in A \} $ forms a
$\Delta$-system, i.e., so that for $i \neq j \in A$, $w^{p_i} \cap
w\pj$ is some constant value $w$ which is an initial segment
of both. (Note we can assume that for $i \in A$, $w_i \cap i = w$,
and for some fixed $\ell(*) \in \{0,1\}$, for every $i \in A$,
$p_i \force ``\dot r^{\ell(*)}_i \not\in \dot{\cal D}$''.)
Also, by GCH in $V_1$ for cardinals $\k \ge \d$,
$|{[{\cal P}(w)]}^\d| = |{[\d^+]}^\d| = \d^+$.
Therefore, since $|\dom(Z^{p_i})| \le \d$ for each
$i < \l$ and $\l > \d^+$, we can assume in addition that
for all $i \in A$, dom$(Z^{p_i})
\cap {\cal P}(w) =$ dom$(\G^{p_i}) \cap {\cal P}(w)$ is
 some constant
 value $Z$.
Hence, since each $Z^{p_i}$ is a function from a set of
cardinality $\d$ into $\{0, 1\}$, each $\G^{p_i}$ is a function
from a set of cardinality $\d$ into ${[\d]}^{< \d}$
which has cardinality $\d$, and $\l > \d^+$, GCH in
$V_1$ for cardinals $\k \ge \d$ allows us to assume that for
$i \neq j \in A$, $Z^{p_i} \vert Z = Z^{p_j} \vert Z$ and
$\G^{p_i} \vert Z = \G^{p_j} \vert Z$.
Further, since each $\a{^{p_i}} < \d$, we can assume that $\a{^{p_i}}$ is some
 constant $\a^0$
 for $i \in A$.
Then, since any $\bar r{^{p_i}} = \langle r^{p_i}_j: j \in w{^{p_i}} \rangle$ for $i
 \in A$
is composed of a sequence of functions from $\a_0$ to $2$, $\a_0 < \d,$ and $|w| \le
 \d$, GCH in $V_1$ for cardinals $\k \ge \d$ again
allows us to assume that for $i \neq j \in A$, $
\bar r{^{p_i}} \vert w = \bar r\pj \vert w$.
And, since $i \in w{^{p_i}}$, we know that we can also assume (by thinning
 $A$ if necessary) that $B = \{ \sup (w{^{p_i}}): i \in A \}$ is so that $i < j
 \in A$
implies $i \le \sup (w{^{p_i}}) <$ min$(w^{p_j} - w) \le \sup (w\pj)$.
We know in addition by the choice of $X = \langle x_\beta : \beta
\in S \rangle$ that for some $\gamma \in S$, $x_\gamma
\subseteq A$. Let $x_\gamma = \{i_\beta : \beta < \delta \}$.
 
We are now in a position to define the condition $q$ referred to earlier.
We proceed by defining each of the five coordinates of $q$.
First, let $w^q	 = \underset \b < \delta \to{\cup} w^{p_{i_{\beta}}}$.
As $\d$ is regular, $\delta < \l$, and each
$w^{p_{i_\beta}} \in
[\l  ]^{<\d^+}$, $ w^q$
is well-defined and in ${[\lambda  ]}^{<\d^+}$.
Second, let $\a^q = \a^0$.
Third, let $\bar r^q
= \langle r^q_i : i \in w^q \rangle$ be defined by  $r^q_i =
 r^{p_{i_{\beta}}}_i$
 if $i \in w^{p_{i_{\beta}}}$.
The property of the $\Delta$-system that $i \neq j \in A$ implies $\bar
 r{^{p_i}} \vert w
= \bar r\pj \vert w$ tells us $\bar r^q$ is well-defined.
Finally, to define $Z^q$ and $\G^q$,
let $Z^q = \underset  \b < \delta \to{\cup} Z^{i_{\beta}} \cup
\{ \la\{ i_\beta : \beta < \delta \}, \ell(*) \ra \}$ and
$\G^q = \underset \b < \d \to{\cup} \G^{i_\b} \cup
\{ \la\{ i_\beta : \beta < \delta \}, \emptyset \ra \}$.
By the preceding paragraph and our construction,
$ \{ i_\beta: \beta < \delta \}$ generates a new set
which can be included in dom$(Z^q)$ and dom$(\G^q)$.
Therefore, since $Z^{p_i} \vert Z = Z^{p_j} \vert Z$
and $\G^{p_i} \vert Z = \G^{p_j} \vert Z$ for
$i \neq j \in A$,
$Z^q$ and $\G^q$ are well-defined.
 
We claim now that $q \ge p$ is so that $q \force ``
\dot {\cal D}$ is non-normal''. To see this, assume the
claim fails. Since $p \force ``\dot {\cal D}$ is a normal
ultrafilter over $\d$'' and by construction
$\forall \b < \d[q \ge p_{i_\b} \ge p]$,
$q \force ``{\dot r}^{\ell(*)}_{i_\b} \not \in \dot {\cal D}$''
for $\b < \d$. It must thus be the case that $q \force ``{\dot F}_0
= \{\g < \d : \g \in \underset \b < \g \to{\cup}
{\dot r}^{\ell(*)}_{i_\b}\} \not \in \dot {\cal D}$ and
${\dot F}_1 = \{\g < \d : \g$ is not inaccessible$\}
\not \in \dot {\cal D}$''. As $q \force ``{\dot K}^q
= \cup \{\G^s(\{i_\b : \b < \d \}) : \exists s =
\la w^s, \a^s, \bar r^s, Z^s, \G^s \ra \ge q
[s \in \dot G_1] \}$ is club in $\d$'', $q \force ``{\dot F}_2
= \{\g < \d : \g \not \in {\dot K}^q \} \not \in
\dot {\cal D}$'', so $q \force ``{\dot F} =
{\dot F}_0 \cup {\dot F}_1 \cup {\dot F}_2 \not \in
\dot {\cal D}$''.
 
We show that $q \force ``\dot F = \d$''. If $\g < \d$ is an
arbitrary inaccessible, then by the definition of ${\dot F}$,
it suffices to show that for some $s \ge q$ so that
$s \| ``\g \in \dot F$'', $s \force ``\g \in \dot F$''.
If $s \ge q$ is so that $s \force ``\g \in {\dot F}_2$'',
then we're done, so assume $s \force ``\g \not \in
{\dot F}_2$'', i.e., $s \force ``\g \in {\dot K}^q$''. But
then, by the definition of $\le$ and clause 5b) of the
definition of the forcing,
$s \force ``\g \in {\dot F}_0$''. Thus, $q \force ``\dot F =
\d$'', i.e., $q \force ``\dot F \in \dot {\cal D}$'',
meaning $q \ge p$ is so that $q \force ``\dot {\cal D}$ is
both a normal and non-normal ultrafilter over
$\d$''.
This proves Lemma 3.
 
\noindent
\hfill $\square $  Lemma 3
 
It is clear from the proof of Lemma 3 that since forcing with
$P^1_{\d, \l}[S]$ destroys the measurability of $\d$,
$P^1_{\d, \l}[S]$ can't be $\d$-directed closed. (Otherwise,
since $P^0_{\d, \l}$ is $\d$-directed closed, if $\d$
were supercompact and Laver [L] indestructible and
$P^1_{\d, \l}[S]$ were $\d$-directed closed, then the forcing
$P^0_{\d, \l} \ast P^1_{\d, \l}[S]$ would be $\d$-directed
closed and hence would preserve the supercompactness of
$\d$.) Note, however, that if $\g < \d$ and $\la p_i =
\la w^{p_i}, \a^{p_i}, \bar r^{p_i}, Z^{p_i}, \break \G^{p_i} \ra
: i < \g \ra$ is a directed sequence of conditions in $P^1_{\d, \l}
[S]$, it is possible
to define $\G^0 = \underset i < \g \to{\cup}
\G^{p_i}$, where $z \in {\hbox{\rm dom}}(\G^0)$ if
$z \in {\hbox{\rm dom}}(\G^{p_i})$ for some $i < \g$
and $\G^0(z)$ is the closure of
$\underset i < \g \to{\cup} \G^{p_i}(z)$.
($\G^{p_i}(z) = \emptyset$ if $z \not\in
\dom(\G^{p_i})$.) Then
$p = \la
\underset i < \g \to{\cup} w^{p_i},
\underset i < \g \to{\cup} \a^{p_i},
\underset i < \g \to{\cup} \bar r^{p_i},
\underset i < \g \to{\cup} Z^{p_i},
\G^0 \ra$,
where if $\bar r^{p_i} = \la r^{p_i}_j : j \in w^{p_i} \ra$,
then $r_j \in
\underset i < \g \to{\cup} \bar r^{p_i}$ if $j \in
\underset i < \g \to{\cup} w^{p_i}$ and $r_j =
\underset i < \g \to{\cup} r^{p_i}_j$ ($r^{p_i}_j =
\emptyset$ if
$j \not \in w^{p_i}$), is almost a condition.
The trouble occurs when for $z \in$ dom$(
\G^0)$,
$\G^0(z)$ contains a new element which is inaccessible. If, however,
we can guarantee that for any $z \in$ dom$(
\G^0)$, $\G^0(z)$ contains no new element which is inaccessible,
then $p$ as just defined is
a condition. Therefore,
we define a new partial ordering $\le^\g$ on
$P^1_{\d, \l}[S]$ by $p_1 =
\la w^{p_1}, \a^{p_1}, \bar r^{p_1}, Z^{p_1}, \G^{p_1} \ra
\le^\g p_2 =
\la w^{p_2}, \a^{p_2}, \bar r^{p_2}, Z^{p_2}, \G^{p_2} \ra
$
iff $p_1 = p_2$ or $p_1 < p_2$ and for
$z \in \dom(\G^{p_1})$, if $\G^{p_1}(z) \neq
\G^{p_2}(z)$, then $\g < \max(\G^{p_2}(z))$. If the sequence
$\la p_i =
\la w^{p_i}, \a^{p_i}, \bar r^{p_i}, Z^{p_i}, \G^{p_i} \ra
: i < \g \ra$
is a directed sequence of conditions in
$P^1_{\d, \l}[S]$ with respect to $\le^\g$, then since
$\sup(\underset i < \g \to{\cup} \G^{p_i}(z))$ must be an
ordinal $> \g$ of cofinality $\g$, the upper bound $p$
as defined earlier exists. Further, if
$p_1 =
\la w^{p_1}, \a^{p_1}, \bar r^{p_1}, Z^{p_1}, \G^{p_1} \ra
\le p_2 =
\la w^{p_2}, \a^{p_2}, \bar r^{p_2}, Z^{p_2}, \G^{p_2} \ra
$,
for any $z \in \dom(\G^{p_1})$ so that
$\G^{p_1}(z) \neq \G^{p_2}(z)$, we can define a function
$\G$ having domain $\G^{p_2}$ so that
$\G \vert (\dom(\G^{p_2}) - \dom(\G^{p_1})) =
\G^{p_2} \vert (\dom(\G^{p_2}) - \dom(\G^{p_1}))$
and such that for $z \in \dom(\G^{p_1})$,
$\G(z) = \G^{p_2}(z) \cup \{\eta_z\}$, where
$\eta_z < \d$ is the least cardinal above
$\max(\max(\G^{p_2}(z)), \g)$. If $q =
\la w^{p_2}, \a^{p_2}, \bar r^{p_2}, Z^{p_2}, \G \ra$,
then $q \in P^1_{\d, \l}[S]$ is a valid condition so that
$p_1 \le^\g q$ and $p_2 \le^\g q$. This easily implies that
$G$ is (appropriately) generic with respect to
$\la P^1_{\d, \l}[S], \le^\g \ra$, an ordering that is
$\g^+$-directed closed, iff $G$ is (appropriately) generic
with respect to
$\la P^1_{\d, \l}[S], \le \ra$, i.e., forcing with
$\la P^1_{\d, \l}[S], \le^\g \ra$ and
$\la P^1_{\d, \l}[S], \le \ra$ are equivalent. This key
observation will be critical in the proof of Theorem 3.
 
Recall we mentioned prior to the proof of Lemma 3 that $P^1_{\delta, \l}[S]$ is
designed so that a further forcing with $P^2_{\delta, \l} [S]$ will resurrect
 the $< \l$
supercompactness of $\d$, assuming the correct iteration has been done.
That this is so will be shown in the next section.
In the meantime, we give an idea of why this will happen by showing that the
 forcing
 $P^0_{\delta, \lambda} * (P^1_{\delta, \l} [\dot S] \times P^2_{\delta, \l}
[\dot S])$ is rather nice.
First, for $\k_0 \le \k_1$ regular cardinals, let
$\K(\k_0, \k_1) = \{\la w, \a, \bar r \ra :
w \subseteq \k_1$ is so that $|w| = \k_0$, $\a < \k_0$,
and $\bar r = \la r_i : i \in w \ra$ is a sequence of
functions from $\a$ to $\{0, 1\}\}$, ordered by
$\la w^1, \a^1, \bar r^1 \ra \le \la w^2, \a^2, \bar r^2 \ra$
iff $w^1 \subseteq w^2$, $\a^1 \le \a^2$, $r^1_i
\subseteq r^2_i$ if $i \in w^1$, and
$|\{i \in w^1 : r^2_i \vert (\a_2 - \a_1)$
is not constantly $0 \}| < \k_0$.
Given this definition, we now have the following lemma.
\proclaim{Lemma 4} $P^0_{\delta, \l} * (P^1_{\delta, \l}[\dot S] \times
 P^2_{\delta, \l}
 [\dot S])$ is equivalent to $\C(\l) * \dot \C(\d^+, \l)
* \dot \K(\d, \l)$.
\endproclaim
\demo{Proof of Lemma 4} Let $G$ be $V$-generic over $P^0_{\delta, \l}
* (P^1_{\dell} [\dot S] \times P^2_{\dell} [\dot S])$, with
$G        ^0_{\dell}$, $G        ^1_{\dell}$,
and $G        ^2_{\dell}$ the projections onto $P^0_{\dell}$, $P^1_{\dell} [S]$,
 and $P^2_{\dell} [S]$ respectively.
Each $G        ^i_{\dell}$ is appropriately generic.
So, since $P^1_{\dell} [S] \times P^2_{\dell} [S]$ is a product in
 $V[G^0_{\dell}]$,
 we can rewrite the forcing in $V[G^0_{\dell}]$ as $P^2_{\dell} [S]
\times P^1_{\dell} [S]$ and rewrite $V[G]$ as $V[G^0_{\dell}] [G^2_{\dell}]
[G^1_{\dell}]$.
 
It is well-known (see [MS]) that the forcing $P^0_{\dell} * P^2_{\dell} [\dot S]$ is
equivalent to $\C(\l)$.
That this is so can be seen from the fact that $P^0_{\dell} * P^2_{\dell} [\dot
 S]$ is non-trivial, has cardinality
$\l  $, and is such that $D=\{ \langle p, q
\rangle \in P^0_{\dell} * P^2_{\dell} [\dot S]$ : For some
$\alpha$, $\dom (p) = \dom (q)  = \a + 1$, $p \force `` \a \notin
\dot S"$, and $q \force `` \a \in \dot C "  \}$ is dense
 in $P^0_{\dell} * P^2_{\dell} [\dot S]$ and is $< \l$-closed.
This easily implies the desired equivalence.
Thus, $V$ and $V[G^0_{\dell}] [G^2_{\dell}]$ have the same cardinals and
cofinalities, and the proof of Lemma 4 will be complete once we show that in
$V[G^0_{\dell}] [G^2_{\dell}]$, $P^1_{\dell}  [S]$ is equivalent to
$\C(\d^+, \l) \ast \dot \K(\d, \l)$.
 
To this end, working in $V[G^0_{\d, \l}][G^2_{\d, \l}]$, let
$R = \{Z : Z$ is a function from $\{  x_\b : \b \in S \} $
into $\{0, 1\}$ so that $|{\hbox{\rm dom}}(Z)| \le \d \}$,
ordered by inclusion. Since $|\{  x_\b : \b \in S \} |
= \l$ and $R$ is $\d$-closed, it is clear $R$ is
equivalent to $\C(\d^+, \l)$. Further, the following
facts are easy to see.
\item{1.} If $p = \la w^p, \a^p, \bar r^p, Z^p, \G^p \ra \in
P^1_{\d, \l}[S]$, then $Z^p \in R$.
\item{2.} If $p_1 = \la w^{p_1}, \a^{p_1}, \bar r^{p_1}, Z^{p_1}, \G^{p_1} \ra$,
$p_2 = \la w^{p_2}, \a^{p_2}, \bar r^{p_2}, Z^{p_2}, \G^{p_2} \ra$ are so that
$p_1, p_2 \in P^1_{\d, \l}[S]$ and $p_1 \le p_2$, then
$Z^{p_1} \subseteq Z^{p_2}$.
\item{3.} If $p_1 = \la w^{p_1}, \a^{p_1}, \bar r^{p_1}, Z^{p_1}, \G^{p_1} \ra
\in P^1_{\d, \l}[S]$ is so that $Z^{p_1} \subseteq Z^{p_2}$ for
some $Z^{p_2} \in R$, then there exists $p_2 \in
P^1_{\d, \l}[S]$ with $p_1 \le p_2$, $p_2 =
\la w^{p_2}, \a^{p_2}, \bar r^{p_2}, Z^{p_2}, \G^{p_2} \ra$.
 
\noindent From these three facts, it then easily follows that
$H = \{Z \in R : \exists p \in G^1_{\d, \l}[Z = Z^p] \}$
is $V[G^0_{\d, \l}][G^2_{\d, \l}]$-generic over $R$.
This means we can rewrite $P^1_{\d, \l}[S]$ in
$V[G^0_{\d, \l}][G^2_{\d, \l}]$ as
$R \ast (\dot P^1_{\d, \l}[S] / R)$, which is isomorphic to
$\C(\d^+, \l) \ast (\dot P^1_{\d, \l}[S] / R)$. We will
thus be done if we can show in
$V[G^0_{\d, \l}][G^2_{\d, \l}][H]$ (which has the same
cardinals and cofinalities as $V$ and
$V[G^0_{\d, \l}][G^2_{\d, \l}]$) that
$P^1_{\d, \l}[S] / R$ is equivalent to $\K(\d, \l)$.
 
Working now in $V[G^0_{\dell}][G^2_{\dell}][H]$, we first note that
as $S \subseteq \l  $ is in
$V[G^0_{\d, \l}][G^2_{\d, \l}]$ a non-stationary set
all of whose initial
 segments are
 non-stationary,
by Lemma 2, for the sequence $\langle x_\beta: \beta \in S \rangle$,
 there
must be a sequence $\langle y_\beta: \beta \in S \rangle 
\in V[G^0_{\d, \l}][G^2_{\d, \l}] \subseteq
V[G^0_{\d, \l}][G^2_{\d, \l}][H]$ so that for every
$\beta \in S$, $y_\beta \subseteq x_\beta$, $x_\beta - y_\beta$ is bounded in
 $\b$,
 and if $\beta_1 \neq \beta_2 \in S$, then $y_{\beta_{1}} \cap y_{\beta_{2}}
= \emptyset$.
Given this fact, it is easy to observe that $P^1 = \{ \langle w, \a, \bar r,
\G \rangle \in P^1_{\dell} [S] / R:$ For every $\beta \in S$,
 either $y_\beta \subseteq w $ or $ y_\beta \cap w = \emptyset\}$
is dense in $P^1_{\dell} [S] / R$.
To show this, given $\langle w, \a, \bar r,   \G \rangle \in P^1_{\dell} [S] / R$,
$ \bar r = \langle r_i : i \in w \rangle$,
let $Y_w = \{ y \in \langle y_\beta: \beta \in S \rangle: y \cap w \neq
 \emptyset \}$.
As $|w| \le \d$ and $y_{\beta_{1}} \cap y_{\beta_{2}}  = \emptyset $ for $\b_1
 \neq \b_2 \in S$, $| Y_w | \le \d$.
Hence, as $|y| = \delta < \l$ for $y \in Y_w$,
$|w'|\le \d$ for $ w' = w \cup
( \cup Y_w)$.
This means $\langle w', \a, \bar r',    \G\rangle $ for $\bar r'= \langle
 r'_i: i \in w' \rangle$ defined
 by $r'_i = r_i$ if $ i \in w $ and $r'_i$ is the empty function
if $i \in w' - w$
is a well-defined condition extending $\langle w, \a, \bar r
, \G \rangle$.
Thus, $P^1$ is dense in $P^1_{\dell} [S] / R$, so to analyze the forcing properties
 of $P^1_{\dell} [S] / R$, it suffices to analyze the forcing properties of $P^1$.
 
For $\b \in S$, let $Q_\b = \{ \langle w, \a, \bar r,    \G \rangle \in P^1: w =
 y_\b \}$,
and let $Q^* = \{ \langle w, \a, \bar r,    \G \rangle  \in P^1: w \subseteq \l   -
\underset \b \in S \to{\cup} y_\beta \}$.
Let $Q$ be those elements of
$\underset \b \in S \to{\prod} Q_\b \times Q^*$ of support
$\d$ so that for $p = \la \la w^{p_i}, \a^{p_i},
\bar r^{p_i}, \G^{p_i} \ra_{i < \d},
\la w^p, \a^p, \bar r^p, \G^p \ra \ra \in Q$,
$\a^{p_i} = \a^{p_j} = \a^p$ for $i < j < \d$.
Let $\le_Q$ on $Q$ be defined by
$p = \la \la w^{p_i}, \a, \bar r^{p_i}, \G^{p_i}
\ra_{i < \d}, \la w^p, \a, \bar r^p, \G^p \ra \ra
\le_Q q = \la \la w^{q_i}, \b, \bar r^{q_i}, \G^{q_i}
\ra_{i < \d}, \la w^q, \b, \bar r^q, \break \G^q \ra \ra$
iff the following hold.
\item{1.} $\la w^p, \a, \bar r^p, \G^p \ra \le
\la w^q, \b, \bar r^q, \G^q \ra$.
\item{2.} $q$ can be written in the form
$\la \la w^{q_i}, \b, \bar r^{q_i}, \G^{q_i}
\ra_{i < \d}, \la u^{q_i}, \b, \bar s^{q_i},
\Delta^{q_i} \ra_{i < i_0 \le \d},
\la w^q, \b, \bar r^q, \break \G^q \ra \ra$ so that
$\forall i < \d[w^{p_i} = w^{q_i}$ and
$\la w^{p_i}, \a, \bar r^{p_i}, \G^{p_i} \ra \le
\la w^{q_i}, \b, \bar r^{q_i}, \G^{q_i} \ra]$.
\item{3.} $|\{j \in \underset i < \d \to{\cup}
w^{p_i}$ : For the unique $i$ so that $j \in
w^{p_i} = w^{q_i}$, $r^{p_i}_j \neq r^{q_i}_j \}|
< \d$, where $\bar r^{p_i} = \la r^{p_i}_j : j \in
w^{p_i} \ra$ and $\bar r^{q_i} = \la r^{q_i}_j :
j \in w^{q_i} \ra$.
\item{4.} $|\{z \in \underset i < \d \to{\cup}
\dom(\G^{p_i}) : \G^{p_i}(z) \neq \G^{q_i}(z) \}| < \d$.
 
\noindent Then, for $p = \la \la w^{p_i}, \a,
\bar r^{p_i}, \G^{p_i} \ra_{i < \d},
\la w^p, \a, \bar r^p, \G^p \ra \ra \in Q$, as
$w^{p_i} \cap w^{p_j} = \emptyset$ for $i < j < \d$
(each $w^{p_i} = y_{\b_i}$ for some unique $\b_i \in S$
and $y_{\b_i} \cap y_{\b_j} = \emptyset$ for
$\b_i \neq \b_j$), $w^{p_i} \cap w^p = \emptyset$ for
$i < \d$, $\dom(\bar r^{p_i}) \cap
\dom(\bar r^{p_j}) = \emptyset$ for $i < j < \d$,
$\dom(\bar r^{p_i}) \cap \dom(\bar r^p)
= \emptyset$ for $i < \d$, $\dom(\G^{p_i}) \cap
\dom(\G^{p_j}) = \emptyset$ for $i < j < \d$ (since if
$z \in \dom(\G^{p_i})$, $z = x_\b$ for some
$\b \in S$, meaning $w^{p_i} = y_\b$ by the definitions of
$P^1_{\d, \l}[S]$, $P^1_{\d, \l}[S] / R$, and $Q$), and
$\dom(\G^p) = \emptyset$ (since for every $\b \in S$,
$w^p \cap y_\b = \emptyset$, $y_\b \subseteq x_\b$, and
$x_\b - y_\b$ is bounded in $\b$), conditions 3) and 4)
above on the definition of $\le_Q$ show the function
$F(p) = \la \underset i < \d \to{\bigcup} w^{p_i}
\cup w^p, \a, \underset i < \d \to{\bigcup} \bar r^{p_i}
\cup \bar r^p, \underset i < \d \to{\bigcup} \G^{p_i} \ra$
yields an isomorphism between $Q$ and $P^1$. Thus, over
$V[G^0_{\d, \l}][G^2_{\d, \l}][H]$, forcing with
$P^1$, $P^1_{\d, \l}[S] / R$, and $Q$ are all equivalent.
 
We examine next in more detail the exact nature of
$\la Q, \le_Q \ra$. For
$\b \in S$, note that if
$p = \la w^p, \a^p, \bar r^p, \G^p \ra \in Q_\b$ and
$\dom(\G^p) \neq \emptyset$, then $\dom(\G^p) =
\{x_\b\}$. We can therefore define an ordering
$\le_\b$ on $Q_\b$ by
$p_1 = \la w^{p_1}, \a^{p_1}, \bar r^{p_1}, \G^{p_1} \ra \le_\b
p_2 = \la w^{p_2}, \a^{p_2}, \bar r^{p_2}, \G^{p_2} \ra$ iff
$p_1 = p_2$ or $p_1 < p_2$ and $\a^{p_1} \le
\max(\G^{p_2}(x_\b))$, and we can reorder $Q$ by replacing
each occurrence of $\le$ on $Q_\b$ by $\le_\b$. If we
call the new ordering on $Q$ thus obtained
$\le'_Q$, then by an argument virtually identical to the
one given in the remark following the proof of Lemma 3,
if $p_1 \le_Q p_2$, $p_1, p_2 \in Q$, there is some
condition $q \in Q$ so that $p_1 \le'_Q q$ and
$p_2 \le'_Q q$. It is hence once more the case that
$\la Q, \le_Q \ra$ and $\la Q, \le'_Q \ra$ are forcing
equivalent, i.e., $I$ is (appropriately) generic
with respect to $\la Q, \le_Q \ra$ iff $I$ is
(appropriately) generic with respect to
$\la Q, \le'_Q \ra$, so without loss of generality, we
analyze the forcing properties of $\la Q, \le'_Q \ra$.
 
We examine now $\la Q_\b , \le_\b \ra$ for $\b \in S$.
We first note that by the definition of $\le_\b$, for
$\g < \d$ any fixed but arbitrary cardinal, if
$\la p_i =
\la w^{p_i}, \a^{p_i}, \bar r^{p_i},          \G^{p_i} \ra
: i < \g \ra$
is a directed sequence of conditions with respect to
$\le_\b$, then (using the notation in the remark immediately
following the proof of Lemma 3) the 4-tuple $p' = \la
\underset i < \g \to{\cup} w^{p_i},
\underset i < \g \to{\cup} \a^{p_i},
\underset i < \g \to{\cup} \bar r^{p_i},
\underset i < \g \to{\cup} \G^{p_i} \ra$ can be
extended to a condition $p \in Q_\b$.
This is since the definition of $\le_\b$ ensures
$\underset i < \g \to{\cup} \a^{p_i} = \a' \le \eta$ for
$\eta = \max(\underset i < \g \to{\cup}
\G^{p_i})$. Thus, if we let $\a'' = \max(\a', \eta) + 1$,
we can extend each $r \in \underset i < \g \to{\cup}
\bar r^{p_i}$ to a function $s$ having domain $\a''$ by
letting $s \vert \a' = r \vert \a'$, and for
$\a \in [\a', \a'')$, $s(\a) =
(\bigcup H)(x_\b)$. If we call the sequence of all
such extensions $\bar r'$, $p =
\la \underset i < \g \to{\cup} w^{p_i}, \a'',
\bar r', \underset i < \g \to{\cup} \G^{p_i} \ra$
is a well-defined element of $Q_\b$ so that
$p_i \le_\b p$ for all $i < \g$. This just means
$\la Q_\b, \le_\b \ra$ is $\d$-directed closed.
 
Now, let $R_\b = \{ \la y_\b, \a, \bar r \ra : \a < \d$
and $\bar r = \la r_i : i \in y_\b \ra$ is a sequence of
functions from $\a$ to $\{0,1\} \}$, ordered by
$p_1 = \la y_\b, \a^1, \bar r^1 \ra \le_{R_\b}
p_2 = \la y_\b, \a^2, \bar r^2 \ra$ iff
$\a^1 \le \a^2$, $r^1_i \subseteq r^2_i$ for
$r^1_i \in \bar r^1$, $r^2_i \in \bar r^2$, and
$|\{i \in y_\b : r^1_i \neq r^2_i \}| < \d$.
Further, if $I_\b$ is
$V[G^0_{\d, \l}][G^2_{\d, \l}][H]$-generic
over $R_\b$, define in
$V[G^0_{\d, \l}][G^2_{\d, \l}][H][I_\b]$ an ordering
$T_\b$ having field $\{\G : \G$ is a function having
domain $\{x_\b\}$ and range $\{C \subseteq \d :
C$ is closed and bounded$\}\}$, ordered by
$\G^1 \le_{T_\b} \G^2$ iff
$\exists \la y_\b, \a^1, \bar r^1 \ra \in I_\b
\exists \la y_\b, \a^2, \bar r^2 \ra \in I_\b
[\la y_\b, \a^1, \bar r^1, \G^1 \ra \le_\b
\la y_\b, \a^2, \bar r^2, \G^2 \ra]$.
Since $R_\b$ is $\d$-directed closed,
$\{C \subseteq \d : C$ is closed and bounded$\}$ is
the same in either
$V[G^0_{\d, \l}][G^2_{\d, \l}][H]$ or
$V[G^0_{\d, \l}][G^2_{\d, \l}][H][I_\b]$.
This means $R_\b \ast \dot T_\b$ is isomorphic to
$\la Q_\b, \le_\b \ra$.
 
It is easy to see that the definition of $R_\b$ implies
$R_\b$ is isomorphic to $\K(\d, \d)$. Further, since
$\la Q_\b, \le_\b \ra$ is $\d$-directed closed and
$\la Q_\b, \le_\b \ra$ is isomorphic to
$R_\b \ast \dot T_\b$, $T_\b$ is $\d$-directed closed in
$V[G^0_{\d, \l}][G^2_{\d, \l}][H][I_\b]$.
Also, by its definition, $T_\b$ has cardinality $\d$ in
$V[G^0_{\d, \l}][G^2_{\d, \l}][H][I_\b]$, i.e.,
since $T_\b$ is $\d$-directed closed, $R_\b \ast \dot T_\b$ is
isomorphic to $\K(\d, \d) \ast \dot \C(\d)$, i.e.,
$\la Q_\b, \le_\b \ra$ is isomorphic to
$\K(\d, \d)$.
Since $|\{  x_\b : \b \in S \} | = \l$, conditions 3) and 4)
on the definition of $\le_Q$ ensure the ordering composed of
those elements of $\underset \b \in S \to{\prod} Q_\b$
having support $\d$ ordered by $\le'_Q \vert
\underset \b \in S \to{\prod} Q_\b$ is isomorphic to
$\K(\d, \l)$. Then, if $\la w, \a, \bar r, \G \ra \in Q^*$,
since we have already observed $\dom(\G)
= \emptyset$, $Q^*$ can easily be seen to be isomorphic to
$\K(\d, \l)$. Putting all of this together yields $Q$
ordered by $\le'_Q$ is isomorphic to $\K(\d, \l)$.
Thus, since $\la Q, \le_Q \ra$ and
$\la Q, \le'_Q \ra$ are forcing equivalent,
this proves Lemma 4.
 
\noindent
\hfill $\square $  Lemma 4
 
We remark that in what follows, it will frequently be the
case that a partial ordering $P$ is forcing equivalent to a
partial ordering $P'$ in the sense that a generic object for
one generates a generic object for the other. Under these
circumstances, we will often abuse terminology somewhat by saying
that either $P$ or $P'$ satisfies a certain chain condition,
a certain amount of closure, etc$.$ when this is true of at
least one of these partial orderings. We will then as
appropriate further
compound the abuse by using this property
interchangeably for either partial ordering.
 
As the definition of $\K(\k_0, \k_1)$ indicates, without
the last coordinates $Z^p$ and $\G^p$ of a condition
$p \in P^1_{\d, \l}[S]$ and the associated restrictions
on the ordering, $P^1_{\d, \l}[S]$ is essentially
$\K(\d, \l)$. These last coordinates and change in the
ordering are necessary to destroy the measurability of
$\d$ when forcing with $P^1_{\d, \l}[S]$. Once the fact
$S$ is stationary has been destroyed by forcing with
$P^2_{\d, \l}[S]$, Lemma 4 shows that these last two
coordinates $Z^p$ and $\G^p$ of a condition $p \in
P^1_{\d, \l}[S]$ can be factored out to produce the ordering
$\C(\d^+, \l) \ast \dot \K(\d, \l)$.
 
$\K(\d, \l)$, although somewhat similar in nature to
$\C(\d, \l)$ (e.g., $\K(\d, \l)$ is $\d$-directed closed),
differs from $\C(\d, \l)$ in a few very important aspects.
In particular, as we shall see presently, forcing with
$\K(\d, \l)$ will collapse $\d^+$. An indication that this
occurs is provided by the next lemma.
 
\proclaim{Lemma 5} $\K(\d, \l)$ satisfies
$\d^{++}$-c.c. whenever $2^\d = \d^+$.
\endproclaim
 
\demo{Proof of Lemma 5}
Suppose $\la \la w^\b, \a^\b, \bar r^\b \ra :
\b < \d^{++} \ra$ is a sequence of $\d^{++}$ many
incompatible elements of $\K(\d, \l)$. Since $\l > \d^+$,
each $w^\b \in {[\l]}^{< \d^+}$, $2^\d = \d^+$,
and $\l$ is either inaccessible or the successor of
a cardinal of cofinality $> \d$,
we can find some $A \subseteq \d^{++}$,
$|A| = \d^{++}$ so that $\{w^\b : \b \in A\}$
forms a $\Delta$-system, i.e.,
$\{w^\b : \b \in A\}$ is so that for
$\b_1 \neq \b_2 \in A$, $w^{\b_1} \cap w^{\b_2}$ is
some constant value $w$.
Since each $\a^\b < \d$, let $B \subseteq A$,
$|B| = \d^{++}$ be so that for $\b_1 \neq \b_2 \in B$,
$\a^{\b_1} = \a^{\b_2} = \a$. Since for $\b \in B$,
$\bar r^\b \vert w$ is a sequence of functions from
$\a < \d$ into $\{0, 1\}$, the facts $|w| \le \d$
and $2^\d = \d^+$ together imply there is
$C \subseteq B$, $|C| = \d^{++}$ so that for
$\b_1 \neq \b_2 \in C$, $\bar r^{\b_1} \vert w =
\bar r^{\b_2} \vert w$. It is thus the case that
$\la \la w^\b, \a^\b, \bar r^\b \ra : \b \in C \ra$
is now a sequence of $\d^{++}$ many compatible
elements of $\K(\d, \l)$, a contradiction.
This proves Lemma 5.
 
\noindent
\hfill $\square$ Lemma 5
 
It is clear from Lemmas 4 and 5 and the definition of
$\K(\d, \l)$ that since GCH holds in $V$ for cardinals
$\k \ge \d$, $P^0_{\d, \l} \ast
(P^1_{\d, \l}[\dot S] \times P^2_{\d, \l}[\dot S])$,
being equivalent to $\C(\l) \ast
\dot C(\d^+, \l) \ast \dot K(\d, \l)$, preserves
cardinals and cofinalities $\le \d$ and $\ge \d^{++}$,
has a dense subset which is $\d$-directed closed,
satisfies $\l^+$-c.c., and is so that
$V^{P^0_{\d, \l} \ast (P^1_{\d, \l}[\dot S] \times
P^2_{\d, \l}[\dot S])} \models ``$For every cardinal
$\k \in [\d, \l)$, $2^\k = \l$''. Our next lemma shows
that the forcing $P^0_{\d, \l}
\ast P^1_{\d, \l}[\dot S]$ is also rather nice, with the
exception that it collapses $\d^+$. By Lemma 4, this
has as an immediate consequence that the forcing
$\K(\d, \l)$ also collapses $\d^+$.
\proclaim{Lemma 6} $P^0_{\dell} * P^1_{\dell} [\dot S]$
preserves cardinals and cofinalities
$\le \d$ and $\ge \d^{++}$, collapses $\d^+$,
is $< \d$-strategically closed, satisfies $\l^+$-c.c.,
and is so that
$V^{                P^0_{\dell} * P^1_{\dell} [\dot S]} \models
``2^\k = \l$ for all cardinals $\k \in [\d, \l)$''.
\endproclaim
\demo{Proof of Lemma 6}
Let $G        ' = G        ^0_{\dell} * G        ^1_{\dell}$ be $V$-generic
over $P^0_{\dell} * P^1_{\dell}[\dot S]$, and let $G        ^2_{\dell}$
 be $V[G']$-generic over $P^2_{\dell} [S]$.
Thus, $G        ' *       G  ^2_{\dell} =       G$ is $V$-generic over $P^0_{\dell} * (P^1_{\dell}
[\dot S] * P^2_{\dell} [\dot S]) =
P^0_{\dell} * (P^1_{\dell} [\dot S] \times P^2_{\dell} [\dot S])$.
 
By Lemmas 4 and 5 and the remarks immediately following,
since GCH holds in $V$ for cardinals $\k \ge \d$,
$V[G] \models ``2^\k = \l$ for all cardinals $\k \in
[\d, \l)$'' and has the same cardinals and cofinalities
as $V$ $\le \d$ and $\ge \d^{++}$.
Hence, since $V[G'] \subseteq V[G]$, forcing with $P^0_{\dell}
 * P^1_{\dell} [\dot S]$ over $V$ preserves cardinals and
cofinalities $\le \d$ and $\ge \d^{++}$
and is so that
$V^{                P^0_{\dell} * P^1_{\dell} [\dot S]} \models
``2^\k = \l$ for all cardinals $\k \in [\d, \l)$''.
 
We now show forcing with $P^0_{\dell} \ast
P^1_{\dell}[\dot S]$ over $V$ collapses $\d^+$. Since forcing
with $P^0_{\dell}$ over $V$ collapses no cardinals and
preserves GCH for cardinals $\k \ge \d$, we assume without
loss of generality our ground model is $V[G^0_{\dell}] =
V_1$.
 
Using the notation of Lemma 3, i.e., that for $i < \l$,
$\dot r'_i$ is a term for $\cup\{r^p_i : \exists p =
\la w^p, \a^p, \bar r^p, Z^p, \G^p \ra \in
G^1_{\dell}[r^p_i \in \bar r^p]\}$, we can now define a
term $\dot \b_\zeta$ for $\zeta < \d$ by
$\dot \b_\zeta = \min(\{\b : \b < \d^+$ and $\d$
is the order type of
$\{i < \b : \dot r_i'(\zeta) = 1 \})$. To see that
$\dot \b_\zeta$ is well-defined, let $p =
\la w^p, \a^p, \bar r^p, Z^p, \Gamma^p \ra
\in P^1_{\dell}[S]$ be a condition. Without loss of
generality, we can assume that $\a^p > \zeta$.
Further, since $|w^p| = \d$,
$\sup(w^p \cap \d^+) < \d^+$, so we can let
$\g < \d^+$ be so that $\g > \sup(w^p \cap \d^+)$.
We can then define $q \ge p$, $q =
\la w^q, \a^q, \bar r^q, Z^q, \Gamma^q \ra
$ by letting $\a^q = \a^p$, $Z^q = Z^p$,
$\Gamma^q = \Gamma^p$, $w^q = w^p \cup
[\g, \g + \d)$, and $\bar r^q$ by $\bar r^q =
\la r^q_i : i \in w^q \ra$ where $r^q_i$ is
$r^p_i$ if $i \in w^p$, and $r^q_i$ is the function
having domain $\a^q$ which is constantly 1 if
$i \in w^q - w^p$. Clearly, $q$ is well-defined, and
$p \le q$. Also, by the definition of $q$, since
$\zeta < \a^q$, $q \force ``\forall i \in
[\g, \g + \d) [\dot r_i'(\zeta) = 1]$''. This means
$q \force ``\dot \b_\zeta \le \g + \d$'', so
$\dot \b_\zeta$ is well-defined.
 
We will be done if we can show
$\force_{P^1_{\dell}[S]} ``\la \dot \b_\zeta : \zeta < \d \ra$
is unbounded in $\d^+$''. Assume now towards a contradiction
that $p =
\la w^p, \a^p, \bar r^p, Z^p, \Gamma^p \ra$ is so that
$p \force ``\sup(\la \dot \b_\zeta : \zeta <
\d \ra) = \sigma < \d^+$''. If we define $q =
\la w^q, \a^q, \bar r^q, Z^q, \Gamma^q \ra$ by
$w^q = w^p \cup \{i < \sigma : i \in \sigma -
w^p \}$, $\a^q = \a^p$, $\bar r^q = \la
r^q_i : i \in w^q \ra$ where $r^q_i = r^p_i$ if
$i \in w^p$ and $r^q_i$ is the function having domain
$\a^q$ which is constantly $0$ if $i \in w^q - w^p$,
$Z^q = Z^p$, and $\Gamma^q = \Gamma^p$, then as above,
$q$ is well-defined, and $p \le q$. We claim that for
$\zeta = \a^q$, $q \force ``\dot \b_\zeta > \sigma$''.
 
To see that the claim is true, let
$s \ge q$, $s =
\la w^s, \a^s, \bar r^s, Z^s, \Gamma^s \ra$ be so that
$\a^s > \a^q$. By clause 3 in the definition of
$\le$ on $P^1_{\dell}[S]$, $s \force ``
|\{i < \sigma : \dot r_i'(\zeta) = 1 \}| < \d$''.
This means $q \force ``\dot \b_\zeta > \sigma$'', thus
proving our claim and showing that $\d^+$ is collapsed.
 
We next show the $<\d$-strategic closure of
$P^0_{\delta, \l} * P^1_{\delta, \l}[\dot S]$. We first note that
as $(P^0_{\delta, \l} * P^1_{\delta, \l}[\dot S]) *
P^2_{\delta, \l}[\dot S] = P^0_{\delta, \l} *
(P^1_{\delta, \l}[\dot S] * P^2_{\delta, \l}[\dot S])$ has by
Lemma 4 a dense subset which is $<\d$-closed, the desired fact
follows from the more general fact that if $P * \dot Q$ is a
partial ordering with a dense subset $R$ so that $R$ is
$<\d$-closed, then $P$ is $<\d$-strategically closed. To show
this more general fact, let $\gamma < \l$ be a cardinal.
Suppose I and II play to build an increasing chain of elements of
$P$, with $\langle p_\beta : \beta \le \alpha + 1 \rangle$
enumerating all plays by I and II through an odd stage
$\alpha + 1$ and $\langle \dot{q_\beta} : \beta < \alpha +
1$ and $\beta$ is even or a limit ordinal$\rangle$ enumerating a
set of auxiliary plays by II which have been chosen so that
$\langle \langle p_\beta, \dot{q_\beta} \rangle : \beta <
\alpha + 1$ and $\beta$ is even or a limit ordinal$\rangle$
enumerates an increasing chain of elements of the dense
subset $R \subseteq  P * \dot Q$. At stage $\alpha + 2$, II
chooses $\langle p_{\alpha +2}, \dot q_{\alpha + 2}
\rangle$ so that $\langle p_{\alpha + 2}, \dot q_{\alpha+2}
\rangle \in R$ and so that $\langle p_{\alpha + 2},
\dot q_{\alpha+2} \rangle \ge \langle p_{\alpha + 1},
\dot{q_\alpha} \rangle$; this makes sense, since
inductively, $\langle p_\alpha, \dot q_\alpha \rangle
\in R \subseteq P * \dot Q$, so as I has chosen
$p_{\alpha + 1} \ge p_\alpha$, $\langle p_{\alpha + 1},
\dot q_\alpha \rangle \in P * \dot Q$. By the $<\d$-closure of
$R$, at any limit stage $\eta \le \gamma$, II can choose
$\langle p_\eta, \dot{q_\eta} \rangle$ so that
$\langle p_\eta, \dot{q_\eta} \rangle$ is an upper bound to
$\langle \langle p_\beta, \dot{q_\beta} \rangle : \beta <
\eta$ and $\beta$ is even or a limit ordinal$\rangle$. The
preceding yields a winning strategy for II, so $P$ is
$<\d$-strategically closed.
 
Finally, to show $P^0_{\dell} * P^1_{\dell} [\dot S]$ satisfies
$\l^+$-c.c.,
 we simply note that
this follows from the general fact about iterated forcing (see
[Ba]) that if $P * \dot Q$ satisfies $\l^+   $-c.c., then $P$
satisfies $\l^+   $-c.c. (Here, $P = P^0_{\delta, \l} *
P^1_{\delta, \l}[\dot S]$ and $Q = P^2_{\delta, \l}
[\dot S]$.)
This proves Lemma 6.
 
\noindent
\hfill $\square $  Lemma 6
 
We remark that
$\force_{P^0_{\dell}} ``P^1_{\dell} [\dot S]$ is $\d^{++}$-c.c.".
Otherwise,
if ${\cal A} = \langle p_\a : \a < \d^{++} \rangle$ were a size
$\d^{++}$ antichain of elements of $P^1_{\delta, \l}[S]$ in
$V[G^0_{\delta, \l}]$, then (using the notation of Lemma 4)
since $P^1_{\d, \l}[S]$ is isomorphic to
$\C(\d^+, \l) \ast (\dot P^1_{\d, \l}[S] / R)$ and
$P^1_{\d, \l}[S] / R$ has a dense subset which is isomorphic to
$\la Q, \le_Q \ra$, without loss of generality, $\A$ can be
taken as an antichain in $\la Q, \le_Q \ra$. Since
$\le'_Q \ \subseteq \ \le_Q$, $\A$ must also be an antichain
with respect to $\le'_Q$, and as
$V[G^0_{\d, \l}]$,
$V[G^0_{\d, \l}][G^2_{\d, \l}]$, and
$V[G^0_{\d, \l}][G^2_{\d, \l}][H]$
all have the same cardinals, $\A$ must be a size
$\d^{++}$ antichain with respect to $\le'_Q$ in
$V[G^0_{\d, \l}][G^2_{\d, \l}][H]$. The fact that
$V[G^0_{\d, \l}][G^2_{\d, \l}][H] \models ``2^\d =
\d^+$ and $\la Q, \le'_Q \ra$ is isomorphic to
$\K(\d, \l)$'' then tells us that $\A$ is isomorphic to
a size $\d^{++}$ antichain with respect to $\K(\d, \l)$ in
$V[G^0_{\d, \l}][G^2_{\d, \l}][H]$. Lemma 5, which says that
$\K(\d, \l)$ is $\d^{++}$-c.c$.$ in any model in which
$2^\d = \d^+$, now yields an immediate contradiction.
 
We conclude this section with a lemma that will be used
later in showing that it is possible to extend certain
elementary embeddings witnessing the appropriate degree of
supercompactness.
 
\proclaim{Lemma 7}
For $V_1 = V^{P^0_{\d, \l}}$, the models
$V^{P^1_{\d, \l}[S] \times P^2_{\d, \l}[S]}_1$ and
$V^{P^1_{\d, \l}[S]}_1$ contain the same $< \l$ sequences
of elements of $V_1$.
\endproclaim
 
\demo{Proof of Lemma 7}
By Lemma 4, since
$P^0_{\d, \l} \ast P^2_{\d, \l}[\dot S]$
is equivalent to the forcing $\C(\l)$ and
$V \subseteq
V^{P^0_{\d, \l}}
\subseteq
V^{P^0_{\d, \l} \ast P^2_{\d, \l}[\dot S]}
$, the models $V$, $
V^{P^0_{\d, \l}}
$, and $
V^{P^0_{\d, \l} \ast P^2_{\d, \l}[\dot S]}
$ all contain the same $< \l$ sequences of elements of $V$.
Thus, since a $< \l$ sequence of elements of $V_1 =
V^{P^0_{\d, \l}}
$ can be represented by a $V$-term which is
actually a function $h : \g \to V$
for some $\g < \l$, it immediately follows that
$
V^{P^0_{\d, \l}}
$ and $
V^{P^0_{\d, \l} \ast P^2_{\d, \l}[\dot S]}$
contain the same $< \l$ sequences of elements of $
V^{P^0_{\d, \l}}
$.
 
Let now $f : \g \to V_1$ for $\g < \l$ be so that $f \in (
V^{P^0_{\d, \l} \ast P^2_{\d, \l}[\dot S]}
){}^{P^1_{\d, \l}[S]} =
 V^{P^1_{\d, \l}[S] \times P^2_{\d, \l}[S]}_1$, and let
$g : \g \to V_1$, $g \in
V^{P^0_{\d, \l} \ast P^2_{\d, \l}[\dot S]}
$ be a term for $f$. By the previous paragraph, $g \in
V^{P^0_{\d, \l}}
$. Since Lemma 5 shows that $P^1_{\d, \l}[S]$ is
$\d^{++}$-c.c$.$ in $
V^{P^0_{\d, \l} \ast P^2_{\d, \l}[\dot S]}
$ and $\d^{++} \le \l$, for each $\a < \g$, the antichain $\A
_\a$ defined in $
V^{P^0_{\d, \l} \ast P^2_{\d, \l}[\dot S]}
$ by $\{p \in P^1_{\d, \l}[S] : p$ decides a value for
$g(\a) \}$ is so that $
V^{P^0_{\d, \l} \ast P^2_{\d, \l}[\dot S]}
\models ``|\A_\a| <   \l$''. Hence, by the preceding
paragraph, since $\A_\a$ is a set of elements of $
V^{P^0_{\d, \l}}
$, $\A_\a \in
V^{P^0_{\d, \l}}
$ for each $\a < \g$. Therefore, again by the preceding
paragraph, the sequence $\la \A_\a : \a < \g \ra \in
V^{P^0_{\d, \l}}
$. This just means that the term $g \in
V^{P^0_{\d, \l}}
$ can be evaluated in $
V^{P^1_{\d, \l}[S]}_1
$, i.e., $f \in
V^{P^1_{\d, \l}[S]}_1
$. This proves Lemma 7.
 
\no \hfill $\square$ Lemma 7
 
\S 2 The Proof of Theorem 1
 
We turn now to the proof of Theorem 1. Recall that we are assuming
our ground model $V \models ``$ZFC + GCH + $\k$ is $< \l$
supercompact + $\l > \k^+$ is regular
and is either inaccessible or is the successor of a cardinal
of cofinality $> \k$ +
$h: \k \to \k$ is so that
for some elementary embedding $j:V \to M$ witnessing the $< \l$
supercompactness of $\k$, $j(h)(\k) = \l$''. By reflection,
we can assume without loss of generality that for every
inaccessible $\d < \k$, $h(\d) > \d^+$ and $h(\d)$ is
regular. Given this,
we are now in a position to define the partial ordering $P$
used in the proof of Theorem 1.
We define a $\k$ stage Easton support iteration $P_\k = \langle \langle P_\a,
\dot Q_\a \rangle : \a < \k \rangle$, and then define $P = P_{\k + 1} = P_\k *
\dot Q_\k
$ for a certain partial ordering $Q_\k$ definable in $V^{P_\k}$.
The definition is as follows:
\item{1.} $P_0$ is trivial.
\item{2.} Assuming $P_\a$ has been defined for $\a < \k$,
let $\d_\a$ be so that $\d_\a$ is the least cardinal
$\ge \underset \b < \a \to{\cup} \d_\b$ such that
$\force_{P_\a} ``\d_\a$ is inaccessible'', where $\d_{-1} = 0$.
Then $P_{\a + 1} = P_\a \ast \dot Q_\a$, with
$\dot Q_\a$ a term for $P^0_{\d_\a, h(\d_\a)} \ast
P^1_{\d_\a, h(\d_\a)}[\dot S_{h(\d_\a)}]$, where
$\dot S_{h(\d_\a)}$ is a term for the non-reflecting
stationary subset of $h(\d_\a)$ introduced by
$P^0_{\d_\a, h(\d_\a)}$.
\item{3.} $\dot Q_\k$ is a term for
$P^0_{\k, \l} \ast (P^1_{\k, \l}[\dot S_\l] \times
P^2_{\k, \l}[\dot S_\l])$, where again, $\dot S_\l$ is a term for
the non-reflecting stationary subset of $\l$ introduced
by $P^0_{\k, \l}$.
 
 The intuitive motivation behind the above definition is that below $\k$
 at any inaccessible $\d$, we must force to ensure that
$\d$ becomes non-measurable and is so that $2^\d = h(\d)$.
At $\k$, however, we must force so as simultaneously to make
$2^\k = \l$ while first destroying and then resurrecting the
$\d$ supercompactness of $\k$ for all regular $\d \in [\k, \l)$.
 
\proclaim{Lemma 8}
$V^P \models ``$For all
inaccessible $\d < \k$ and all cardinals $\g \in [\d,
h(\d))$, $2^\g = h(\d)$, for all cardinals $\g \in
[\k, \l)$, $2^\g = \l$,
and no cardinal $\d < \k$
is measurable''.
\endproclaim
\demo{Proof of Lemma 8}
We first note that Easton support iterations of
$\delta$-strategically
closed partial orderings are $\delta$-strategically closed for $\delta$ any
regular cardinal.
The proof is via induction.
If $R_1$ is $\delta$-strategically closed and $\force_{R_1} ``\dot R_2$ is
$\delta$-strategically
 closed", then let $ p \in R_1$ be so that $p \force ``\dot g$ is a
strategy for
player II ensuring that the game which produces an increasing chain of
elements of
$\dot R_2$ of length $\delta $ can always be continued for
 $\a \le \delta"$.
If II begins by picking $r_0 = \langle p_0, \dot q_0 \rangle \in R_1 \ast
\dot R_2$ so that $p_0 \ge p$ has been chosen
according to the strategy $f$ for $R_1$ and $p_0 \force
``\dot q_0$ has been chosen according to $\dot g"$, and at even stages
$\a + 2$
picks $r_{\a + 2} = \langle p_{\a + 2}, \dot q_{\a+ 2} \rangle$ so that
$p_{\a + 2}$
 has been chosen according to $f$ and is so that $p_{\a + 2} \force
 ``\dot q_{\a + 2}$ has been chosen according to $\dot g$",
 then at limit stages $\l \le \delta$, the chain
$r_0 = \langle p_0, \dot q_0 \rangle \le r_1 = \langle p_1, \dot q_1 \rangle
\le
 \cdots \le r_\a = \langle p_\a, \dot q_\a \rangle \le \cdots (\a < \l)$ is so
that II can find an upper bound $p_\l$ for
$\langle p_\a: \a < \l \rangle $ using $f$.
By construction, $p_\l \force `` \langle \dot q_\a : \a < \l \rangle  $ is so
 that at
 limit and even stages, II has played according to $\dot g$",  so for some
$\dot q_\l$, $p_\l \force `` \dot q_\l$ is an upper bound to $\langle \dot q_\a:
 \a
< \l \rangle "$,  meaning the condition $\langle p_\l, \dot q_\l \rangle$ is as
 desired.
These methods, together with the usual proof at limit stages (see [Ba], Theorem
 2.5)
 that the Easton support iteration of $\delta$-closed partial orderings is
 $\delta$-closed,
 yield that $\delta$-strategic closure is preserved at limit stages
 of any Easton support iteration of $\delta$-strategically closed
partial orderings.
In addition, the ideas of this paragraph will also show
that Easton support iterations of $\prec \d  $-strategically
closed partial orderings are $\prec \d  $-strategically
closed for $\d$ any regular cardinal.
 
Given this fact, it is now easy to prove by induction that
$V^{P_\k} \models ``$For all inaccessible $\d < \k$ and
all cardinals $\g \in [\d, h(\d))$, $2^\g = h(\d)$,
and no cardinal $\d < \k$ is measurable''.
Given $\a < \k$, we assume inductively that
$V^{P_\a} \models ``$For all inaccessible $\d < \d_\a$ and
all cardinals $\g \in [\d, h(\d))$, $2^\g = h(\d)$,
and no cardinal $\d < \d_\a$ is measurable'',
where $\d_\a$ is as in the definition of $P$.
By Lemmas 3 and 6 and the definition of $P$, since inductively
$\force_{P_\a} ``$GCH holds for all cardinals $\d \ge \d_\a$'',
$V^{P_\a \ast \dot Q_\a} = V^{P_{\a + 1}} \models ``$For
all inaccessible $\d \le \d_\a$ and
all cardinals $\g \in [\d, h(\d))$, $2^\g = h(\d)$,
and no cardinal $\d \le \d_\a$ is measurable''.
 
If now $\b \le \k$ is a limit ordinal, then we know by
induction that for all $\a < \b$,
$V^{P_\a} \models ``$For all inaccessible $\d < \d_\a$ and
all cardinals $\g \in [\d, h(\d))$, $2^\g = h(\d)$,
and no cardinal $\d < \d_\a$ is measurable''.
If we write $P_\b = P_\a \ast \dot R$, then by the definition
of $P$, the proof of Lemma 4, Lemma 6, and the fact contained in
the first paragraph of the proof of this lemma, $\force_{P_\a}
``\dot R$ is forcing equivalent to a $< \d_\a$-strategically
closed partial ordering'', so
$V^{P_\a \ast \dot R} = V^{P_\b} \models ``$For all
inaccessible $\d < \d_\a$ and all cardinals
$\g \in [\d, h(\d))$, $2^\g = h(\d)$, and no cardinal
$\d < \d_\a$ is measurable''.
If we let $\d_\k = \k$, since $\a < \b$ is arbitrary
in the preceding, it
thus follows by the definition of $\d_\a$ for $\a < \k$ that
$V^{P_\k} \models ``$For all inaccessible $\d < \d_\k$
and all cardinals $\g \in [\d, h(\d))$, $2^\g = h(\d)$,
and no cardinal $\d < \d_\k$ is
measurable''.
 
The proof of Lemma 8 will be complete once we show
$V^{P_\k \ast \dot Q_\k} = V^P$
is so that $V^P \models ``$For all
inaccessible $\d < \k$ and all cardinals $\g \in [\d, h(\d))$,
$2^\g = h(\d)$, for all cardinals $\g \in [\k, \l)$,
$2^\g = \l$,
and no cardinal $\d < \k$ is measurable''. By the last
paragraph, $\force_{P_\k} ``\k$ is inaccessible'', and by
the definition of $P_\k$, $|P_\k| = \k$, meaning
$\force_{P_\k} ``$GCH holds for all cardinals $\d \ge
\k$''. Therefore, by Lemma 4 and the definition of
$\dot Q_\k$, $V^{P_\k} \models ``Q_\k$ is equivalent to
$\C(\l) \ast \dot \C(\k^+, \l) \ast \dot \K(\k, \l)$'',
so
$V^P \models ``$For all inaccessible
$\d < \k$ and all cardinals $\g \in [\d, h(\d))$,
$2^\g = h(\d)$, for all cardinals $\g \in [\k, \l)$,
$2^\g = \l$, and no cardinal
$\d < \k$ is measurable''.
This proves Lemma 8.
\pbf \hfill $\square $ Lemma 8
 
We now show that the intuitive motivation for the definition of $P$ as set forth
in the paragraph immediately preceding the statement of Lemma 8 actually
works.
\proclaim{Lemma 9} For $G$ $V$-generic over $P$,
$V[G] \models ``\k$ is $<\l$ supercompact''.
\endproclaim
\demo{Proof of Lemma 9}
Let $j:V \to M$ be an elementary embedding witnessing the
$< \l$ supercompactness of $\k$ so that
$j(h)(\k) = \l$.
We will actually show that for $G =    G _{\k} \ast G        '_{\k}$
our $V$-generic object over $P = P_{\k} \ast \dot Q_{\k}$, the
embedding $j$ extends to $k : V[G_\k \ast G'_\k] \to
M[H]$ for some
$    H   \subseteq j(P)$.
As $\langle j (\a): \a < \gamma \rangle \in M$
for every $\g < \l$, this will be enough
to allow for every $\g < \l$
the definition of the ultrafilter $x \in {\cal U}_\g$ iff $
\langle j (\a ) : \a < \gamma \rangle \in  k(x) $ to be given in $V[G_{\k}
 \ast
       G  '_{\k}]$,
thereby showing $V[G] \models ``\k$ is $<\l$ supercompact''.
 
We construct $H$ in stages.
In $M$, as $\k$ is the critical point of
 $j$, $j(P_{\k} \ast \dot Q_{\k}) = P_{\k} \ast \dot R'_{\k} \ast
\dot R''_{\k} \ast \dot R'''_{\k}$, where $\dot R'_{\k} $ will be a term
for $        P^0_{\k, \l} \ast
 P^1_{\k,     \l} [\dot S_\l]$
(note that as
$M^{< \l} \subseteq M$, $M \models ``\d_\k = \k$'',
$j(\k) > \k$, and $j(h)(\k) = \l$,
$\dot R'_{\k}$ is indeed as just stated),
$ \dot R''_{\k}$ will be a term for
the rest of the portion of $j(P_{\k})$ defined below $j(\k)$, and $\dot R'''_\k$
will be a term for
 $j (\dot Q_{\k})$.
This will allow us to define $    H$ as $H        _{\k} \ast H        '_{\k} \ast       H  ''_{\k}
\ast       H  '''_{\k}$.
Factoring $G        '_{\k} $ as $ G^0_{\k, \l} \ast (       G  ^1_{\k, \l} \times       G  ^2_{\k,
 \l})$,
we let $    H    _{\k} =       G  _{\k}$
and
$H'_\k = G^0_{\k, \l} \ast G^1_{\k, \l}$.
Thus, $    H    '_{\k}$ is the same as $G        '_{\k}$, except, since
$M \models ``\k < j(\k)$ and $j(h)(\k) = \l$'', we omit the generic object
$G        ^2_{\k, \l        }$.
 
To construct $    H    ''_{\k}$, we first note that the definition of
$P$ ensures $|P_{\k}| = \k$ and, since $\k$ is necessarily Mahlo, $P_{\k}$ is $\k$-c.c.
As $V[G_{\k}]$ and $M[G_{\k}]$ are both  models of GCH
for cardinals $\g \ge \k$, the definition of
$R'_{\k}$ in $M[H_{\k}]$ and the remark following Lemma 6
then ensure that
$M[H_{\k}] \models ``R'_{\k}$ is
a $< \l$-strategically closed partial ordering followed by a
$\k^{++}$-c.c$.$ partial ordering and
$\k^{++} \le \l$''.
Since $M^{< \l} \subseteq M$
implies cardinals in $V$ $\le \l$ are the same
as cardinals in $M$ $\le \l$ and
$P_\k    $ is $\k    $-c.c., Lemma 6.4 of [Ba] shows
$V[G_\k]$ satisfies these facts as well.
This means $< \l$-strategic closure and the argument of
Lemma 6.4 of [Ba] can be applied to show $M[H_{\k} *
H_{\k}'] = M[G_{\k} * H_{\k}']$ is
closed under $ < \l $ sequences with respect to
$V[G_{\k} * H_{\k}']$,
i.e., if $\g < \l$ is a cardinal,
$f : \gamma \to M[H_\k \ast       H'_\k]  $, $f \in V[G_\k \ast       H  '_\k]$,
 then $f \in M[H_{\k} \ast       H  '_{\k}].$
Therefore, as Lemma 8 tells us $M[H_\k \ast H'_\k] \models
``R''_\k$ is forcing equivalent to a $\prec \l$-strategically
closed partial ordering'',
this fact is true in $V[G_{\k} \ast
       H  '_{\k}]$  as well.
 
Observe now that GCH in $V$ allows us to assume $\l < j (\k) < j(\k^+)
< \l^+       $.
Since $M[H_{\k} \ast       H  '_{\k} ] \models `` |R''_\k    | = j (\k)
$ and $|{\cal P} (R''_{\k}) | = j(\k^+)$" (this last fact follows from GCH
 in $M[H_{\k } \ast       H  '_\k      ] $ for cardinals $\g \ge \l$), in $V[G_{\k} \ast       H  '_{\k}]$,
 we can let $\langle D_\a: \a < \l \rangle $ be an enumeration of the
dense open subsets of $R''_{\k} $ present in $M[H_{\k} \ast       H  '_{\k}]$.
The $\prec\l$-strategic closure of $R''_{\k}$
in both $M[H_{\k} \ast       H  '_{\k}]$ and $V[G_{\k} \ast       H  '_{\k}]$
 now allows us to meet all of these  dense subsets as follows.
Work in $V[G_{\k} \ast        H  '_{\k}]$.
Player I picks $p_\a \in D_\a $ extending $\sup ( \langle q_\b: \b < \a \rangle
 )$
 (initially, $q_{-1}$ is the trivial condition), and player II responds by
 picking
$q_\a \ge p_\a$ (so $q_\a \in D_\a$).
By the $\prec\l$-strategic closure of
$R''_{\k}$ in $V[G{}_{\k} \ast       H
{}'_{\k}]$,
player II has a winning strategy for this game,
so $\langle q_\a: \a < \l \rangle$ can be taken as
an increasing sequence of conditions with $q_\a \in D_\a$ for
$\a <
 \l$.
Clearly, $    H    ''_{\k} = \{ p \in R''_{\k} : \exists \a < \l [q_\a \ge p]
 \}$
 is our $M [H_{\k} \ast       H  '_{\k}]$-generic object over $R''_{\k}$ which has
 been
constructed in $V[G_{\k} \ast       H  '_{\k}] \subseteq V[G_{\k} \ast
       G  '_{\k}]$,
 so $    H    ''_{\k} \in V[G_{\k} \ast       G  '_{\k}]$.
 
By the  above construction, in $V[G{}_{\k} \ast       G    '_{\k}]$,
the embedding $j$ extends to an
embedding $j^* : V[G{}_{\k}] \to M[H_{\k} \ast H'_{\k} \ast H''_{\k}]$.
We will be done once we have constructed in $V[G{}_{\k} \ast       G    '_{\k}]$
 the appropriate generic object for $R'''_{\k} = P^0_{j(\k), j(\l)} * (P^1_{j(\k), j
 (\l)}
[ \dot S_{j(\l)}] \times P^2_{j(\k), j(\l)} [\dot S_{j(\l)} ]) =
 (P^0_{j(\k), j(\l)} \ast
 P^2_{j(\k), j(\l)} [\dot S_{j(\l)}]) \ast P^1_{j(\k), j (\l)} [\dot S_{j
 (\l)}]$.
To do this, first rewrite $G        '  _{\k}$ as $(G{}^0_{\k, \l} \ast       G    ^2_{\k, \l})
 \ast       G    ^1_{\k , \l}$.
By the nature of the forcings, $G^0_{\k, \l} \ast G^2_{\k, \l}$
is $V[G_\k]$-generic over a partial ordering which is
$(< \l, \infty)$-distributive. Thus, by a general fact about
transference of generics via elementary embeddings
(see [C], Section 1.2, Fact 2, pp$.$ 5-6), since
$j^* : V[G_\k] \to
M[H_\k \ast H'_\k \ast H''_\k]
$
is so that every element of $
M[H_\k \ast H'_\k \ast H''_\k]
$ can be written $j^*(F)(a)$ with $\dom(F)$ having
cardinality $< \l$, ${j^{*}}''G^0_{\k, \l} \ast
G^2_{\k, \l}$ generates an $
M[H_\k \ast H'_\k \ast H''_\k]
$-generic set $H^4_\k$.
 
It remains to construct $H^5_{\k}$, our $M                        [H_{\k} \ast H'_{\k} \ast
 H''_{\k}
\ast H^4_{\k}]$-generic object over \break $P^1_{j(\k), j(\l)}
 [S_{j(\l    )}]$.
To do this, first recall that in $M[H_\k \ast H'_\k]$, as
previously noted, $R''_\k$ is $\prec \l$-strategically
closed. Since $M[H_\k \ast H'_\k]$ has already been observed to
be closed under $< \l$ sequences with respect to
$V[G_\k \ast H'_\k]$, and since any $\g$ sequence of elements
for $\g < \l$ a cardinal of $M[H_\k \ast H'_\k \ast H''_\k]$
can be represented, in $M[H_\k \ast H'_\k]$, by a term which is
actually a function $f: \g \to M[H_\k \ast H'_\k]$,
$M[H_\k \ast H'_\k \ast H''_\k]$ is closed under $< \l$
sequences with respect to $V[G_\k \ast H'_\k]$, i.e., if
$f : \g \to M[H_\k \ast H'_\k \ast H''_\k]$ for $\g < \l$ a
cardinal, $f \in V[G_\k \ast H'_\k]$, then
$f \in M[H_\k \ast H'_\k \ast H''_\k]$.
 
Choose in $V[G_\k \ast G'_\k]$ an enumeration
$\la p_\a : \a < \l   \ra$ of $G^1_{\k, \l}$. Adopting the
notation of Lemma 4 and working now in
$V[G_\k \ast G'_\k]$, first note that by Lemma 4, there is
an isomorphism between $P^1_{\k, \l}[S_\l]$ and
$R \ast (\dot P^1_{\k, \l}[S_\l] / R)$. Again by Lemma 4, since
$P^1$ is dense in $P^1_{\k, \l}[S_\l] / R$ and
$\la Q, \le_Q \ra$ and $P^1$ are isomorphic, there is an
isomorphism $g_0$ between a dense subset of
$P^1_{\k, \l}[S_\l]$ and $R \ast \la \dot Q, \le_Q \ra$.
Therefore, as $R$ is isomorphic to $\C(\k^+, \l)$ and
$\la Q, \le'_Q \ra$ is isomorphic to $\K(\k, \l)$, we
can let $g_1 : R \ast \la \dot Q, \le'_Q \ra \to
\C(\k^+, \l) \ast \dot \K(\k, \l)$ be an isomorphism.
It is then the case that $g = g_1 \circ g_0$ is a bijection
between a dense subset of $P^1_{\k, \l}[S_\l]$ and
$\C(\k^+, \l) \ast \dot \K(\k, \l)$.
This gives us a sequence $I =
 \la g(p_\a) : \a < \l \ra$ of $\l$ many compatible elements of
$\C(\k^+, \l) \ast \dot \K(\k, \l)$. By Lemma 7,
$V[G_\k \ast G^0_{\k, \l} \ast G^1_{\k, \l} \ast
G^2_{\k, \l}] = V[G_\k \ast G^0_{\k, \l} \ast
G^2_{\k, \l} \ast G^1_{\k, \l}] = V[G_\k \ast G'_\k]$
and $V[G_\k \ast G^0_{\k, \l} \ast G^1_{\k, \l}] =
V[G_\k \ast H'_\k]$ have the same $< \l$ sequences of
elements of $V[G_\k \ast G^0_{\k, \l}]$
and hence of $V[G_\k \ast H'_\k]$.
Thus, any $< \l$ sequence of
elements of
$M[H_\k \ast H'_\k \ast H''_\k]$
present in $V[G_\k \ast G'_\k]$ is actually an element of
$V[G_\k \ast H'_\k]$ (so
$M[H_\k \ast H'_\k \ast H''_\k]$
is really closed under $< \l$ sequences with respect to
$V[G_\k \ast G'_\k]$.)
 
For $\a \in (\k^+, \l)$, if $p \in
 \C(\k^+, \l)$, let
$p \vert \a = \{\la \la \rho, \sigma \ra, \eta \ra \in p :
\sigma < \a \}$, and if $q \in \K(\k, \l)$, $q =
\la w, \sigma, \bar r \ra$, let $q \vert \a = \la
w \cap \a, \sigma, \bar r \vert (w \cap \a) \ra$.
Call $w$ the domain of $q$. Since
$\force_{C(\k^+, \l)} ``$There are no new $\k$ sequences of
ordinals'', we can assume without loss of generality that
for any condition $p = \la p^0, p^1 \ra \in
 \C(\k^+, \l) \ast \dot \K(\k, \l)$, $p^1$ is an actual
condition and not just a term for a condition. Thus,
for $\a \in (\k^+, \l)$ and
$p = \la p^0, p^1 \ra \in
 \C(\k^+, \l) \ast \dot \K(\k, \l)$, the definitions
$p \vert \a = \la p^0 \vert \a, p^1 \vert \a \ra$ and
$I \vert \a = \{p \vert \a : p \in I \}$ make sense.
And, since GCH holds in
$V[G_\k \ast G^0_{\k, \l} \ast G^2_{\k, \l}]$ for cardinals
$\g \ge \k$, it is clear $V[G_\k \ast G'_\k] \models
``|I \vert \a| < \l$ for all $\a \in (\k^+, \l)$''.
Thus, since
$\C(j(\k^+), j(\l)) \ast \dot \K(j(\k), j(\l)) \in
M[H_\k \ast H'_\k \ast H''_\k]$
($\C(\k^+, \l) \ast \dot \K(\k, \l) \in V[G_\k]$) and
$M[H_\k \ast H'_\k \ast H''_\k] \models ``
 \C(j(\k^+), j(\l)) \ast \dot \K(j(\k), j(\l))$ is
$j(\k)$-directed closed'', the
facts $M[H_\k \ast H'_\k \ast H''_\k]$ is closed under
$< \l$ sequences with respect to $V[G_\k \ast G'_\k]$
and $I$ is compatible imply that $q_\a = \la
q^0_\a, q^1_\a \ra = \bigcup \{j^*(p) : p \in
I \vert \a \}$ for $\a \in (\k^+, \l)$ is well-defined
and is an element of
$\C(j(\k^+), j(\l)) \ast \dot \K(j(\k), j(\l))$.
 
Letting $q^1_\a = \la w^\a, \sigma^\a, \bar r^\a \ra$,
the definition of $I \vert \a$ and the elementarity of
$j^*$ easily imply that if $\rho \in w^\a -
\underset \b < \a \to{\cup} w^\b =
\dom(q^1_\a) - \dom(\underset \b < \a \to{\cup}
q^1_\b)$, then $\rho \in [
\underset \b < \a \to{\cup} j(\b), j(\a))$. Also, by
the fact $M[H_\k \ast H'_\k \ast H''_\k]$ is closed
under $< \l$ sequences with respect to $V[G_\k \ast G'_\k]$,
$\underset \b < \a \to{\cup} q^1_\b \in
\K(j(\k), j(\l))$ and $\underset \b < \a \to{\cup}
q^0_\b \in  \C(j(\k^+), j(\l))$, i.e.,
$\underset \b < \a \to{\cup} q_\b \in
 \C(j(\k^+), j(\l)) \ast \dot \K(j(\k), j(\l))$.
Further, if $\la \rho, \sigma \ra \in
\dom(q^0_\a) - \dom(\underset \b < \a \to{\cup}
q^0_\b)$, then again as before, $\sigma \in
[\underset \b < \a \to{\cup} j(\b), j(\a))$. This is
since if $\sigma < \underset \b < \a \to{\cup}
j(\b)$, then let $\b$ be minimal so that $\sigma <
j(\b)$, and let $\rho$ and $\sigma$ be so that
$\la \rho, \sigma \ra \in \dom(q^0_\a)$. It follows
that for some $r = \la r^0, r^1 \ra \in
I \vert \a$, $\la \rho, \sigma \ra \in
\dom(j^*(r^0))$. Since by elementarity and the
definitions of $I 
\vert \b$ and $I \vert \a$, for $r^0 \vert \b =
s^0$ and $r^1 \vert \b = s^1$, $\la s^0, s^1 \ra \in
I \vert \b$ and $j^*(s^0) = j^*(r^0) \vert j(\b) =
j^*(r^0 \vert \b)$, it must be the case that
$\la \rho, \sigma \ra \in \dom(j^*(s^0))$. This means
$\la \rho, \sigma \ra \in \dom(q^0_\b)$, a contradiction.
 
We define now an
$M[H_\k \ast H'_\k \ast H''_\k \ast H^4_\k]$-generic object
$H^{5,0}_\k$ over
$\C(j(\k^+), j(\l)) \ast \dot \K(j(\k), j(\l))$
so that $p \in g''G^1_{\k, \l}$ implies
$j^*(p)  \in H^{5,0}_\k$. First, define in
$M[H_\k \ast H'_\k \ast H''_\k]$ for $\b \in
(j(\k^+), j(\l))$ the partial ordering
$\C(j(\k^+), \b   ) \ast \dot \K(j(\k), \b) =
\{p \vert \b : p \in
 \C(j(\k^+), j(\l)) \ast \dot \K(j(\k), j(\l)) \}$,
ordered analogously to
$\C(j(\k^+), j(\l)) \ast \dot \K(j(\k), j(\l))$.
($p \vert \b$ essentially has the same meaning as when
$p \in
 \C(\k^+, \l) \ast \dot \K(\k, \l)$ and
$\b \in (\k^+, \l)$.) Next, note that since GCH holds in
$M[H_\k \ast H'_\k \ast H''_\k \ast H^4_\k]$ for
cardinals $\g \ge j(\k)$, $j(\k^{++}) \le
j(\l)$, and $j(\l)$ is regular, Lemma 5 implies
$M[H_\k \ast H'_\k \ast H''_\k \ast H^4_\k] \models
``\C(j(\k^+), j(\l)) \ast \dot \K(j(\k), j(\l))$ is
$j(\k^{++})$-c.c.,
$\C(j(\k^+), j(\l)) \ast \dot \K(j(\k), j(\l))$
has $j(\l)$ many antichains, and if $\A \subseteq
 \C(j(\k^+), j(\l)) \ast \dot \K(j(\k), j(\l))$
is a maximal antichain, then $\A \subseteq
 \C(j(\k^+), \b   ) \ast \dot \K(j(\k), \b   )$
for some $\b \in (j(\k^+), j(\l))$''. As
$V \models ``|j(\l)| = \l$'' and by Lemma 5 and the fact
$\k^{++} \le \l$,
$V[G_\k \ast G'_\k] \models
``\l$ is a cardinal'',
we can let $\la \A_\a : \a \in (\k^+, \l) \ra \in
V[G_\k \ast G'_\k]$ be an enumeration of the maximal
antichains of
$\C(j(\k^+), j(\l)) \ast \dot \K(j(\k), j(\l))$
present in
$M[H_\k \ast H'_\k \ast H''_\k \ast H^4_\k]$.
 
Working in $V[G_\k \ast G'_\k]$, we define now an
increasing sequence
$\la r_\a : \a \in (\k^+, \l) \ra$
so that
$\forall \a < \l[r_\a \ge q_\a$ and $r_\a \in
 \C(j(\k^+), j(\a)) \ast \dot \K(j(\k), j(\a))]$
and so that
$\forall \A \in \la \A_\a : \a \in (\k^+, \l) \ra
\exists \b \in (\k^+, \l) \exists r \in \A[r_\b \ge r]$.
Assuming we have such a sequence, $H^{5,0}_\k = \{p \in
 \C(j(\k^+), j(\l)) \ast \dot \K(j(\k), j(\l)) :
\exists r \in
 \la r_\a : \a \in (\k^+, \l) \ra[r \ge p]\}$ is our
desired generic object. To define
$\la r_\a : \a \in (\k^+, \l) \ra$, if $\a$ is a limit
and each $r_\b$ for $\b < \a$ is written as
$\la r^0_\b, r^1_\b \ra$, we let $r_\a = \la
\underset \b < \a \to{\cup} r^0_\b,
\underset \b < \a \to{\cup} r^1_\b \ra =
\underset \b < \a \to{\cup} r_\b$. By the facts
$\la q_\b : \b \in (\k^+, \l) \ra$
is (strictly) increasing and
$M[H_\k \ast H'_\k \ast H''_\k]$ is closed under
$< \l$ sequences with respect to $V[G_\k \ast G'_\k]$,
this definition is valid. Assuming now $r_\a$ has been
defined and we wish to define $r_{\a + 1}$, let
$\la \B_\b : \b < \eta < \l \ra$ be the subsequence of
$\la \A_\b : \b \le \a + 1 \ra$ containing each antichain
$\A$ so that $\A \subseteq
 \C(j(\k^+), j(\a + 1)) \ast \dot \K(j(\k), j(\a + 1))$.
Since $q_\a, r_\a \in
 \C(j(\k^+), j(\a)) \ast \dot \K(j(\k), j(\a))$,
$q_{\a + 1} \in
 \C(j(\k^+), j(\a + 1)) \ast \dot \K(j(\k), j(\a + 1))$,
and $j(\a) < j(\a + 1)$, if we write as before
$r_\a$, $q_\a$, and $q_{\a + 1}$ as
$r_\a = \la r^0_\a, r^1_\a \ra$,
$q_\a = \la q^0_\a, q^1_\a \ra$, and
$q_{\a + 1} = \la q^0_{\a + 1}, q^1_{\a + 1} \ra$,
then the condition $r'_{\a + 1} = \la r^0_\a \cup q^0_{\a + 1},
r^1_\a \cup q^1_{\a + 1} \ra =
r_\a \cup q_{\a + 1}$ is well-defined. This is because,
as our earlier observations show, any new elements of
$\dom(q^i_{\a + 1})$ won't be present in either
$\dom(q^i_\a)$ or $\dom(r^i_\a)$ for
$i \in \{0,1\}$. We can thus, using the fact
$M[H_\k \ast H'_\k \ast H''_\k]$ is closed under
$< \l$ sequences with respect to $V[G_\k \ast G'_\k]$,
define by induction an increasing sequence
$\la s_\b : \b < \eta \ra$ of elements of
$\C(j(\k^+), j(\a + 1)) \ast \dot \K(j(\k), j(\a + 1))$
so that
$s_0 \ge r'_{\a + 1}$, $s_\rho =
\underset \b < \rho \to{\cup} s_\b$ if $\rho$ is a limit,
and $s_{\b + 1} \ge s_\b$ is so that $s_{\b + 1}$
extends some element of $\B_\b$. The just mentioned closure
fact implies $r_{\a + 1} =
\underset \b < \eta \to{\cup} s_\b$ is a well-defined condition
in
$\C(j(\k^+), j(\a + 1)) \ast \dot \K(j(\k), j(\a + 1))$.
 
In order to show $H^{5,0}_\k$ is
$M[H_\k \ast H'_\k \ast H''_\k \ast H^4_\k]$-generic over
$\C(j(\k^+), j(\l)) \ast \break\dot \K(j(\k), j(\l))$,
we must show
$\forall \A \in \la \A_\a : \a \in (\k^+, \l) \ra
\exists \b \in (\k^+, \l) \exists r \in \A[r_\b \ge r]$.
To do this, we first note that $\la j(\a) : \a < \l \ra$
is unbounded in $j(\l)$. To see this, if $\b < j(\l)$ is
an ordinal, then for some $\g < \l$ and some
$f:\g \to M$ representing $\b$, we can assume that for
$\rho < \g$, $f(\rho) < \l$. Thus, by the regularity of
$\l$ in $V$, $\b_0 = \underset \rho < \g \to{\cup}
f(\rho) < \l$, and $j(\b_0) > \b$. This means by our
earlier remarks that if $\A \in \la
\A_\a : \a < \l \ra$, $\A = \A_\rho$, then we can let
$\b \in (\k^+, \l)$ be so that $\A \subseteq
 \C(j(\k^+), j(\b)) \ast \dot \K(j(\k), j(\b))$.
By construction, for $\eta > \max(\b, \rho)$, there is
some $r \in \A$ so that $r_\eta \ge r$. Finally,
since any $p \in
 \C(\k^+, \l) \ast \dot \K(\k, \l)$ is so that
for some $\a \in (\k^+, \l)$, $p = p \vert \a$,
$H^{5,0}_\k$ is so that if $p \in g''G^1_{\k, \l}$,
then $j^*(p) \in H^{5,0}_\k$.
 
Note now that our earlier work ensures $j^*$ extends to
$j^{**} :             V[G_\k \ast G^0_{\k, \l} \ast
G^2_{\k, \l}] \to
M                 [H_\k \ast H'_\k \ast H''_\k \ast H^4_\k]$.
The proof of Lemma 4 can be given in
$V[G_\k \ast G^0_{\k, \l} \ast G^2_{\k, \l}]$,
meaning we can assume without loss of generality that
$g = g_1 \circ g_0 \in
 V[G_\k \ast G^0_{\k, \l} \ast G^2_{\k, \l}]$.
Thus, by elementarity, $j^{**}(g^{-1}_1)$ is an
order isomorphism between
$\C(j(\k^+), j(\l)) \ast \dot \K(j(\k), j(\l))$
and $j^{**}(R) \ast j^{**}(\la Q, \le'_Q \ra)$.
Our earlier observations on the forcing equivalence between
$\la Q, \le_Q \ra$ and $\la Q, \le'_Q \ra$ (a generic for
one is a generic for the other)
combined with elementarity hence show
$H^{5,1}_\k = \{j^{**}(g^{-1}_1)(p) : p \in H^{5,0}_\k \}$
is an
$M[H_\k \ast H'_\k \ast H''_\k \ast H^4_\k]$-generic object over
$j^{**}(R) \ast j^{**}(\la \dot Q, \le_Q \ra)$.
Since the elementarity of $j^{**}$ implies
$j^{**}(g^{-1}_0)$ is an order isomorphism between
$j^{**}(R) \ast j^{**}(\la \dot Q, \le_Q \ra)$
and a dense subset of
$P^1_{j(\k), j(\l)}[S_{j(\l)}]$,
$H^{5}  _\k = \{j^{**}(g^{-1}_0)(p) : p \in H^{5,1}_\k \}$
is an
$M[H_\k \ast H'_\k \ast H''_\k \ast H^4_\k]$-generic object
over a dense subset of
$P^1_{j(\k), j(\l)}[S_{j(\l)}]$ so that
$p \in $ (a dense subset of) $P^1_{\k, \l}[S_\l]$ implies
 $j^{** }
 (p) \in H^5_{\k}$.
Therefore, for $H'''_\k = H^4_{\k} \ast H^5_{\k}$ and
$H = H_{\k} \ast H'_{\k} \ast H''_{\k} \ast H'''_{\k}$,
$j : V \to M$ extends to $k : V[G_\k \ast G'_\k] \to
M[H]$, so $V[G] \models ``\k$ is
$< \l$ supercompact''.
This proves Lemma 9.
\pbf
\hfill $\square $ Lemma 9
 
Lemmas 1-9 complete the proof of Theorem 1.
 
\pbf
\hfill $\square $ Theorem 1
 
We remark that as opposed to the statement of Theorem 1
in Section 0, the proof just given is not so that
$V$ and $V[G]$ share the same cardinals and cofinalities.
By Lemma 5, at each stage of the iteration, one cardinal
is collapsed. We will outline in Section 4 a notion of
forcing similar to the one found in [AS] that can be
used to give a proof of Theorems 1 and 2 (although not of
Theorem 3) preserving cardinals and cofinalities.
 
In conclusion to this section, we note that it is tempting to
think the above proof of Lemma 9 contains some hidden mistake,
i.e., that it can be extended to show $\k$ remains $\g$
supercompact for cardinals $\g \ge 2^\k$. If we tried to prove
Lemma 9 for some $\g \ge 2^\k$, then we would run into trouble
when we tried to construct the generic object $H^5_\k$. In
the above proof, the construction of
$H^5_\k$ depends heavily on Lemma 7, more specifically,
on the fact that $V[G_\k \ast H'_\k]$ and
$V[G_\k \ast G'_\k]$ have the same $< \l$ sequences of
elements of $V[G_\k \ast H'_\k]$.
Since $V[G_\k \ast G'_\k]$
contains a subset of $\l$ not present in
$V[G_\k \ast H'_\k]$, i.e., the generic object
$G^2_{\k, \l}$, if $\g \ge 2^\k$, it will be false that
$V[G_\k \ast H'_\k]$ and $V[G_\k \ast G'_\k]$
contain the same $\g$ sequences of elements of
$V[G_\k \ast H'_\k]$.
 
\S3 The Proofs of Theorems 2 and 3
 
We turn now to the proof of Theorem 2. Recall that we
are assuming the existence of a class function $\l$ so
that for any infinite cardinal $\d$,
$\l(\d) > \d^+$ is a regular cardinal
which is either inaccessible or is the successor of
a cardinal of cofinality $> \d$ with $\l(\d)$
below the least inaccessible $> \d$ if $\d$ is
singular and $\l(0) = 0$ and that our ground model $V$ is
such that $V \models ``$ZFC + GCH + $A$ is a proper class
of cardinals so that for each $\k \in A$, $h_\k : \k \to \k$
is a function and $j_\k : V \to M$ is an elementary embedding
witnessing the $< \l(\k)$ supercompactness
of $\k$ with $j_\k(h_\k)(\k) = \l(\k) < \k^*$ for $\k^*$
the least element of $A$ $> \k$''. Without loss of
generality (by ``cutting off'' $V$ if necessary), we
assume that for each $\k \in A$, $\sup(\{\d \in A : \d
< \k \}) = \rho_\k < \k$ and isn't inaccessible if
order type$(\{\d \in A : \d < \k \})$ is a limit ordinal.
 
For each $\k \in A$, let $P(\k, \l(\k))$ be the version of
the partial ordering $P$ used in the proof of Theorem 1
which ensures each inaccessible $\d \in (\l(\rho_\k), \k)$
isn't measurable yet $\k$ is $< \l(\k)$ supercompact, i.e.,
using the notation of Section 2, $P(\k, \l(\k))$ is the
$\k + 1$ stage Easton support iteration $\la \la
P_\a, \dot Q_\a \ra : \a \le \k \ra$ defined as follows:
\item{1.} $P_0$ is trivial.
\item{2.} Assuming $P_\a$ has been defined for $\a < \k$,
let $\d_\a$ be so that $\d_\a$ is the least cardinal
$\ge \underset \b < \a \to{\cup} \d_\a$ such that
$\force_{P_\a} ``\d_\a$ is inaccessible'', where
$\d_{-1} = \l(\rho_\k)$ and $\d_0$ is the least inaccessible
in $(\l(\rho_\k), \k)$. Then $P_{\a + 1} = P_\a \ast
\dot Q_\a$, with $\dot Q_\a$ a term for
$P^0_{\d_\a, h_\k(\d_\a)} \ast
P^1_{\d_\a, h_\k(\d_\a)}[\dot S_{h_\k(\d_\a)}]$.
\item{3.} $\dot Q_\k$ is a term for $P^0_{\k, \l(\k)} \ast
(P^1_{\k, \l(\k)}[\dot S_{\l(\k)}] \times
 P^2_{\k, \l(\k)}[\dot S_{\l(\k)}])$.
 
\noindent We then define the first partial ordering $P$ used
in the proof of Theorem 2 as $P = \{p \in
\underset \k \in A \to{\prod} P(\k, \l(\k))$ : support$(p)$
is a set$\}$, ordered by componentwise extension.
 
Note that for each $\k \in A$, we can write $P$ as
$T_\k \times T^\k$, where $T_\k =
\underset \{\d \in A : \d \le \k \} \to{\prod}
P(\d, \l(\d))$, $T^\k = \{p \in \underset
\{\d \in A : \d > \k \} \to{\prod} P(\d, \l(\d))$ :
support$(p)$ is a set$\}$, and $T_\k$ and $T^\k$ are both
ordered componentwise. For each $\k \in A$, the definition of
the component partial orderings of $T_\k$ and $T^\k$ ensures
that $T_\k$ is ${\l(\k)}^+$-c.c. and $T^\k$ is
${\l(\k)}^+$-strategically closed. This allows us to conclude in
the manner of [KiM] that $V^P \models$ ZFC.
 
\proclaim{Lemma 10}
$V^P \models ``2^\d = \l(\k)$ if $\k \in A$
and $\d \in [\k, \l(\k))$ + Every $\k \in A$ is $< \l(\k)$
supercompact + $\forall \k[\k$ is measurable iff $\k$ is
$< \l(\k)$ strongly compact iff $\k$ is $< \l(\k)$
supercompact]''.
\endproclaim
\demo{Proof of Lemma 10}
Using the notation above, for each $\k \in A$, write
$P = T_\k \times T^\k$; further, write $T_\k$ as
$T_{< \k} \times P(\k, \l(\k))$ for $T_{< \k} =
\underset \{\d \in A : \d < \k \} \to{\prod}
P(\d, \l(\d))$. Since each component partial ordering of
$T^\k$ is at least $< \sigma(\k)$-strategically closed for
$\sigma(\k)$ the least inaccessible $> \l(\k)$,
$T^\k$ is $< \sigma(\k)$-strategically closed, meaning by
the remark immediately following Lemma 5 and Lemma 9 that
$V^{T^\k \times P(\k, \l(\k))} \models
``2^\d = \l(\k)$ if $\d \in [\k, \l(\k))$, $\k$ is $< \l(\k)$
supercompact, and no cardinal $\d \in [\l(\rho_\k), \k)$ is
measurable''.
As our assumptions and the definition of $T_{< \k}$
ensure $V^{T^\k \times P(\k, \l(\k))} \models
``|T_{< \k}|$ $<$ the least inaccessible
$\theta(\k)$ above $\l(\rho_\k)$'', by the L\'evy-Solovay
arguments [LS], $V^{T^\k \times P(\k, \l(\k)) \times
T_{< \k}} = V^P \models ``2^\k = \l(\k)$, $\k$ is
$< \l(\k)$ supercompact, and no cardinal $\d \in
[\l(\rho_\k), \k)$ is measurable''.
Finally, as the definitions
of $\l(\k)$, $\rho_\k$, $\sigma(\k)$, and $\theta(\k)$
guarantee every inaccessible cardinal $\d$ in $V^P$ must
be so that $\d \in [\l(\rho_\k), \k]$ for some $\k \in A$
and every ordinal $\d$ $\ge$
the least inaccessible must be so that $\d \in [\theta(\k),
\sigma(\k)]$ for some $\k \in A$, $V^P \models ``\k$ is
measurable iff $\k$ is $< \l(\k)$ strongly compact iff
$\k$ is $< \l(\k)$ supercompact''.
This proves Lemma 10.
\pbf
\hfill $\square$ Lemma 10
 
Note it may be the case that $V^P \models ``$Some
cardinal $\k$ is $\l(\k)$ strongly compact''. In order
to ensure this doesn't happen, we must for each $\k \in A$
force with the
partial ordering $P^0_{\omega, \l(\k)}$ of Section 1, i.e.,
for each $\k \in A$, we add a non-reflecting stationary set
of ordinals of cofinality $\omega$ to $\l(\k)$. By a
theorem of [SRK], this guarantees $\k$ isn't $\l(\k)$
strongly compact.
 
Keeping the preceding paragraph in mind, we take as our ground
model $V^P = V_1$. Working in $V_1$, we let $R =
\{p \in \underset \k \in A \to{\prod} P^0_{\omega,
\l(\k)}$ : support$(p)$ is a set$\}$, ordered by componentwise
extension. As before, we can write for each $\k \in A$
$R = R_\k \times R^\k$, where $R_\k =
\underset \{\d \in A : \d \le \k \} \to{\prod}
P^0_{\omega, \l(\d)}$ and $R^\k = \{p \in
\underset \{\d \in A : \d > \k \} \to{\prod}
P^0_{\omega, \l(\d)}$ : support$(p)$ is a set$\}$. We can also
write $R_\k$ as $R_{< \k} \times P^0_{\omega, \l(\k)}$, where
$R_{< \k} = \underset \{\d \in A : \d < \k \} \to{\prod}
P^0_{\omega, \l(\d)}$. Again, for each $\k \in A$, the fact
that $V_1 \models ``2^\d = \l(\k)$ for each cardinal
$\d \in [\k, \l(\k))$'' and the definitions of $R_\k$ and
$R^\k$ ensure that $R_\k$ is ${\l(\k)}^+$-c.c. and
$R^\k$ is ${\l(\k)}^+$-strategically closed. Hence, once
more, $V^R_1 \models$ ZFC.
 
\proclaim{Lemma 11} $V_1$ and $V^R_1$ have the same cardinals
and cofinalities and $V^R_1 \models ``2^\d =
\l(\k)$ if $\k \in A$ and $\d \in [\k, \l(\k))$ + Every
$\k \in A$ is $< \l(\k)$ supercompact + $\forall\k[\k$ is
measurable iff $\k$ is $< \l(\k)$ strongly compact iff $\k$
is $< \l(\k)$ supercompact$]$ + No cardinal $\k$ is
$\l(\k)$ strongly compact''.
\endproclaim
\demo{Proof of Lemma 11}
We mimic the proof of Lemma 10. Using the notation of Lemma 10,
for $\k \in A$, since each component partial ordering
of $R^\k$ is at least $\sigma(\k)$-strategically closed and
$V^{R^\k}_1 \models ``P^0_{\omega, \l(\k)}$ is
$< \l(\k)$-strategically closed and is ${\l(\k)}^+$-c.c.'',
$V^{R^\k \times P^0_{\omega, \l(\k)}}_1 \models ``2^\d =
\l(\k)$ if $\d \in [\k, \l(\k))$, $\k$ is $< \l(\k)$
supercompact, no cardinal $\d \in [\l(\rho_\k), \k)$ is
measurable, cardinals and cofinalities for any $\d \le
\sigma(\k)$ are the same as in $V_1$, and $\k$ isn't
$\l(\k)$ strongly compact''. Since analogously to Lemma 10
$V^{R^\k \times P^0_{\omega, \l(\k)}}_1 \models
``|R_{< \k}| < \theta(\k)$'', as in Lemma 10,
$V^{R^\k \times P^0_{\omega, \l(\k)} \times R_{< \k}}_1
= V^R_1 \models ``2^\d = \l(\k)$ if $\d \in [\k, \l(\k))$,
$\k$ is $< \l(\k)$ supercompact, no cardinal $\d \in
[\l(\rho_\k), \k)$ is measurable, cardinals and cofinalities
for any $\d \in [\theta(\k), \sigma(\k)]$ are the same as in
$V_1$, and $\k$ isn't $\l(\k)$ strongly compact''. Once more,
every inaccessible $\d$ in $V^P$ must be so that
$\d \in [\l(\rho_\k), \k)$ for some $\k \in A$ and every
ordinal $\d$ $\ge$ the least inaccessible
must be so that $\d \in [\theta(\k), \sigma(\k)]$ for some
$\k \in A$, so $V^R_1$ and $V_1$ have the same cardinals and
cofinalities, and $V^R_1 \models ``2^\d = \l(\k)$ if
$\k \in A$ and $\d \in [\k, \l(\k))$ + Every $\k \in A$ is
$< \l(\k)$ supercompact + $\forall\k[\k$ is measurable iff
$\k$ is $< \l(\k)$ strongly compact iff $\k$ is $< \l(\k)$
supercompact$]$ + No cardinal $\k$ is $\l(\k)$ strongly
compact''. This proves Lemma 11.
\pbf
\hfill $\square$ Lemma 11
 
Lemmas 10 and 11 complete the proof of Theorem 2.
\pbf
\hfill $\square$ Theorem 2
 
We turn now to the proof of Theorem 3. We begin by giving
a proof of Menas' theorem that the least measurable limit
$\k$ of strongly compact or supercompact cardinals is not
$2^\k$ supercompact.
 
\proclaim{Lemma 12 (Menas [Me])}
If $\k$ is the least measurable limit of either
strongly compact or supercompact cardinals, then
$\k$ is strongly compact but isn't $2^\k$ supercompact.
\endproclaim
\demo{Proof of Lemma 12}
We assume without loss of generality that $\k$ is
the least measurable limit of strongly compact
cardinals. As readers will easily see, the proof given
works equally well if $\k$ is the least measurable limit
of supercompact cardinals.
 
Let $\la \k_\a : \a < \k \ra$ enumerate in increasing
order the strongly compact cardinals below $\k$.
Fix $\l > \k$ an arbitrary cardinal.
Let $\mu$ be any measure (normal or non-normal) over
$\k$, and let $\la \mu_\a : \a < \k \ra$ be a
sequence of fine, $\k_\a$-complete measures over
$P_{\k_\a}(\l)$. The set ${\cal U}_\l$ given by
$X \in {\cal U}_\l$ iff $X \subseteq P_\k(\l)$ and
$\{\a < \k : X \vert \k_\a \in \mu_\a \} \in \mu$,
where for $X \subseteq P_\k(\l)$, $\a < \k$,
$X \vert \k_\a = \{p \in X : p \in
P_{\k_\a}(\l) \}$. It can easily be verified that
${\cal U}_\l$ is a $\k$-additive, fine measure over
$P_\k(\l)$. Since $\l > \k$ is arbitrary, $\k$
is strongly compact.
 
Assume now that $\k$ is $2^\k$ supercompact, and let
$k : V \to M$ be an elementary embedding with critical
point $\k$ so that $M^{2^\k} \subseteq M$. By the fact
that $\k$ is the critical point of $k$, if $\d < \k$
is strongly compact, then $M \models ``k(\d) =
\d$ is strongly compact''. By the fact
$M^{2^\k} \subseteq M$, $M \models ``\k$ is
measurable''. Thus, $M \models ``\k$ is a measurable
limit of strongly compact cardinals'', contradicting
the fact that $M \models ``k(\k) > \k$ is the
least measurable limit of strongly compact
cardinals''. This proves Lemma 12.
\pbf
\hfill $\square$ Lemma 12
 
We return to the proof of Theorem 3. Recall that we are
assuming our ground model $V \models ``$ZFC + GCH + $\k$
is the least supercompact limit of supercompact cardinals +
$\l > \k^+$ is a regular cardinal which is either
inaccessible or is the successor of a cardinal of
cofinality $> \k$
and $h : \k \to \k$ is a function so that for some elementary
embedding $j : V \to M$ witnessing the $< \l$ supercompactness
of $\k$, $j(h)(\k) = \l $''. As in the proof of Theorem 1,
we assume without loss of generality that for every
inaccessible $\d < \k$, $h(\d) > \d^+$ and $h(\d)$ is
regular.
 
We also assume without loss of generality
that $h(\d) = 0$ if $\d$ isn't inaccessible and
that if $\d < \k$ is inaccessible and $\d'$
is the least supercompact cardinal $> \d$, then
$h(\d) < \d'$.
To see that the conditions on $h$ imply this last restriction,
assume that $h(\d) \ge \d'$ on the set of inaccessibles
below $\k$. It must then be the case that
$M \models ``$For some cardinal $\rho \le \l$,
$\rho$ is supercompact''. By the closure properties of
$M$, $M \models ``\k$ is $\zeta$ supercompact for all
$\zeta < \rho$''. It is a theorem of Magidor [Ma2] that
if $\a$ is $< \b$ supercompact and $\b$ is supercompact,
then $\a$ is supercompact. It is thus the case that
$M \models ``\k$ is supercompact''. Since $V \models
``\k$ is the least supercompact limit of supercompact
cardinals'' and $\k$ is the critical point of $j$, if
$V \models ``\a < \k$ is supercompact'', $M \models
``j(\a) = \a$ is supercompact''. Putting these last
two sentences together yields a contradiction to the
fact that $M \models ``j(\k) > \k$ is the least
supercompact limit of supercompact cardinals''.
 
We next show the following fact about Laver indestructibility
[L] which will play a critical role in the proof of
Theorem 3.
 
\proclaim{Lemma 13}
If $\d$ is a supercompact cardinal, then the
definition of the partial ordering which makes $\d$
Laver indestructible under $\d$-directed closed
forcings can be reworked so that for any fixed
$\g < \d$, there are no strongly compact cardinals in
the interval $(\g, \d)$.
\endproclaim
\demo{Proof of Lemma 13}
Let $f : \d \to V_\d$ be a Laver function, i.e.,
$f$ is so that for every $x \in V$ and every
$\sigma \ge {\hbox{\rm TC}}(x)$, there is a fine,
$\d$-complete, normal ultrafilter ${\cal U}_\sigma$
over $P_\d(\sigma)$ so that for $j_\sigma$ the
elementary embedding generated by ${\cal U}_\sigma$,
$(j_\sigma f)(\delta) = x$. The Laver partial ordering
$P^*$ which makes $\d$ Laver indestructible under
$\d$-directed closed forcings and destroys all strongly
compact cardinals in the interval $(\g, \d)$ is as
usual defined as a $\d$ stage Easton support iteration
$\la \la P^*_\a, \dot Q^*_\a \ra : \a < \d \ra$.
As in the original definition, at each stage $\a < \d$,
an ordinal $\rho_\a < \d$ is chosen, where at limit
stages $\a$, $\rho_\a =
\underset \b < \a \to{\cup} \rho_\b$. We define
$P^*_{\a + 1} = P^*_\a \ast \dot Q^*_\a$, where
$\dot Q^*_\a$ is a term for the trivial partial
ordering $\{\emptyset\}$ and $\rho_{\a + 1} =
\rho_\a$, unless for all $\b < \a$, $\rho_\b < \a$
and $f(\a) = \la \dot R, \sigma \ra$, where
$\sigma \neq \g$ is a regular cardinal $\ge \max(\g, \a)$ and
$\dot R$ is a term so that
$\force_{P^*_\a} ``\dot R$ is $\max(\g, \a)$-directed
closed''. Under these circumstances, $\dot Q^*_\a$
is a term for the partial ordering
$\dot R \ast \dot P^0_{\g, \sigma}$ and $\rho_\a
= \sigma^+$, where
$P^0_{\g, \sigma}$ has the same meaning as it did
in Section 1.
 
To show $P^*$ is as desired, let $\dot Q$ be a term
in the forcing language with respect to $P^*$ so that
$\force_{P^*}``\dot Q$ is $\d$-directed closed'', and
let $\eta \ge \d$. Let $\sigma > |{\hbox{\rm TC}}
(\dot Q)|$ be a regular cardinal so that
$\force_{P^* \ast \dot Q} ``\sigma > \zeta$ for
$\zeta = \max(|{\hbox{\rm TC}}(\dot Q)|,
2^{|{[\eta]}^{< \d}|})$''.
Take ${\cal U}_\sigma$ and the associated elementary
embedding $j_\sigma : V \to M$ so that
$(j_\sigma f)(\d) = \la \dot Q, 
\sigma \ra$, and call $j_\sigma$ $k$. By the definition
of $P^*$, in $M$ we have
$P^*_{\d + 1} = P^*_\d \ast \dot Q \ast
\dot P^0_{\g, \sigma} =
{(P^*)}^V \ast \dot Q \ast P^0_{\g, \sigma}$.
Thus, for $G \ast H$ $V$-generic over
$P^* \ast \dot Q$, since $P^*_\d$ is $\d$-c.c.,
the usual arguments (see Lemma 6.4 of [Ba]) show that
$M[G \ast H]$ is closed under $\zeta$ sequences
(in the sense of Lemma 9) with respect to
$V[G \ast H]$. Further, since $\force_{P^* \ast \dot Q}
``\sigma > \zeta$'',
by Lemma 8, the closure properties of $M[G \ast H]$
with respect to $V[G \ast H]$, and the definitions of
$P^*$ and $k(P^*)$ (including the fact both
$P^*$ and $k(P^*)$ are Easton support iterations of
partial orderings satisfying a certain degree of
strategic closure), the partial ordering $T \in
M[G \ast H]$ so that $P^*_\d \ast \dot Q \ast
\dot P^0_{\g, \sigma} \ast \dot T^* =
P^* \ast \dot Q \ast \dot T =
P^*_{k(\d)}$ (where $\dot T^*$ is a term for the
appropriate partial ordering) is $\zeta$-strategically
closed in both $M[G \ast H]$ and $V[G \ast H]$.
As in Lemma 9, this means that if $H'$ is
$V[G \ast H]$-generic over $T$, then
$M[G \ast H \ast H']$ is closed under $\zeta$ sequences
with respect to $V[G \ast H \ast H']$, and the embedding $k$
extends in $V[G \ast H \ast H']$ to
$k^* : V[G] \to M[G \ast H \ast H']$.
The definition of $\zeta$ and the fact
$k(Q)$ is $\zeta$-directed closed in both
$M[G \ast H \ast H']$ and $V[G \ast H \ast H']$
allow us to find in $V[G \ast H \ast H']$ a
master condition $q$ extending each $p \in
{k^*}''H$. If $H''$ is now a $V[G \ast H \ast H']$-generic
object over $k^*(Q)$ containing $q$, then
$k^*$ extends further in $V[G \ast H \ast H' \ast H'']$
to $k^{**} : V[G \ast H] \to
M[G \ast H \ast H' \ast H'']$. By the fact that
$T \ast k(\dot Q)$ is $\zeta$-strategically closed
in either $V[G \ast H]$ or $M[G \ast H]$ and the
definition of $\zeta$, the ultrafilter over
${(P_\d(\eta))}^{V[G \ast H]}$ definable via
$k^{**}$ in $V[G \ast H \ast H' \ast H'']$
is present in $V[G \ast H]$. Since $\eta$ was arbitrary,
$V[G \ast H] \models ``\d$ is supercompact''.
 
It remains to show that $V^{P^*} \models ``$No cardinal
in the interval $(\g, \d)$ is strongly compact''. To
see this, choose an embedding $k' : V \to M$ so that
$(k'f)(\d) = \la \dot Q, \d^+ \ra$, where
$\dot Q$ is a term with respect to $P^*$ for the
trivial partial ordering $\{\emptyset\}$. By the
definition of $P^*$, it will then be the case that
$M \models ``P^*_{\d + 1} = P^*_\d \ast \dot Q \ast
\dot P^0_{\g, \d} = {(P^*)}^V \ast \dot Q \ast
\dot P^0_{\g, \d}$ and the $\dot T$ so that
$P^*_{\d + 1} \ast \dot T = k'(P^*)$ is such that
$\force_{k'(P^*)} ``\d$ contains a non-reflecting
stationary set of ordinals of cofinality $\g$''.
By reflection, $\{\b < \d : P^*_{\b + 1} =
P^*_\b \ast \dot Q \ast \dot P^0_{\g, \b}$ and the
$\dot T$ so that $P^*_{\b + 1} \ast \dot T = P^*$
is such that
$\force_{P^*} ``\b$ contains a non-reflecting
stationary set of ordinals of cofinality
$\g$''$\}$ is unbounded in $\d$, where $\dot Q$
is a term with respect to $P^*_\b$ for the trivial
partial ordering $\{\emptyset\}$. By a theorem of
[SRK], if $\b$ contains a non-reflecting stationary
set of ordinals of cofinality $\g$, then there are
no strongly compact cardinals in the interval
$(\g, \b)$. Thus, the last two sentences immediately
imply $V^{P^*} \models ``$No cardinal in the interval
$(\g, \d)$ is strongly compact''. This proves Lemma 13.
\pbf
\hfill $\square$ Lemma 13
 
We note that in the proof of Lemma 13, GCH is not
assumed.  If GCH were assumed, then as in Lemma 9, we could
have taken the generic object $H'$ as an element of
$V[G \ast H]$.
 
We can now define in the manner of Theorem 1 the partial
ordering $P$
used in the proof of Theorem 3 by defining a $\k$ stage
Easton support iteration $P_\k = \la \la P_\a,
\dot Q_\a \ra : \a < \k \ra$ and then defining $P =
P_{\k + 1} = P_\k \ast \dot Q_\k$ for a certain partial
ordering $Q_\k$ definable in $V^{P_\k}$.
We first let $\la \d_\a : \a < \k \ra$ be the
continuous, increasing enumeration of the set
$\{\d < \k : \d$ is a supercompact cardinal$\} \cup
\{\d < \k : \d$ is a limit of supercompact
cardinals$\}$. We take as an inductive hypothesis that
the field of $P_\a$ is $\{\d_\b : \b < \a \}$ and that
if $\d_\a$ is supercompact, then $|P_\a| < \d_\a$.
The definition is then as follows:
\item{1.} $P_0$ is trivial.
\item{2.} Assuming $P_\a$ has been defined for $\a < \k$,
we consider the following three cases.
 
\itemitem{(1):} $\d_\a$ is a supercompact cardinal.
By the inductive hypothesis, since $|P_\a| < \d_\a$,
the L\'evy-Solovay results [LS] show that
$V^{P_\a} \models ``\d_\a$ is supercompact''.
It is thus possible in $V^{P_\a}$ to make $\d_\a$
Laver indestructible under $\d_\a$-directed closed
forcings. We therefore let $\dot Q_\a$ be a term for
the partial ordering of Lemma 13 making $\d_\a$ Laver
indestructible so that $\force_{P_\a} ``\dot Q_\a$ is
defined using partial orderings that are at least
${(2^{\g_\a})}^+$-directed closed and add non-reflecting
stationary sets of ordinals of cofinality
${(2^{\g_\a})}^+$ for ${\g_\a} =
\max(\sup(\{\d_\b : \b < \a \}),
h(\sup(\{\d_\b : \b < \a \})))$'',
and we define $P_{\a + 1} = P_\a \ast \dot Q_\a$.
(If $\a = 0$, then $\dot Q_\a$ is a term for the
Laver partial ordering of Lemma 13 where $\g = \omega$,
and $P_{\a + 1} = P_\a \ast \dot Q_\a$.)
Since $\dot Q_\a$ can be chosen so that
$|P_{\a + 1}| = \d_\a < \d_{\a + 1}$, the inductive
hypothesis is easily preserved.
 
\itemitem{(2):} $\d_\a$ is a regular limit of
supercompact cardinals.
Then $P_{\a + 1} = P_\a \ast \dot Q_\a$, with $\dot Q_\a$
a term for $P^0_{\d_\a, h(\d_\a)} \ast
P^1_{\d_\a, h(\d_\a)}[\dot S_{h(\d_\a)}]$, where
$\dot S_{h(\d_\a)}$ is a term for the non-reflecting
stationary subset of $h(\d_\a)$ introduced by
$P^0_{\d_\a, h(\d_\a)}$.
Since $\d_{\a + 1}$ must be supercompact, by the
conditions on $h$,
$|P_{\a + 1}| < \d_{\a + 1}$, so the inductive hypothesis
is once again preserved.
 
\itemitem{(3):} $\d_\a$ is a singular limit
of supercompact cardinals. Then $P_{\a + 1} =
P_\a \ast \dot Q_\a$, where $\dot Q_\a$ is a term
for the trivial partial ordering $\{\emptyset\}$.
 
\item{3.} $\dot Q_\k$ is a term for $P^0_{\k, \l} \ast
(P^1_{\k, \l}[\dot S_\l] \times P^2_{\k, \l}[\dot S_\l])$,
where again, $\dot S_\l$ is a term for the non-reflecting
stationary subset of $\l$ introduced by $P^0_{\k, \l}$.
 
Note that if $\a < \k$ is a limit ordinal, then since
$\k$ is the least supercompact limit of supercompact
cardinals, $\sup(\{\d_\b : \b < \a \}) = \d < \d'$,
where $\d'$ is the least supercompact cardinal in the
interval $[\d, \k)$. It is this fact that preserves the
inductive hypothesis at limit ordinals $\a < \k$
and at successor stages $\a + 1$ when $\d_\a$
is a singular limit of supercompact cardinals.
 
The intuitive motivation behind the above definition is much
the same as in Theorem 1. Specifically, below $\k$ at any
inaccessible limit $\d$ of supercompact cardinals, we must
force to ensure that $\d$ becomes non-measurable and is
so that $2^\d = h(\d)$. At $\k$, however, we must force so
as simultaneously to make $2^\k = \l$ while first
destroying and then resurrecting the $< \l$ supercompactness
of $\k$. The forcing will preserve the supercompactness
of every $V$-supercompact cardinal below $\k$
and will ensure there are no measurable limits of supercompacts
below $\k$. In addition, the forcing will guarantee that the
only strongly compact cardinals below $\k$ are those that
were supercompact in $V$.
Thus, $\k$ will have become the least measurable limit of
supercompact and strongly compact cardinals in the generic
extension.
 
\proclaim{Lemma 14} $V^P \models ``$If $\d < \k$ is
supercompact in $V$, then $\d$ is supercompact''.
\endproclaim
\demo{Proof of Lemma 14}
Let $\d < \k$ be a
$V$-supercompact cardinal. Write $P = R_\d \ast \dot R^\d$,
where $R_\d$ is the portion of $P$ whose field is all cardinals
$\le \d$ and $\dot R^\d$ is a term for the rest of $P$.
By case 1 in clause 2 of the inductive definition of $P$,
$V_1 = V^{R_\d} \models ``\d$ is supercompact
and is indestructible under $\d$-directed closed forcings''.
 
Assume now that $V^P = V^{R_\d}_1 \models ``\d$ isn't
supercompact'', and let $p = \la \dot p_\a : \a \le \k
\ra \in R^\d$ be so that over $V_1$, $p \force_{R_\d}
``\d$ isn't supercompact''.
By the remark after the proof of Lemma 3,
case 1 in clause 2 of the inductive definition of $P$,
and the fact each
$P^0_{\d_\a, h(\d_\a)}$ is $\d_\a$-directed closed for
$\a < \k$, it inductively follows that if $H$ is a
$V_1$-generic object over $R^\d$ so that $p \in H$, then
$H$ must be $V_1$-generic over a partial ordering $T^\d
\in V_1$ so that $p \in T^\d$ and so that $V_1 \models
``T^\d$ is $\d$-directed closed''. This means
$V_1[H] \models ``\d$ is supercompact''.
This
contradicts that over $V_1$, $p \force_{R_\d} ``\d$ isn't
supercompact''. This proves Lemma 14.
\pbf
\hfill $\square$ Lemma 14
 
\proclaim{Lemma 15} $V^P \models ``$No inaccessible $\d < \k$
which is a limit of $V$-supercompact cardinals is
measurable''.
\endproclaim
\demo{Proof of Lemma 15}
If $\d < \k$ is in $V^P$ an inaccessible limit of
$V$-supercompact cardinals, then since
$V^{P_\d} \subseteq V^P$, this fact must be true in
$V^{P_\d}$ as well. Hence, since $\d$ is so that
$P_\d$ is the direct limit of the system
$\la \la P_\a, \dot Q_\a \ra : \a < \d \ra$ and
$V \models {\hbox{\rm GCH}}$, $V^{P_\d} \models
``$All cardinals $\g \ge \d$ are the same as in
$V$ and GCH holds for all cardinals $\g \ge \d$''.
Therefore,
the same arguments as in Lemmas 3 and 6 show that
$V^{P_{\d + 1}} \models ``\d$ isn't measurable and
$2^\g = h(\d)$ if $\g \in [\d, h(\d))$ is a cardinal''.
(The same argument as in Lemma 8 also tells us
that
$V^P \models ``2^\g = \l$
if $\g \in [\k, \l)$ is a cardinal.)
It then follows by case 1 in clause 2 of the inductive
definition of $P$
that $V^P \models ``\d$ isn't measurable and
$2^\g = h(\d)$ if $\g \in [\d, h(\d))$ is a cardinal''.
This proves Lemma 15.
\pbf
\hfill $\square$ Lemma 15
 
\proclaim{Lemma 16}
$V^P \models ``$For any $\d < \k$,
$\d$ is supercompact in $V$ iff $\d$ is
supercompact iff $\d$ is strongly compact''.
\endproclaim
\demo{Proof of Lemma 16}
Let $\d < \k$ be strongly compact and not $V$-supercompact,
and let $\d' \in (\d, \k)$ be the least $V$-supercompact
cardinal $> \d$. Since Lemma 15 shows that no inaccessible
limit of $V$-supercompact cardinals is measurable,
$\sup(\{\b < \d : \b$ is a $V$-supercompact cardinal$\})
= \d_\a < \d$, where $\a < \k$ and $\d_\a$ is as in
the definition of $P$.
(If $\d < \d_0$, then let $\a = -1$ and $\d_\a = 0$.)
Thus, $\d \in (\d_\a, \d')$ and $\d' = \d_{\a + 1}$.
By the definition of $P$,
$P_{\a + 2} = P_{\a + 1} \ast \dot Q_{\a + 1}$,
where $\dot Q_{\a + 1}$ is so that
$\force_{P_{\a + 1}} ``\dot Q_{\a + 1}$ destroys all
strongly compact cardinals in the interval
$(\d_\a, \d_{\a + 1})$ by adding non-reflecting stationary
sets of ordinals of cofinality $  {(2^{\g_{\a + 1}})}^+$
to unboundedly many in $\d_{\a + 1}$ cardinals,
where as before, $\g_{\a + 1} = \max(\sup(\{\d_\b : \b <
\a + 1\}), h(\sup(\{\d_\b : \b < \a + 1\}))) =
\max(\d_\a, h(\d_\a)) $''. (If $\a = -1$,
then $\dot Q_{\a + 1}$ is so that
$\force_{P_{\a + 1}} ``\dot Q_{\a + 1}$ destroys all
strongly compact cardinals in the interval
$(\d_\a, \d_{\a + 1})$ by adding non-reflecting stationary
sets of ordinals of cofinality $\omega$
to unboundedly many in $\d_{\a + 1}$ cardinals''.)
Again by the definition of $P$, for the $\dot T$ so
that $P_{\a + 2} \ast \dot T = P$,
$\force_{P_{\a + 2}} ``\dot T$ is
${(2^{\g_{\a + 2}})}^+ =
{(2^{h(\d_{\a + 1})})}^+$-strategically closed'', so
$\force_P ``$There are unboundedly many in $\d_{\a + 1}$
cardinals in the interval $(\d_\a, \d_{\a + 1})$ containing
non-reflecting stationary sets of ordinals of either
cofinality ${(2^{\g_{\a + 1}})}^+$ or $\omega$''.
This means $V^P \models ``$No cardinal in the interval
$(\d_\a, \d_{\a + 1})$ is strongly compact'', a
contradiction. This, combined with Lemma 14, proves
Lemma 16.
\pbf
\hfill $\square$ Lemma 16
 
\proclaim{Lemma 17}
$V^P \models ``\k$ is the least measurable limit of
either strongly compact or supercompact cardinals''.
\endproclaim
\demo{Proof of Lemma 17}
By Lemma 16, if $\d < \k$ is strongly compact, then
$\d$ must be supercompact in both $V$ and $V^P$. By
Lemma 15, there are no measurable limits of $V$-supercompact
cardinals in $V^P$ below $\k$. This proves Lemma 17.
\pbf
\hfill $\square$ Lemma 17
 
Lemmas 13-17, together with the observation that the
same arguments as in Lemma 9 yield $V^P \models ``\k$ is
$< \l$ supercompact'', complete the proof of Theorem 3.
\pbf
\hfill $\square$ Theorem 3
 
\S4 Concluding Remarks
 
In conclusion to this paper, we outline an alternate
notion of forcing that can be used to construct
models witnessing Theorems 1 and 2 in which
cardinals and cofinalities are the same as in the
ground model. The forcing we use
is a slight variation of the forcing used in [AS].
Specifically, as in Section 1, we let $\d < \l$
be cardinals with $\d$ inaccessible, $\l > \d^+$
regular, and $\l$ either inaccessible or the successor
of a cardinal of cofinality $> \d$. We also assume as
in Section 1 that our ground model $V$ is so that GCH
holds in $V$ for all cardinals $\k \ge \d$, and we fix
$\g < \d$ a regular cardinal. As before, we define three
notions of forcing. $P^0_{\dell}$ is just the standard
notion of forcing for adding a non-reflecting stationary
set of ordinals of cofinality $\g$ to $\l$, i.e.,
$P^0_{\dell}$ is defined as in Section 1, only replacing
$\d$ in the definition with $\g$. If $\dot S$ is a term
for the non-reflecting stationary set of ordinals of
cofinality $\g$ introduced by $P^0_{\dell}$, then
$P^2_{\dell}[     S] \in V_1 =
V^{P^0_{\dell}}$ is the standard notion of forcing for
introducing a club set $C$ which is disjoint to $S$,
i.e., $P^2_{\dell}[     S]$ essentially has the same
definition as in Section 1.
 
To define $P^1_{\dell}[     S]$ in $V_1$, as in [AS]
or Section 1,
we first fix in $V_1$ a $\clubsuit(S)$ sequence
$X = \la x_\b : \b \in S \ra$. (Since each element of
$S$ has cofinality $\g$, either Lemma 1 of [AS] or our
Lemma 1 shows each $x \in X$ can be assumed to be
so that order type$(x) = \g$.) Then, in analogy
to the definition given in Section 1 of [AS],
$P^1_{\dell}[     S]$ is defined as the set of all
4-tuples $\la w, \a, \bar r, Z \ra$ satisfying the
following properties.
\item{1.} $w \in {[\l]}^{< \d}$.
\item{2.} $\a < \d$.
\item{3.} $ \bar r = \la r_i : i \in w \ra$ is a
sequence of functions from $\a$ to $\{0,1\}$, i.e.,
a sequence of subsets of $\a$.
\item{4.} $Z \subseteq \{x_\b : \b \in S\}$
is a set so that if $z \in Z$, then for some
$y \in {[w]}^\g$, $y \subseteq z$ and $z - y$
is bounded in the $\b$ so that $z = x_\b$.
 
\noindent As in [AS], the definition of $Z$ implies
$|Z| < \d$.
 
The ordering on $P^1_{\dell}[S]$ is given by
$\la w^1, \a^1, \bar r^1, Z^1 \ra \le
\la w^2, \a^2, \bar r^2, Z^2 \ra$ iff the following hold.
\item{1.} $w^1 \subseteq w^2$.
\item{2.} $\a^1 \le \a^2$.
\item{3.} If $i \in w^1$, then $r^1_i
\subseteq r^2_i$.
\item{4.} $Z^1 \subseteq Z^2$.
\item{5.} If $z \in Z^1 \cap {[w^1]}^\g$ and
$\a^1 \le \a < \a^2$, then $|\{i \in z :
r^2_i(\a) = 0\}| = |\{i \in z : r^2_i(\a) = 1\}| = \g$.
 
The intuition behind the definition of $P^1_{\dell}[S]$
just given is essentially the same as in [AS] or in
the remarks immediately following the definition of
$P^1_{\dell}[S]$ in Section 1 of this paper. Specifically,
we wish to be able simultaneously to make $2^\d = \l$,
destroy the measurability of $\d$, and be able to
resurrect the $< \l$ supercompactness of $\d$ if
necessary. $P^1_{\dell}[S]$ has been designed so as to
allow us to do all of these things.
 
The proof that $V^{P^1_{\dell}[S]}_1 \models
``\d$ is non-measurable'' is as in Lemma 3 of [AS].
In particular, the argument of Lemma 3 of [AS] will show
that $\d$ can't carry a $\g$-additive uniform ultrafilter.
We can then carry through the proof of Lemma 4 of [AS] to show
$P^0_{\dell} \ast (P^1_{\dell}[\dot S] \times
P^2_{\dell}[\dot S])$ is equivalent to
$\C(\l) \ast \dot C(\d, \l)$. The proofs of Lemma 5 of [AS]
and Lemma 6 of this paper will then show
$P^0_{\dell} \ast P^1_{\dell}[\dot S]$ preserves cardinals
and cofinalities, is $\l^+$-c.c., and is so that
$V^{P^0_{\dell} \ast P^1_{\dell}[\dot S]} \models
``2^\k = \l$ for every cardinal $\k \in [\d, \l)$''.
Then, if we assume $\k$ is $< \l$ supercompact with
$\l$ and $h : \k \to \k$ as in Theorem 1 and define in
$V$ an iteration $P$ as in Section 2 of this paper, we
can combine the arguments of Lemma 8 of [AS] and Lemma 8
of this paper to show $V$ and $V^P$ have the same cardinals
and cofinalities and $V^P \models ``$For all inaccessible
$\d < \k$ and all cardinals $\g \in [\d, h(\d))$,
$2^\g = h(\d)$, for all cardinals $\g \in [\k, \l)$,
$2^\g = \l$, and no cardinal $\d < \k$ is measurable''.
We can then prove as in Lemma 9 of this paper that
$V^P \models ``\k$ is $< \l$ supercompact''.
The proof now of Theorem 2 is as before, this time using the
just described iteration as the building blocks of the
forcing. This will allow us to conclude that the model
witnessing the conclusions of Theorem 2 thereby
constructed is so that it and the ground model contain the
same cardinals and cofinalities.
 
We finish by explaining our earlier remarks that it is
impossible to use the just described definitions of
$P^0_{\dell}$, $P^1_{\dell}[S]$,
and $P^2_{\dell}[S]$ to give a proof of
Theorem 3 of this paper. This is since if $V \models
``\k$ is a supercompact limit of supercompact cardinals,
$\mu < \k$ is supercompact, $\d$ and $\l$ are both
regular cardinals, and $\g < \mu < \d < \l$'', then
forcing with either $P^0_{\dell}$ or
$P^0_{\dell} \ast P^1_{\dell}[\dot S]$ will kill the
$\l$ strong compactness of $\mu$, as $S$ will be a
non-reflecting stationary set of ordinals of cofinality
$\g$ through $\l$ in either $V^{P^0_{\dell}}$ or
$V^{P^0_{\dell} \ast P^1_{\dell}[\dot S]}$.
(See [SRK] or [KiM] for further details.)
This type of forcing must of necessity occur if we use the
iteration described in the proof of Theorem 3.
 
\hfill\break\vfill\eject\frenchspacing
\centerline{References}\vskip .5in
 
\item{[AS]} A. Apter, S. Shelah, {\sl ``On the Strong Equality
between Supercompactness and Strong Compactness''}, to appear in
\underbar{Transactions AMS}.
\item{[Ba]} J. Baumgartner, {\sl ``Iterated Forcing"}, in: A. Mathias, ed.,
\underbar{Surveys in Set Theory},
 Cambridge University Press, Cambridge, England, 1--59.
\item{[Bu]} J. Burgess, {\sl ``Forcing"}, in: J. Barwise, ed.,
\underbar{Handbook of Mathematical Logic}, North-Holland, Amsterdam, 1977,
403--452.
\item{[C]} J. Cummings, {\sl ``A Model in which GCH Holds at
Successors but Fails at Limits''}, \underbar{Transactions
AMS} {\bf 329}, 1992, 1-39.
\item{[CW]} J. Cummings, H. Woodin, {\sl Generalised Prikry Forcings},
circulated manuscript of a forthcoming book.
\item{[J]} T. Jech, {\sl Set Theory}, Academic Press, New York, 1978.
\item{[Ka]} A. Kanamori, {\sl The Higher Infinite},
Springer-Verlag, New York and Berlin, 1994.
\item{[KaM]} A. Kanamori, M. Magidor, {\sl ``The Evolution of Large
Cardinal Axioms in Set Theory"}, in: \underbar{Lecture Notes in Mathematics}
 {\bf 669}, Springer-Verlag, Berlin, 1978, 99--275.
\item{[KiM]} Y. Kimchi, M. Magidor, {\sl ``The Independence between
the Concepts of Compactness and Supercompactness''}, circulated
manuscript.
\item{[L]} R. Laver, {\sl ``Making the Supercompactness of $\k$
Indestructible under $\k$-Directed Closed Forcing''},
\underbar{Israel J. Math.} {\bf 29}, 1978, 385--388.
\item{[LS]} A. L\'evy, R. Solovay, {\sl ``Measurable Cardinals
and the Continuum Hypothesis''}, \underbar{Israel J.}\hfil\break
\underbar{Math.} {\bf 5}, 1967, 234--248.
\item{[Ma1]}  M. Magidor, {\sl ``Changing Cofinality of
Cardinals''}, \underbar{Fundamenta Mathematicae} {\bf 99},\hfil\break
1978, 61--71.
\item{[Ma2]} M. Magidor, {\sl ``On the Role of Supercompact and Extendible
 Cardinals
 in Logic"},
\underbar{Israel J. Math.} {\bf 10}, 1971, 147--157.
\item{[Me]} T. Menas, {\sl ``On Strong Compactness and Supercompactness"},
 \underbar{Annals Math. Logic} {\bf 7}, 1975, 327--359.
\item{[MS]} A. Mekler, S. Shelah, {\sl ``Does $\k$-Free Imply Strongly
$\k$-Free?"}, in:
\item{} \underbar{Proceedings of the Third Conference on Abelian Group Theory},
 Gordon and
Breach, Salzburg, 1987, 137--148.
\item{[SRK]} R. Solovay, W. Reinhardt, A. Kanamori, {\sl ``Strong Axioms of
 Infinity
 and Elementary Embeddings"},  \underbar{Annals Math. Logic} {\bf 13}, 1978,
 73--116.
 
\bye